\documentclass[a4paper,12pt,reqno]{amsart}
\usepackage{amssymb}
\usepackage[dvips]{graphicx}
\usepackage{hyperref}

\numberwithin{equation}{section}
\input lipschitz2.def

\title{Local Properties of $J$-complex Curves\\[3pt] 
in Lipschitz-Continuous Structures}
\author{S. Ivashkovich, V. Shevchishin}

\address{
Universit\'e de Lille-1, UFR de Math\'ematiques, 59655 Villeneuve
d'Ascq, France.} \email{ivachkov@math.univ-lille1.fr}
\address{IAPMM Nat. Acad. Sci. Ukraine
Lviv, Naukova 3b,
79601 Ukraine.}
\address{
Department of Mathematics, University of Hamburg, 
Bundesstrasse 55, D-20146 Hamburg, Germany.}
\email{shevchishin@googlemail.com }

\date{\today}

\begin{document}

\begin{abstract}
We prove the existence of primitive curves and positivity of
intersections of $J$-complex curves for Lipschitz-continuous almost
complex structures. These results are deduced from the
Comparison Theorem for $J$-holomorphic maps in Lipschitz structures,
previously known for $J$ of class $\calc^{1, Lip}$. We also give the
optimal regularity of curves in Lipschitz structures. It occurs to
be $\calc^{1,LnLip}$, \ie the first derivatives of a $J$-complex curve
for Lipschitz $J$ are Log-Lipschitz-continuous. A simple example
that nothing better can be achieved is given. Further we prove the
Genus Formula for $J$-complex curves and determine their principal
Puiseux exponents (all this for Lipschitz-continuous $J$-s).
\end{abstract}

\maketitle

\setcounter{tocdepth}{1}

\tableofcontents

\newsect[sectINT]{Introduction}

An almost complex structure on a real manifold $X$ is a section of $\endr{TX}$
such that $J^2=-\id$. In this paper we are interested in the case when $J$ is
Lipschitz-continuous. A $J$-holomorphic curve in an almost complex manifold
$(X,J)$ is a $\calc^1$-map $u:S\to X$ from a complex curve $(S,j)$ to $X$ such
that $du$ commutes with complex structures, \ie for every $s\in S$ one has the
equality
\[
du(s)\circ j(s) = J(u(s))\circ du(s)
\]
of mappings $T_sS\to T_{u(s)}X$. In a local $j$-holomorphic
coordinate $z=x+iy$ on $S$ and local coordinates
$u=(u_1,...,u_{2n})$ on $X$ this writes as

\begin{equation}
\frac{\d u}{\d x} + J(u)\frac{\d u}{\d y} = 0, \eqqno(1.1)
\end{equation}
\ie as the Cauchy-Riemann equation.

\medskip The goal of this paper is to prove that $J$-complex curves for
Lipschitz-continuous $J$ possess all nice properties of the usual complex
curves. 

\newprg[prgINT.prim]{Existence of primitive parameterizations}
Recall (see also Definition \ref{prim-map}) that a $J$-holomorphic
map $u:S\to X$  is called {\slsf primitive} if there are no disjoint 
non-empty open sets $U_1, U_2$ in $S$ that $u(U_1)=u(U_2)$. Our first 
result states that every non-primitive $J$-holomorphic map factorizes 
through a primitive one, provided $J$ is Lipschitz-continuous.

\bigskip \smallskip\noindent%
{\bf Theorem A.} {\it Let $(S, j)$ be a smooth connected complex curve and $u:
 (S,j) \to (X,J)$ a non-constant $J$-holomorphic map with $J$ being
 Lipschitz-continuous. Then there exists a smooth {connected} complex curve
 $(\tilde S, \ti \jmath)$, a \emph{primitive} $J$-holomorphic map $\ti u: (\wt
 S,\ti \jmath) \to (X,J)$ and a surjective holomorphic map $\pi : (S,j) \to (\wt S,\ti
 \jmath)$ such that $u = \ti u\circ \pi$. }

\begin{exmp} \rm We would like to underline here that $\pi$ in general is not a
 \emph{covering}.  Let us give a simple, but instructive example. As a
 parameterizing complex curve $S$ consider the interior of the ellipse $\{
 \frac{1}{4}\cos\phi + i \sin\phi : \phi \in (-\pi,\pi]\}$.  The structure $j$ on $S$ is
 standard, \ie $j=J\st$ is the multiplication by $i$. The almost complex
 manifold in this example is $(\cc,J\st)$. The $J$-holomorphic map $u:S\to \cc$
 (\ie the usual holomorphic function) is taken as follows:
\[
u(z) = \Big(\frac{1}{2}z + 1\Big)^7.
\]

\begin{figure}[h]
\centering
\includegraphics[width=2.0in]{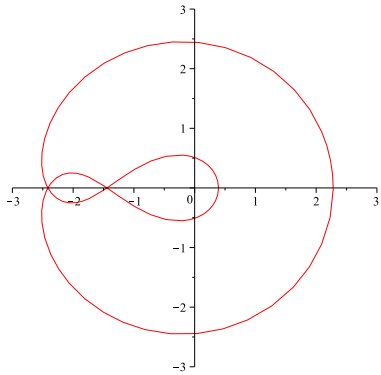}
\caption{On this picture we see the image $u(S)$. It overlaps around the point $-2$.} \label{examp-fig}
\end{figure}

Since we have an overlapping in the image, this map is not primitive. The
Theorem A for this example states that if one takes as $\tilde S$ the image
$u(S)$, as $\tilde u =\id$ and as the projection $\pi = u$ then $\tilde u
:\tilde S\to \cc$ is primitive and $u=\tilde u \circ \pi$.
\end{exmp}

\begin{rema} \rm 
\label{covering}
Let us notice that in the case when $S$ is closed the curve $\tilde S$ is also
closed and $\pi :S\to \tilde S$ is a ramified covering, see Corollary
\ref{closed}.
\end{rema}

\newprg[prgINT.posit]{Positivity of intersections} We denote by $\Delta_r$ the disc
of radius $r>0$ in $\cc$, $\Delta$ stands for $\Delta_1$. Further, let $u_i:\Delta\to(\cc^2,
J)$, $i=1, 2$ be two distinct (see Definition \ref{dist-map}) primitive
$J$-complex discs such that $u_1(0)=u_2(0)$.  Set $M_i\deff u_i (\Delta)$. Remark
that for Lipschitz-continuous $J$ the notion of multiplicity of zero of a
$J$-holomorphic map is well defined. Namely, due to the Corollary 3.1.3 from
\cite{IS1} every non-constant $J$-holomorphic map $u:(\Delta ,0)\to (\cc^n,J)$,
$J(0)=J\st$, has the form
\begin{equation}
\eqqno(repres1)
u(z) = v_0z^{\mu} + O(|z|^{\mu +\alpha}),
\end{equation}
where $0\not= v_0\in \cc^n$ is called the tangent vector to $u(\Delta)$ at the origin,
$\mu \geq 1$ is a natural number, called the multiplicity of $u$ at $0$, and $0<\alpha<1$. 
Our second result is the following

\smallskip\noindent{\bf Theorem B.}
{\it  Let  $J$ be a Lipschitz-continuous almost complex structure in $\cc^2$ and let
$u_1,u_2:\Delta\to \cc^2$ be two distinct $J$-holomorphic mappings. Then the following
holds:

\smallskip
\sli For every $0<r<1$ the set $\{\, (z_1,z_2) \in \Delta^2_r : u_1(z_1)= u_2(z_2)\,\}$ is
finite.

\smallskip
\slii If $\mu_1$ and $\mu_2$ are the~multiplicities of
$u_1$ and $u_2$ at $z_1$ and $z_2$, respectively, with
$u_1(z_1)=u_2(z_2)=p$, then the intersection index $\delta_p$  of
branches of $M_1$ and $M_2$ at $z_1$ and $z_2$ is at least
$\mu_1\cdot \mu_2$. In particular, $\delta_p$ is always strictly positive.

\smallskip \sliii $\delta_p=1$ if and only if $M_1$ and $M_2$ intersect at $p$
transversally. }

\newprg[prgINT.compar]{Comparison Theorem}
Both results of Theorems A and B are obtained using the following
statement, which should be considered as the main result of this paper.
Let $J$ be a Lipschitz-continuous almost complex structure in the unit ball
$B$ in $\cc^n$ and let $u_1, u_2: \Delta \to B$ be two $J$-holomorphic maps such that
$u_1(0) = u_2(0) =0$. Assume that both maps have the same order and the same
tangent vector at $0$, \ie in the representation \eqqref(repres1) one has
\begin{equation}
\eqqno(repres2)
u_i(z) = v_0 z^{\mu} + O(|z|^{\mu+\alpha}) \quad\text { for } \quad i=1,2.
\end{equation}
Our goal is to compare these mappings, \ie to describe in a best possible
way their difference.

\medskip\noindent{\bf Comparison Theorem.} {\it Let $(X,J)$ be an almost complex 
manifold with Lipschitz-continuous almost complex structure $J$ and let
$u_i:\Delta\to X$ be $J$-holomorphic mappings having the same order and the
same tangent vector at $0$ as in \eqqref(repres2).

\smallskip\noindent {{\bf(a)}\ \ 
There exists a holomorphic function $\psi$ of the form $\psi (z) = z + O(z^2)$, an
integer $\nu > \mu$ and a $\cc^n$-valued function $w$,  which belongs to 
$L^{1,p}_{loc}$ for all $2<p<\infty$, such that for some $r>0$
\begin{equation}
u_2(z) = u_1(\psi(z)) + z^{\nu} w(z) \text{ for $z\in\Delta (r)$}.
\eqqno(compar)
\end{equation}
Moreover, the following alternative holds:

\sli either $w(z)$ vanishes identically and then $u_2(\Delta
(\eps))\subset u_1(\Delta)$ for some $\eps >0$,

\slii or,  the vector $w(0)$ can be chosen orthogonal to $v_0$, in particular,
$w(0)\not= 0$ and
\begin{equation}
\eqqno(log-est)
|\pr_{v_0}w(z)| \leq C\cdot|z|\,\ln\frac{1}{|z|}\cdot|w(z)|.
\end{equation}
\smallskip\noindent
{\bf(b)}\ \ Let  $1\not=d\leq\mu$ be a divisor of $\mu$, 
 and $\eta=e^{2\pi\isl/d}$ be the primitive root of unity of degree $d$.  Let
 $u_1(\eta z)=u_1(\psi(z))+z^{\nu} w(z)$ be the presentation provided by \eqqref(compar)
for the map $u_2(z):=u_1(\eta z)$.  Then there exists 
a holomorphic reparameterization $\phi $ of the form
 $\phi(z)=z+O(z^2)$ such that 
\begin{itemize}
\item[\sli]  $u_1(\phi(\eta z))\equiv u_1(\phi(z))$ in the case when $w(z)\equiv0$;

\item[\slii] $u_1(\phi(\eta z)) = u_1(\phi(z)) + w(0) z^\nu+O(|z|^{\nu + \alpha})$
 otherwise. Moreover, in this case $\nu$ is not a multiple of $d$.
\end{itemize}
} }

\smallskip\noindent
In \eqqref(log-est) $\pr_{v_0}w$ denotes the orthogonal
projection of vector $w$ onto the vector $v_0$. Note that
$|z|\,\ln\frac{1}{|z|}=o(|z|^{\alpha})$ for any $0<\alpha <1$. In fact, from our proof
it follows that the vector $w(0)$ can be taken to belong to a prescribed
$(n-1)$-dimensional complex subspace $E_2$ of $\cc^n$ transverse to $v_0$, see
Remark \ref{rem-transvers}. Of course, the choice of $E_2$ will affect the
reparameterization function $\psi$ and the vector-function $w(z)$. This theorem
is proved in Section \ref{sectLOC}.

\newprg[prgINT.regul]{Optimal regularity of complex curves in Lipschitz structures}
Our next result is about the
precise regularity of $J$-complex curves for Lipschitz-continuous
$J$. Recall that a mapping $f$ from a compact set $B\subset \rr^n$
to a normed space is called Log-Lipschitz-continuous if
\begin{equation}
\norm{f}_{\calc^{LnLip}(B)}\deff \norm{f}_{L^{\infty}(B)} + \sup
\left\lbrace \frac{|f(x)-f(y)|}{|x-y|\cdot
\ln{\frac{1}{|x-y|}}}: x\not= y \in B, |x-y|\leq\frac12 \right\rbrace  < \infty ,
\eqqno(log-lip)
\end{equation}
and in this case $\norm{f}_{\calc^{LnLip}(B)}$ is called its Log-Lipschitz norm. 
Usually one takes $B$ to be the closure of a relatively
compact domain $D$ and then one sets
$\norm{f}_{\calc^{LnLip}(D)}=\norm{f}_{\calc^{LnLip}(\bar D)}$.
Without the logarithm in the right hand side \eqqref(log-lip) gives the
Lipschitz norm of $f$, which is denoted by
$\norm{f}_{\calc^{Lip}(D)}$.

\bigskip\noindent{\bf Theorem C.} {\it
Let $u: \Delta \to (\rr^{2n},J)$ be a $J$-holomorphic map. If
$J\in \calc^{Lip}(\rr^n)$ then $u\in \calc^{1,LnLip}$ \ie
the differential of $u$ is Log-Lipschitz-continuous. }

\medskip We show by a simple example that nothing better can be achieved, in
particular $u$ need not belong to $\calc^{1,Lip}$.

\newprg[prgINT.genus]{Local and global numerical invariants of complex curves} 
We also prove the following useful formula relating the local and
global invariants of a $J$-complex curve, known as Genus or Adjunction
Formula.  Let $M =\bigcup_{j=1}^d M_j$ be a compact $J$-complex curve in an almost
complex surface $(X, J)$ with the~distinct irreducible components $\{M_j\}$,
where $J$ is Lipschitz-continuous. Denote by $g_j$ the genera of parameter
curves $S_j$, \ie each $M_j$ is the image $u_j(S_j)$ of a compact Riemann
surface $S_j$ of genus $g_j$ under a primitive $J$-holomorphic mapping
$u_j:S_j\to X$.  Denote by $[M]^2$ the homological self-intersection of $M$ and
by $c_1(X, J)[M]$ the value of the first Chern class of $(X,J)$ on $M$. These
are the global invariants of $M$. Denote by $\delta$ the sum of all local
intersection indices $\delta_p$ of points $p\in M$.  For any singular local branch of
$M$ through a point $p$ we define the cusp index $\varkappa_p$ as the virtual number
of ordinary double points (see Definition \ref{cusp-ind}) and denote by $\varkappa$
the sum of the cusp-indices of all cusps of $M$. These are the local
invariants of $M$. These invariants are related by the following

\bigskip\noindent{\bf Theorem D.} {\sl (Genus Formula)} {\it If $J$ is
Lipschitz-continuous and $M=\bigcup_{j=1}^dM_j$ is a compact $J$-complex curve, where all
irreducible components $M_j$ of $M$ are distinct, then}
\begin{equation}
\eqqno(genus-form)
\sum_{j=1}^d g_j = \frac{[M]^2 - c_1(X, J)[M]}{2} + d -\delta - \varkappa .
\end{equation}
The novelty here is, of course, in the ability to define the
local invariants and to prove that they possess some nice properties (like
positivity) under the assumption of Lipschitz continuity of $J$ only. The
local intersection indices $\delta_p$ are explained by Theorem B.
The formula (\ref{classic}) below computes the cusp indices $\varkappa_p$.

\newprg[prgINT.puiseux]{Puiseux characteristics of  $J$-complex curves}
In the last part of this paper we provide an analog of the Puiseux series
for a $J$-complex curve in a Lipschitz-continuous structure $J$.

\bigskip\noindent{\bf Theorem E.} {\it Let $J$ be a Lipschitz-continuous 
almost complex structure in the unit ball $B\subset\cc^n$ with $J(0)=J\st$, 
and let $u:\Delta\to B$ be a \emph{primitive} $J$-holomorphic map having 
the form $u(z)=v_0z^\mu+O(|z|^{\mu+\alpha})$, where $\mu\geq2$. Then there 
exist a uniquely defined sequence of natural numbers $p_0=\mu<p_1<\cdots p_l$, 
a sequence of vectors  $v_1,\ldots,v_l$ each orthogonal to $v_0$, 
$J$-holomorphic maps $u_i:\Delta_r\to B$,  $i=0,\ldots,l$, and a complex 
polynomial $\phi(z)=z+O(z^2)$ with the following
 properties:
\begin{itemize}
\item The sequence  $d_i:=\gcd(p_0,\ldots,p_i)$ is
 strictly decreasing, $d_0>d_1>\cdots>d_l$ and $d_l=1$;

\smallskip\item Each map $u_i:\Delta_r\to B$ is primitive;
\item $u_0(z)=v_0z+O(|z|^{1+\alpha})$, $u_i(z) = u_{i-1}(z^{d_{i-1} /d_i})+ v_i\cdot z^{p_i /d_i} +
 O(|z|^{p_i /d_i + \alpha})$ for $i=1,\ldots,l$;

\smallskip\item $u(\phi(z)) - u_i(z^{d_i}) = v_{i+1}z^{p_{i+1}} + O(|z|^{p_{i+1}+\alpha})$ for
 $i=0,\ldots,l-1$ and $u_l(z)=u(\phi(z))$; In particular, if $\eta_i:=e^{2\pi\isl/d_i}$ is the
 primitive root of unity, then
\[
u(\phi(\eta_iz)) - u(\phi(z)) = (\eta_i^{p_{i+1}}-1)
 v_{i+1}z^{p_{i+1}} + O(|z|^{p_{i+1}+\alpha}).
\]
\end{itemize}
}

\smallskip\noindent  We call the sequence of the maps $u_i(z)$ a {\slsf Puiseux
approximation} of the map $u(z)$, the degrees
$p_0=\mu<p_1<\cdots<p_l$ the {\slsf characteristic exponents}, and the numbers
$d_i=\gcd(p_0,\ldots,p_i)$ the {\slsf associated divisors}.
The whole sequence $(p_0,\ldots,p_l)$ is called the {\slsf singularity type} of the
map $u:\Delta\to B$ at $0$ or of the pseudoholomorphic curve $u(\Delta)=M$ at $0$. 
The exponent $p_0$ is called the {\slsf multiplicity} or {\slsf order} of $u$
or of the curve $M$. 

\rm In the classical literature (\cite{Co}, \cite{BK}) the characteristic exponents are
also called {\slsf essential exponents} or even {\slsf Puiseux
characteristics} (\cite{Wl}); the difference $p_0-p_1$ is called the {\slsf
class of the singularity}, see e.g.\ \cite{Co}.

\smallskip Let us illustrate the notions involved in the Theorem E by an example.

\begin{exmp} \rm
\label{ex2}
Consider a (usual) holomorphic map $u:\Delta\to\cc^2$ given by
\[
u(z)=(z^{6}, z^{8}+z^{11}).
\]
Then the $v_0=(1,0)$ is the tangent vector at $z=0$ and $\mu=p_0=6$ is the
multiplicity. Further, its characteristic exponents - where the common divisor drops - are
$\mu=p_0=6,p_1=8,p_2=11$. The corresponding divisors are
$d_0=p_0=6$, $d_1=2$, $d_3=1$. Further, $v_1=e_2$ and $v_2=e_2$. A Puiseux
approximation sequence for $u(z)$ is:\\[4pt]
$\bullet\ $ $u_0(z)=(z,0)$,
\\[4pt]
$\bullet\ $ $u_1(z)=(z^3,z^4)$,
\\[4pt]
$\bullet\ $ $u_2(z)=u(z)$.
\end{exmp}

\smallskip Finally, we prove that the
following classical formula for the index of a cusp of a planar curve
\begin{equation}
\label{classic}
\varkappa_p = \frac{1}{2}\sum_{j=1}^l(d_{j-1} - d_j)(p_j - 1), \quad \text{where} \quad
d_j\deff \gcd(p_0,...,p_j),
\end{equation}
remains valid for $J$-complex curves in Lipschitz-continuous $J$.

\newprg[prgINT.notes]{Notes} {\sl 1.} In the classical case, \ie for algebraic
curves the Genus Formula is due to Clebsch and Gordan in the case when the
curve in question has only nodal singularities, \ie transverse intersections,
see \cite{CG} and Historical Sketch in \cite{S}. For curves with cusps the
Genus Formula is due to Max Noether, see p. 180 in \cite{Fi}.

\smallskip\noindent{\sl 2.} The statements of Theorems A and B and the Genus
Formula where proved in \cite{MW} for $J\in\calc^2$. In \cite{Sh} the description of a
singularity type of a $J$-complex curve for $J\in\calc^2$ was given.  For $J$-s of class
$\calc^{1,Lip}$ the positivity of intersections and the part (a) of the
Comparison Theorem where proved in \cite{Sk2}.

\smallskip\noindent{\sl 3.} Our interest to Lipschitz-continuous structures
comes from the following facts. First, a blowing-up of a general almost
complex manifold $(X,J)$, with $J\in \calc^{\infty}$, results to an almost complex
manifold $(\tilde X,\tilde J)$ with only Lipschitz-continuous $\tilde J$. Such
a blow-up should be performed in a special coordinate system, adapted to $J$,
see \cite{Du}. It is not difficult to see that the ordinary double points and
simple cusps can be resolved by this procedure, as in the classical case, and
give a smooth curve. Now the results of the present paper make possible to
work with such curves as with usual complex ones.

\smallskip\noindent{\sl 4.} Second, the condition of
Lipschitz-continuity cannot be relaxed in any of the statements above. We give
an example of two different $J$-complex curves which coincide by a non-empty
open subset for $J$ in all H\"older classes, or $J$ in all $L^{1,p}$ for all $p<\infty$.
In particular, the unique continuation statement of Proposition 3.1 from \cite{FHS}
fails to be true. In fact in our example $J$ is ``almost'' Log-Lipschitz, \ie
is essentially better than $\bigcap_{p<\infty} L^{1,p}$.

\smallskip\noindent{\sl 5.} At the same time let us point out that even in
continuous almost complex structures pseudoholomorphic curves have certain
nice properties: every two sufficiently close points can be joined by a
$J$-complex curve, a Fatou-type boundary values theorem is still valid, see
\cite{IR2}; Gromov compactness theorem both for compact curves and for curves
with boundaries on immersed totally real (\eg, Lagrangian) submanifolds
hold true for continuous $J$-s, see \cite{IS3,IS4}.

\smallskip\noindent{\sl 6.} To our knowledge the first result about $J$-complex
curves in Lipschitz structures appeared in \cite{NW}, where the existence 
of $J$-complex curves through a given point in a given direction was proved 
for $J\in\calc^{\alpha}$. Further progress is due to J.-C. Sikorav in \cite{Sk1}, 
see more about that in Remark 3.1 after the proof of Lemma \ref{lem3.2}.

\smallskip\noindent{\sl Acknowledgments.} This research  was partially done
during the first authors stay in the Max-Planck-Institute f\"ur Mathematik,
Bonn. He would like to thank this Institution for the hospitality. 

\smallskip We would like to give our thanks to the Referee of this paper
for numerous valuable remarks and suggestions.


\newsect[sectZD]{Zeroes of the Differential of a $J$-Holomorphic Map}

\newprg[prgZD.morrey]{Inner regularity of pseudoholomorphic maps}

Let us first recall few standard facts. For $0< \alpha \leq1$
consider the H{\"o}lder space $\calc^{k,\alpha }(\Delta , \cc^n)$ of
mappings $u:\Delta\to\cc^n$ equipped with the norm
\[
\norm{ u}_{\calc^{k,\alpha }(\Delta)}:= \norm{ u}_{\calc^k(\Delta)}
+ \sup_{z\not= w, |i|=k }\frac{\norm{ D^iu(z)-D^iu(w)}}{|
z-w|^\alpha}<\infty.
\]
For $k=0$ and $\alpha=1$ the space $\calc^{0,1}(\Delta, \cc^n)$ is
the Lipschitz space and  is denoted
by $\calc^{Lip}(\Delta, \cc^n)$. The Lipschitz constant of a map
$u\in \calc^{Lip}(\Delta, \cc^n)$ is defined as
\[
Lip_{\Delta}(u) \deff \sup_{z\not= w\in \Delta}\frac{\norm{u(z) - u(w)}}{|z - w|}.
\]
We also consider Lipschitz continuous (operator valued) functions on
relatively compact subsets of $\rr^{2n}$  with an obvious definitions and
notations for them.  Another scale of functional spaces, which will be used
in this paper, are the Sobolev spaces $L^{k,p}(\Delta,\cc^n)$,  $k\in\nn , 1\leq p\leq +\infty$,
with the norm
\[
\norm{ u}_{L^{k,p}(\Delta)}:= \sum_{0\leq | i| \leq k}\norm{ D^iu}_{L^p(\Delta)},
\]
where $i=(i_1,i_2)$, with $i_1,i_2 \geq 0$, $|i|= i_1 + i_2$, and
$D^iu:=\frac{\d^{|i|} u}{\d x^{i_1} \d y^{i_2}}$.  Let us also notice the
equality $L^{k,\infty}(\Delta,\cc^n)=\calc^{k-1,1}(\Delta,\cc^n)$ and the continuous Sobolev
imbeddings $L^{k,p}(\Delta,\cc^n)\hookrightarrow\calc^{k-1,\alpha}(\Delta,\cc^n)$ for $p>2$ and $\alpha=1-\frac{2}{p}$.
We shall frequently use the following notations: $\d_xu:=\frac{\d u}{\d x}$,
$\d_yu:=\frac{\d u}{\d y}$ and $\dbar u:=\d_xu+i\d_yu$, \ie without $\frac{1}{2}$.

\medskip
Most considerations in this paper are purely local. 
Therefore our framework can be described as follows.
We consider a Lipschitz-continuous matrix valued function $J$ in the unit
ball $B$ of $\cc^n\equiv\rr^{2n}$, \ie $J:B\to ${\slsf Mat}$(2n\times 2n,\rr)$ such that
$J^2(x)\equiv -\id $.  Its Lipschitz constant will be denoted by 
$Lip(J)$. We are studying $J$-holomorphic maps $u:\Delta\to B$. I.e., $u\in \calc^0\cap L^{1,2}(\Delta , B)$
and satisfies 
\begin{equation}
\eqqno(J-holo1)
\dbar_{J\circ u} u \deff \frac{\d u}{\d x} + J(u(z))\frac{\d u}{\d y} =0 
\quad\text{ almost everywhere in }
\Delta .
\end{equation}
We can consider $J(u(z))=(J\circ u)(z)$ as a matrix valued function on the unit disc,
denote it as $J_u(z)$. It satisfies $J_u(z)^2\equiv -\id $ and therefore it can be 
viewed as a complex linear structure on the trivial bundle $E:=\Delta\times\rr^{2n}$.
The mapping $u$ is a section of this bundle. We call the operator $\dbar_{J\circ u}$ the 
$\dbar$-operator for the induced structure $J_u=J\circ u$ on the bundle $E$. 

\smallskip Later in this paper we shall use a similar construction as follows.
In the trivial bundle $E=\Delta\times \rr^{2n} (=\Delta\times\cc^n)$ over the unit disc 
consider a complex structure $J(z)$, 
\ie a continuous $\mat (2n\times 2n,\rr)$-valued function, such that $J(z)^2\equiv -\id$.
It defines on $L^{1,2}_{loc}$ - sections of $E$ a $\dbar$-type operator
\begin{equation}
\eqqno(dbar-type)
\dbar_J u \deff\frac{\d u}{\d x} + J(z)\frac{\d u}{\d y}.
\end{equation}

Therefore we can interpret \eqqref(J-holo1) saying that a $J$-holomorphic map $u$ is a 
section of $E$ , which satisfies \eqqref(dbar-type) with $J(z) = J_u(z)$.

In the Proposition \ref{morrey} below we shall see that $u$ satisfying \eqqref(J-holo1) is,
in fact, of class $\calc^{1,\alpha}$ for all $0<\alpha <1$. 

\begin{prop}
\label{morrey} Let $J$ be  an $\endo (\rr^{2n})$-valued function on $\Delta$
of class $\calc^{k-1,Lip}$, $k\geq1$, and let $R$ be an~$\endo (\rr^{2n})$-valued 
function on $\Delta$ of
class $L^{k,p}$, $1<p<\infty$.  Suppose that $J^2\equiv -\id$ and that $\dbarj u+Ru \in
L^{k,p}(\Delta)$ for some $u\in L^{1, 2}(\Delta,\rr^{2n})$. Then
$u\in L^{k+1,p}_\loc (\Delta,\rr^{2n})$ and for $0<r<1$
\begin{equation}\eqqno(morrey-est)
 \hbox{ } \norm{ u }_{L^{k+1, p}(\Delta(r))}\leq C_{k,p}
\bigl( \norm{ \dbarj u + Ru }_{L^{k, p}(\Delta)} + 
\norm{ u }_{L^p(\Delta)}\bigr),
\end{equation}
where $C_{k,p}=C(\norm{J}_{\calc^{k-1,Lip}}, \norm{ R}_{L^{k,p}}, k, p, r)<\infty$. 
Moreover, there exists an $\eps =\eps (k,p)>0$ such that if
\[
\norm{J-J\st}_{\calc^{k-1,Lip}(\Delta)} +
\norm{R}_{L^{k,p}(\Delta)}<\eps
\]
then the constant $C_{k,p}$ above can be chosen to be independent of
$||J||$ and $||R||$.
\end{prop}
For the proof see \cite{M}, Theorem 6.2.5.  The condition $J^2\equiv -\id$ is 
needed in this statement to insure the ellipticity of the operator $\dbar_J$. 
We shall use in this paper the case $k=1$ only. Remark that our initial 
assumption on $u$ is $u\in L^{1,2}(\Delta)$ which implies that $u\in L^p(\Delta)$ 
for all  $p<\infty$. This proposition implies, in
particular, that a $J$-holomorphic map $u:\Delta\to\rr^{2n}$ is of class
$L^{2,p}_\loc(\Delta)$ for all $p<\infty $ provided $J$ is Lipschitz. In 
particular $u\in \calc^{1,\alpha}_{loc}(\Delta)$ for all $0<\alpha <1$.

\newprg[prgZD.cusp]{Estimation of the differential at cusp-points}

Throughout this subsection we fix some $2<p<\infty$ and make the 
following assumption:

\begin{itemize}
\item[$(*)$]
{\it $J$ is an almost complex structure in $B$ with
 $J(0)=J\st$ such that $\norm{J-J\st}_{\calc^{Lip}(B)}$ is small enough.
$u:\Delta\to B$ is a $J$-holomorphic map such that $u(0)=0$ and such that
$\norm{du}_{L^{1,p}(\Delta)}$ is small enough.}
\end{itemize} 

Let us notice that this assumption is by no means restrictive. Indeed, we can always
replace $J(w)$ by $J_\tau(w):= J(\tau w)$ and $u(z)$ by  $u_{t,\tau}(z):=\tau^{-1} u(tz)$ with
some appropriately chosen $\tau$ and $t$.

\medskip%
By the Corollary 3.1.3 from \cite{IS1} (see also Proposition 3 in \cite{Sk2}
and the corresponding  Corollary 1.4.3 in \cite{IS2}) 
we can assign multiplicity of zero to a $J$-holomorphic map $u:\Delta\to (B,J)$
provided that $J$ is at least Lipschitz. In particular, zeroes of $u$ are
isolated, as for the classical holomorphic functions. Moreover, we can
represent such $u$ (in the neighborhood of its zero point, say $z_0=0$,
provided $J(0)=J\st$) as
\begin{equation}
u(z)=z^\mu P(z) +z^{2\mu - 1}v(z), \eqqno(norm-form-map)
\end{equation}
where $\mu\geq 1$ is an integer (a multiplicity of zero), $P(z)$ is
some (holomorphic) polynomial of degree at most $\mu -1$, $P(0)\not=
0$ and 
$v\in L^{1,p}_{loc}(\Delta, \cc^n)$  for all $2<p<\infty$, and therefore  $v\in \calc^{\alpha
}(\Delta ,\cc^n)$ for all $0<\alpha<1$. 
In addition,  $v(0)=0$. Now we
want to derive from \eqqref(norm-form-map) some properties of the
differential $du$. 

\medskip
Let us start with the following preliminary estimate.  

\begin{lem} 
\label{v-est}
For any integer $\mu \geq 1$ there exists a constant $C= C(\mu,p)<\infty$ with the 
following property: for every $J$-holomorphic map $u:\Delta \to (B,J)$, satisfying 
the assumption (*) and having the form \eqqref(norm-form-map) one has
\begin{equation}\eqqno(v1-p)
\norm{v}_{L^{1,p}(\Delta)} \leq C\cdot \norm{u}_{L^{1,p}(\Delta)}.
\end{equation}
\end{lem}
\proof We use the fact that every $J$-holomorphic map $u:\Delta\to B$ with 
Lipschitz-continuous $J$ satisfies the pointwise estimate 
\begin{equation}\eqqno(dbar-u-du)
|\dbar\st u(z)| \leq Lip(J) \cdot |du(z)|\cdot|u(z)|,
\end{equation} 
see inequality  (1.4.4)  in \cite{IS2}. Following \cite{Sk1} define 
\begin{equation}
\eqqno(deff-H)
H(z) := 
\begin{cases}\textstyle
\ - \frac{\dbar\st u(z)\otimes \bar u(z)}{|u(z)|^2}&  \text{ if } \quad u(z)\not= 0,\\
\ 0 & \text{ if } \quad u(z) = 0.
\end{cases}
\end{equation}
Then $H(z)$ is a measurable function with values in $\mat_{\cc}(n\times n)$,
which satisfies the pointwise estimate 
\begin{equation}
\eqqno(point-est)
| H(z)| \leq Lip(J) \cdot |du(z)|.
\end{equation}
In particular, $\norm{H(z)}_{L^p(\Delta)}$ is bounded by some sufficiently small
constant by the assumption $(*)$.  The function  $u(z)$ in 
its turn satisfies 
\begin{equation}
\eqqno(sikor)
\dbar\st u(z) + H(z)\cdot  u(z)=0.
\end{equation}
Under these conditions Lemma 1.2.3 from
\cite{IS2} insures the existence of a matrix-valued function $F(z)\in
L^{1,p}(\Delta)$ %
which satisfies the equation 
\begin{equation}
\eqqno(iv-sh2)
\dbar\st F= - F\cdot H.
\end{equation} 
with the estimate
\begin{equation}
\eqqno(iv-sh1)
\norm{F(z)-\id}_{L^{1,p}(\Delta)}\leq C\cdot
\norm{H(z)}_{L^p(\Delta)}
\end{equation}

The equations \eqqref(sikor) and \eqqref(iv-sh2) imply that the function
$F(z){\cdot}u(z)$ is holomorphic. Define $u^{(\mu)}(z) := z^{-\mu}u(z)$. It satisfies 
$\dbar\st (Fu^{(\mu)}(z))=0$ and relation \eqqref(norm-form-map) implies
that the function $F(z){\cdot}u^{(\mu)}(z)= z^{-\mu}F(z){\cdot}u(z)$ has no pole at zero
and therefore is holomorphic.

Since $\norm{F(z)-\id}_{L^{1,p}(\Delta)}$ is small we have that for any domain 
$A\subset \Delta$ the following inequality
\[
c_{\mu}\norm{u^{(\mu)}}_{L^{1,p}(\Delta)} \leq \norm{F(z)\cdot u^{(\mu)}(z)}_{L^{1,p}(\Delta)}\leq 
C_{\mu}\norm{u^{(\mu)}}_{L^{1,p}(\Delta)}.
\]
Observe that $F(z)\cdot u^{(\mu)}(z)$ satisfies 
\begin{equation}
\eqqno(max-pr)
\norm{F\cdot u^{(\mu)}}_{L^{1,p}(\Delta)}\leq 
C_1\norm{F\cdot u^{(\mu)}}_{L^{1,p}(\Delta\bs \Delta(\half))}.
\end{equation}
This easily follows from the Cauchy integral formula for the Taylor
coefficients of holomorphic function.  Consequently we get
\begin{equation}
\norm{F\cdot u^{(\mu)}}_{L^{1,p}(\Delta)}\leq C\cdot
\norm{u^{(\mu)}}_{L^{1,p}(\Delta\bs \Delta(\half)) }.
\end{equation}
This gives us the estimate 
\begin{equation}\eqqno(u-mu-p)
\norm{u^{(\mu)}}_{L^{1,p}(\Delta)}\leq C_{\mu}\cdot
\norm{u}_{L^{1,p}(\Delta) },
\end{equation} 
because on the annulus $\Delta\bs \Delta(\half)$ functions $u$ and $u^{(\mu)}$ 
are comparable. The latter estimate implies that 
$d(u(z))= d(z^\mu u^{(\mu)}(z))$ fulfills a.-e. the  estimate 
\begin{equation}\eqqno(du-mu-p)
|du(z)|\leq h(z)\cdot |z^{\mu-1}| 
\end{equation}
for some non-negative $L^p$-function $h$ satisfying 
\begin{equation}
\eqqno(h-lp)
\norm{h}_{L^p(\Delta)}\leq C\cdot\norm{u}_{L^{1,p}(\Delta)}.
\end{equation} 
Substituting the relation \eqqref(du-mu-p) in \eqqref(dbar-u-du) we obtain 
a.-e.\ the pointwise estimate
\[
|\dbar\st u(z)|\leq C\cdot|z^{\mu-1}|\cdot h(z)\cdot|u(z)|
\]
with the same function $h\in L^p(\Delta)$ as above. Multiplying it by $z^{-\mu}$ we obtain 
\begin{equation}
\eqqno(debar-mu)
|\dbar\st u^{(\mu)}(z)|\leq C\cdot|z^{\mu-1}|\cdot h(z)\cdot|u^{(\mu)}(z)|.
\end{equation}

\medskip For $j=1,\ldots,\mu-1$ define functions $u^{(\mu+j)}(z)$ recursively by the relation 
\[
u^{(\mu+j)}(z):= (u^{(\mu+j-1)}(z) -u^{(\mu+j-1)}(0))/z.
\]
Then the coefficients of the polynomial $P(z)$ from \eqqref(norm-form-map) are
given by $a_j=u^{(\mu+j)}(0)$ for $j=0,\ldots,\mu-1$ and $v(z)=u^{(2\mu-1)}(z)
-u^{(2\mu-1)}(0)$.  

We claim that for every $j=0,\ldots,\mu-1$ we have  the  estimation
\begin{equation}
\eqqno(mu+j)
\norm{u^{(\mu+j)}}_{L^{1,p}(\Delta)}\leq C\cdot
\norm{u}_{L^{1,p}(\Delta) },
\end{equation}
and therefore $|a_j| \leq C\cdot \norm{u}_{L^{1,p}(\Delta) }$. The proof is done by induction 
using \eqqref(u-mu-p) as the base for $j=0$.
Thus we assume that for some fixed $j\in\{1,\ldots,\mu-1\}$ we have the estimation of the
form $\norm{u^{(\mu+j-1)}}_{L^{1,p}(\Delta)}\leq C\cdot \norm{u}_{L^{1,p}(\Delta) }$, and in
particular $|a_{j-1}|=|u^{(\mu+j-1)}(0)|\leq C\cdot \norm{u}_{L^{1,p}(\Delta) }$. From the definition 
of $u^{(\mu+j)}(z)$ we obtain a.-e.\ the pointwise differential inequality
\[
\begin{split} 
|\dbar\st u^{(\mu+j)}(z)|\; &= |z\inv \dbar\st u^{(\mu+j-1)}(z)| = 
|z^{-2} \dbar\st u^{(\mu+j-2)}(z)|=\ldots=
\\
&= |z^{-j} \dbar\st u^{(\mu)}(z)|\leq C\cdot h(z)\cdot |z^{\mu-1-j}|\cdot |u^{(\mu)}(z)|.
\end{split} 
\]
This gives the estimate 
\begin{equation}
\norm{\dbar\st u^{(\mu+j)}(z)}_{L^p(\Delta)}\leq C\norm{h}_{L^p(\Delta)}\norm{u^{(\mu)}}_{L^{\infty}(\Delta)}
\leq C_1\norm{u}_{L^{1,p}(\Delta)}
\end{equation}
by \eqqref(u-mu-p) and \eqqref(h-lp). Further 
\[
\norm{ u^{(\mu+j)} }_{L^{1,p}(\Delta\bs \Delta(\half)) }  
\leq C\cdot \big( a_{j-1} +
\norm{ u^{(\mu+j-1)} }_{L^{1,p}(\Delta\bs \Delta(\half)) }\big)
 \leq C_1\cdot \norm{u}_{L^{1,p}(\Delta) }
\]
by inductive assumption. Now the standard inner estimates for $\dbar\st$  provide the desired
estimates
\[
\norm{u^{(\mu+j)}}_{L^{1,p}(\Delta)}\leq C\cdot
\norm{u}_{L^{1,p}(\Delta) } 
\qquad\text{and}\qquad
|a_j| \leq C\cdot \norm{u}_{L^{1,p}(\Delta) }.
\]
The case $j=\mu-1$ yields the 
estimate \eqqref(v1-p) on  $v(z)=u^{(2\mu-1)}(z)-a_{\mu-1}$.  

\smallskip\qed

\smallskip The following lemma will be
used in this paper with various operators $A$.

\begin{lem}
\label{lipschitz-A}
Let $A$ be a Lipschitz continuous $\endo (\rr^{2n})$-valued function on $B$
with $A(0)=0$ and let $u$ be a $J$-holomorphic map with Lipschitz $J$ 
as in \eqqref(norm-form-map). Then for all integers $\nu$ and $\lambda$ satisfying 
$\nu \leq \mu + \lambda - 1$ the function $z^{-\nu}\cdot A(u(z))\cdot
z^{\lambda}$ is Lipschitz-continuous in $\Delta$ with the estimate 
\begin{equation}\eqqno(lip-A)
Lip_{\Delta}( z^{-\nu}\cdot A(u(z))\cdot z^{\lambda} ) \leq C(p)\cdot Lip(A)\cdot \norm{u}_{L^{1,p}(\Delta)}.
\end{equation}
\end{lem}

\proof Remark that we have the following estimate
\begin{equation}
\eqqno(u-mu-est)
\norm{z^{-\mu} u(z)}_{L^\infty(\Delta(\frac23))} \leq C\cdot \norm{u}_{L^{1,p}(\Delta)}.
\end{equation} 
\smallskip\noindent This is clear because in \eqqref(u-mu-p) a $L^{1,p}$ and therefore a $\calc^{\alpha}$ -
norm of $u^{(\mu)}=z^{-\mu}u$ was estimated.

\smallskip 
We continue the proof of the lemma  starting from the remark that  since $u$ is $J$-holomorphic with Lipschitz $J$, it
is of class $\calc^{1,\alpha}$ and therefore itself Lipschitz. 

\smallskip Now we turn to the estimation of the quantity $\frac{\norm{z^{-\nu}_1
A(u(z_1))z^{\lambda}_1-z^{-\nu}_2 A(u(z_2)) z^\lambda_2}}{|z_1-z_2|}$ for $z_1\not= z_2\in \Delta$.  
In order to do  this we consider two cases.

\smallskip\noindent{\slsf Case 1.} $\frac{1}{3}|z_1| \leq  |z_1-z_2|$.  In that case
$|z_1|\leq 3|z_1-z_2|$ and $|z_2|\leq 4|z_1-z_2|$. Therefore

\[
\frac{\norm{z^{-\nu}_1A(u(z_1))z^{\lambda}_1-z^{-\nu}_2A(u(z_2))z^{\lambda}_2}}{|z_1-z_2|}
=
\]
\[
= \frac{\norm{z^{-\nu}_1A(u(z_1))z^{\lambda}_1 -
z^{-\nu}_1A(u(0))z^{\lambda}_1 + z^{-\nu}_2A(u(0))z^{\lambda}_2-
z^{-\nu}_2A(u(z_2))z^{\lambda}_2}}{|z_1-z_2|} \leq
\]
\[
\leq C\cdot Lip(A) \norm{z^{-\mu}u(z)}_{L^\infty}
\frac{|z_1|^{\mu +\lambda -\nu} + |z_2|^{\mu +\lambda -\nu}}{|z_1-z_2|}\leq 
\]
\[
\leq C\cdot Lip(A) \norm{z^{-\mu}u(z)}_{L^\infty}\frac{|z_1|+|z_2|}{|z_1-z_2|}
\leq C\cdot Lip(A) \norm{u(z)}_{L^{1,p}(\Delta)},
\]
because $\mu +\lambda -\nu  \geq 1$. In the second line we silently used the fact that
$A(u(0))=0$. 

\smallskip\noindent{\sl Case 2.} $|z_1-z_2|\leq \frac{1}{3}|z_1|$.  In that case
$\frac{2}{3}|z_1|\leq |z_2|\leq \frac{4}{3}|z_1|$. Therefore

\[
\frac{\norm{z^{-\nu}_1A(u(z_1))z^{\lambda}_1-z^{-\nu}_2A(u(z_2))z^{\lambda}_2}}{|z_1-z_2|}\leq
\frac{\norm{z^{-\nu}_1A(u(z_1))z^{\lambda}_1-z^{-\nu}_1A(u(z_1))z^{\lambda}_2}}{|z_1-z_2|}
+
\]
\[
+ \frac{\norm{z^{-\nu}_1A(u(z_1))z^{\lambda}_2 -
z^{-\nu}_2A(u(z_1))z^{\lambda}_2}}{|z_1-z_2|} +
\frac{\norm{z^{-\nu}_2A(u(z_1))z^{\lambda}_2 -
z^{-\nu}_2A(u(z_2))z^{\lambda}_2}}{|z_1-z_2|} =
\]
\[
=
\frac{\norm{z^{-\nu}_1A(u(z_1))[z^{\lambda}_1-z^{\lambda}_2]}}{|z_1-z_2|}
+ \frac{\norm{[z^{-\nu}_1 - z^{-\nu}_2]
A(u(z_1))z^{\lambda}_2}}{|z_1-z_2|} +
\frac{\norm{z^{-\nu}_2[A(u(z_1)) - A(u(z_2))]
z^{\lambda}_2}}{|z_1-z_2|} =
\]
\[
= I_1 + I_2 + I_3.
\]
Let us estimate these terms separately. First:
\[
I_1 \leq |z_1|^{-\nu}\norm{A(u(z_1))}\sum_{j=0}^{\lambda-1}|z_1^jz_2^{\lambda -j-1}|
\]
\[
\leq C\cdot Lip(A) \norm{z^{-\mu} u(z)}_{L^\infty(\Delta(\frac23))} |z_1|^{-\nu +\mu}|z_1|^{\lambda
-1}\leq C\cdot Lip(A) \norm{u}_{L^{1,p}(\Delta)} .
\]
Second:
\[
I_2\leq C\cdot Lip(A)  \norm{z^{-\mu} u(z)}_{L^\infty(\Delta(\frac23))} \; |z_1|^{\mu } \; |z_2|^{\lambda }
\left|\frac{z_2-z_1}{z_1z_2}\right| \left(\sum_{j=0}^{\nu -1}
\left|z_1^{-j}z_2^{-(\nu -j -1)}\right|\right) \frac{1}{|z_1-z_2|}
\]
\[
\leq
  C\cdot Lip(A) \norm{u}_{L^{1,p}(\Delta)} |z_1|^{\lambda +\mu -\nu -1}\leq C\cdot Lip(A) \norm{u}_{L^{1,p}(\Delta)} .
\]
And, finally, third:
\[
I_3\leq C\cdot Lip(A) \norm{z^{-\mu} u(z)}_{L^\infty(\Delta(\frac23))} |z_2|^{-\nu}|z_1^{\mu}-z_2^{\mu}|\; |z_2|^{\lambda} \;
\frac{1}{|z_1-z_2|}\leq C\cdot Lip(A)\norm{u}_{L^{1,p}(\Delta)}.
\]

\smallskip\qed

We need to produce some extra regularity of the rest term $zv(z)$ in the
representation \eqqref(norm-form-map). In fact we shall prove that 
$z\cdot v\in L^{2,p}_\loc(\Delta)$ for $v$ from \eqqref(norm-form-map) together with the
estimate of its decay at zero.

\begin{lem}
\label{l2p-reg}
Let $J$ be a Lipschitz-continuous almost complex structure in the
unit ball $B \subset \cc^n$ with $J(0)=J\st$ and  $u:\Delta \to B$
be a $J$-holomorphic map written in the form \eqqref(norm-form-map).
Then $zv\in L^{2,p}_{\loc}$ for all $2<p<\infty$. Moreover,  for every
$|z|\leq\frac{1}{2}$ one has
\begin{equation}
| d(zv(z))| \leq C(p) \cdot
|z|^{1-\frac{2}{p}}\cdot \Vert u \Vert_{L^{1,p}(\Delta)} \eqqno(2.2)
\end{equation}
for with the constant $C(p)$ independent of $z$, $u$ and $J$ satisfying the assumption (*).
\end{lem}

\proof Remark that with $\nu = 2\mu - 2$ in Lemma \ref{lipschitz-A} above we have
\begin{equation}
z^{\nu}\cdot (zv(z))=u(z)-z^{\mu}\cdot P(z). \eqqno(2.3)
\end{equation}
So if we apply $\dbar_{J\circ u}$ to the right hand side of \eqqref(2.3) we obtain
$\dbar_{J\circ u}(z^{\mu}P)$ which is Lipschitz continuous, in particular, it
belongs to $L^{1,p}_{loc}$ for all $1<p<\infty$. Elliptic regularity of Proposition
\ref{morrey} gives then the $L^{2,p}$-regularity of the left hand side $z^{\nu}\cdot
(zv(z))$.  
If $\mu =1$ (and therefore $\nu =0$), then $zv(z)=u(z)\in
L^{2,p}_\loc(\Delta)$, and 
the needed $L^{2,p}$-regularity is already proved.
Therefore till the {\sl Step 4} we shall suppose that $\mu\geq 2$.

\smallskip
Let us explain the idea of the proof of this lemma. First we observe that
\begin{equation}
z^{-\nu}\left[ \d_x +
J(u(z))\d_y\right]z^{\nu}(zv)=(\d_x+z^{-\nu}J(u)z^{\nu}\d_y)(zv)+
z^{-\nu}(1+J(u)J\st)\nu z^{\nu -1}(zv). \eqqno(2.4)
\end{equation}
We see \eqqref(2.4) as the equation of the form 
\begin{equation}
\eqqno(de-bar-nu)
f(z) = \left( \dbar_{J^{(\nu)}} + R^{(\nu)} \right ) (zv(z)).
\end{equation}
After establishing the necessary regularity of $J^{\nu)}$, $R^{(\nu)}$ and $f$ we shall apply the
Proposition \ref{morrey} and obtain the desired regularity of the solution $zv$.

\smallskip
Let us start with the right hand side of \eqqref(2.4).

\smallskip\noindent%
{\sl Step 1.} {\it $J^{(\nu)}\deff z^{-\nu}\cdot J(u)\cdot z^{\nu}$ is a Lipschitz
continuous complex structure on $E$ and $Lip (J^{(\nu )})\leq C\cdot Lip (J)\cdot 
\norm{u}_{L^{1,p}(\Delta)}$. }

The proof is straightforward via Lemma \ref{lipschitz-A}: just write
$J^{(\nu)}=z^{-\nu}[J(u)-J\st ] z^{\nu} + J\st$ and apply Lemma \ref{lipschitz-A} to
$A=J-J\st$.

\smallskip\noindent{\sl Step 2.} {\it The endomorphism
 $R^{(\nu)} \deff  z^{-\nu}\cdot(1 + J(u)J\st)\nu z^{\nu - 1}$ of the bundle $E$ is 
Lipschitz continuous, $R^{(\nu)}(0)=0$ and $Lip (R^{(\nu )})\leq C\cdot Lip (J)\cdot 
\norm{u}_{L^{1,p}(\Delta)}$.}

This is again true by Lemma \ref{lipschitz-A} and because $\mu\geq 2$. Note now
that the right hand side of \eqqref(2.4) is of the form $\dbar_{J^{(\nu)}}(zv) +
R^{(\nu)} (zv)$ and that coefficients of this operator are Lipschitz continuous.
Therefore we can apply \eqqref(morrey-est) and obtain 
\begin{equation}\eqqno(morrey-zv)
\norm{ zv }_{L^{2,p}(\Delta(1/2))} \leq C
\big( \norm{\dbar_{J^{(\nu)}}(zv) + R^{(\nu)}(zv) }_{L^{1,p}(\Delta)} 
+ \norm{zv}_{L^p(\Delta)} \big).
\end{equation}
To achieve \eqqref(2.2) for $r=\frac12$ we need to estimate both terms
in the right hand side \eqqref(morrey-zv) by $\norm{u}_{L^{1,p}(\Delta)}$. 
For the term $zv(z)$ it was already done in \eqqref(v1-p). In order to
estimate the first term we shall compute the left
hand side of \eqqref(2.4) in another way. Namely, using
\eqqref(norm-form-map) we write
$$
z^{-\nu}\left[ \d_x + J(u(z))\d_y\right]z^{\nu}(zv)=z^{-\nu}\left[
\d_x + J(u(z))\d_y\right](u(z) - z^{\mu}\cdot P(z)) =
$$
$$
=z^{-\nu}\left[ \d_x + J(u(z))\d_y\right](-z^{\mu}\cdot
P(z))=z^{-\nu}\left[ \d_x + J(u(z))\d_y\right] (\sum_{j=0}^{\mu -1}
a_jz^{\mu +j}) =
$$
\begin{equation}
= z^{-\nu} (1 + J(u)J\st )(\sum_{j=0}^{\mu -1} (\mu +j)a_jz^{\mu
+j-1}) =: f(z). \eqqno(2.5)
\end{equation}

\smallskip\noindent{\sl Step 3.}  {\it The right hand side $f(z)$ of
 \eqqref(2.5) satisfies $f(0)=0$ and is Lipschitz continuous with the
 estimate}
\begin{equation}\eqqno(lip-f)
\norm{f}_{\calc^{Lip}(\Delta)} \leq C
\norm{J}_{\calc^{Lip}(\Delta)} \norm{u}_{L^{1,p}(\Delta)}.
\end{equation}

The worst term is $z^{-\nu} (1 + J(u)J\st )(\mu a_0z^{\mu -1})$, 
but it is still under control of the  Lemma \ref{lipschitz-A}. 
We conclude now by \eqqref(morrey-zv) and \eqqref(lip-f)
the estimate 
\begin{equation}
\eqqno(d-z-v)
\norm{d(zv)}_{L^{1,p}(\Delta (1/2))}\leq C\cdot \norm{u}_{L^{1,p}(\Delta)}.
\end{equation}

\smallskip\noindent{\sl Step 4.}  
{\it The behavior of $z\,v(z)$ under a dilatation.}

For $\tau \in [0,1]$ and we define $\pi_\tau: \Delta\to\Delta$ by $\pi_\tau(z) := \tau\cdot z$. Then for any
function $w(z)$ in the disc $\pi_\tau ^*w(z) = w(\tau z)$ is the dilatation of $w(z)$.
An easy calculation shows the following \emph{dilatation behavior} of the
$L^p$-norms of the derivatives:
\[
\norm{ D^i(\pi_\tau ^*w)}_{L^p(\Delta) } = 
\tau^{|i|- \frac{2}{p}}
\norm{ D^iw }_{L^p(\Delta(\tau)) }.
\]

Estimating $L^p$-norm we use the fact that $v(0)=0$ and that
$\norm{v}_{\calc^\alpha(\Delta)}\leq C\cdot \norm{u(z)}_{L^{1,p}(\Delta)}$, see Lemma \ref{v-est}. 
Hence $|v(z)|\leq C\cdot |z|^\alpha\cdot
\norm{u(z)}_{L^{1,p}(\Delta)}$ and thus (using $\alpha=1-\frac{2}{p}$)
\[
\norm{zv(z)}_{L^p(\Delta(r)) } \leq \norm{zv(z)}_{L^\infty(\Delta(r))
}  (\pi r^2)^{\frac{1}{p} } 
\leq C  \cdot r^{2}\cdot \norm{u(z)}_{L^{1,p}(\Delta) }.
\]
Consequently
\[
\norm{\pi^*_r(zv(z))}_{L^p(\Delta) } 
\leq C  \cdot r^{2-\frac2p}\cdot \norm{u(z)}_{L^{1,p}(\Delta) }.
\]

Now recall that $zv(z)$ satisfies the differential equation
\[
\left( \dbar_{J^{(\nu)}} + R^{(\nu)} \right ) (zv(z)) = f(z)
\]
with $f(z)$ given by the formula \eqqref(2.5). By {\sl Step 3}, $f(z)$
is Lipschitz continuous with the estimate \eqqref(lip-f) and $f(0)= 0$.
This implies that
\[
\norm{\pi_r^* f(z)}_{\calc^{Lip}(\Delta)} \leq C\cdot r\cdot
\norm{J}_{\calc^{Lip}(\Delta)} \norm{u(z)}_{L^{1,p}(\Delta)}.
\]
The same argument yields also
\[
\norm{\pi_r^* (J^{(\nu)} -J\st)}_{\calc^{Lip}(\Delta)} \leq r\cdot
\norm{J^{(\nu)} -J\st}_{\calc^{Lip}(\Delta)}.
\]
Finally, we observe that
\[
\pi^*_r \big( \dbar_{J^{(\nu)}} (zv(z)) \big) = 
r^{-1} \cdot \dbar_{\pi^*_r J^{(\nu)}} (\pi^*_r (zv(z)) ).
\]
Summing up, we see that the rescaled function 
$w_r(z) := \pi^*_r (zv(z))$ satisfies  a $\dbar$-type equation
\[
\big( \dbar_{\pi^*_r J^{(\nu)}}  
 + r\cdot \pi_r^* R^{(\nu)} \big) w_r(z) =
r\cdot \pi_r^* f(z)
\]
in which the norms $\norm{r\cdot \pi_r^* f(z)}_{\calc^{Lip}(\Delta)}$ are bounded by
$C\cdot r^2\cdot\norm{u}_{L^{1,p}(\Delta)}$ uniformly in $r$ and the coefficients
$\pi^*_rJ^{(\nu)}$ and $r\cdot \pi^*_r R^{(\nu)}$ are $\calc^{Lip}$-close to those of
$\dbar\st$ for $r$ close enough to $0$. 

From Proposition \ref{morrey} we
obtain the uniform estimate
\[
\norm{\pi^*_r ( zv(z))}_{L^{2,p}(\Delta)} \leq
 C \cdot r^{2-\frac2p} \cdot \norm{u}_{L^{1,p}(\Delta)},
\]
which implies
\[
\norm{d(\pi^*_r ( zv(z)))}_{L^\infty(\Delta)}
\leq C \cdot r^{2-\frac2p} \cdot \norm{u}_{L^{1,p}(\Delta)}.
\]
After rescaling we obtain 
\[
\norm{d( zv(z))}_{L^\infty(\Delta(r))}
\leq C \cdot r^{\alpha} \cdot \norm{u}_{L^{1,p}(\Delta)},
\]
with $\alpha=1-\frac{1}{p}$. \qed

\smallskip%
The estimate \eqqref(2.2) gives immediately the following

\begin{corol}
For a Lipschitz-continuous $J$ and $J$-holomorphic $u$  in the form
\eqqref(norm-form-map) one has
\begin{equation}
du(z) = d(z^\mu P(z)) + O(|z|^{2\mu -2 +\alpha}).
\eqqno(norm-form-diff)
\end{equation}
for any $0<\alpha <1$. In particular, for a non-constant $u$ zeroes of $du$ are isolated.
\end{corol}

\newprg[prgZD.examp]{An example}

Let us illustrate the statements of this Section by an example.

\begin{exmp} \rm 
Equation \eqqref(1.1) can be rewritten as

\begin{equation}
\frac{\d u}{\d\bar z} - \bar Q(J_u(z))\frac{\d u}{\d z} = 0,
\eqqno(q-bar)
\end{equation}
where
\begin{equation}
\eqqno(du-dbaru)
\frac{\d u}{\d\bar z}=\frac{1}{2}\big(\frac{\d u}{\d x}+J\st \frac{\d
u}{\d y}\big), \qquad\frac{\d u}{\d z}=\frac{1}{2}\big(\frac{\d u}{\d
x}-J\st \frac{\d u}{\d y}\big)
\end{equation}
and
\begin{equation}
\bar Q(J(z))=[J(z)+J\st ]^{-1}[J\st -J(z)]. 
\eqqno(q-bar-j)
\end{equation}
Remark that $\bar Q$ anticommutes with $J\st$ and therefore is a
$\cc$-antilinear operator. Therefore \eqqref(q-bar) can be understood
as an equation for $\cc^n$-valued map (or section) $u$. Usually it
is better to consider the conjugate operator $Q$ and write
\eqqref(q-bar) in the form

\begin{equation}
\frac{\d u}{\d\bar z} - Q(J_u(z))\overline{\frac{\d u}{\d z}} = 0.
\eqqno(J-holo2)
\end{equation}
\noindent Vice versa, given an anti-linear operator $\bar Q$ in $\cc^n\equiv\rr^{2n}$, 
one can reconstruct the corresponding almost complex structure as follows
\begin{equation}
\eqqno(j-q-bar)
J(z) = J\st \big(\id + \bar Q\big)\cdot\big(\id - \bar Q\big)^{-1}.
\end{equation}
After these preliminary considerations (which will be used also in Section \ref{sectOPT}),
we shall turn to the example in question. Remark that the vector-function 
$u(z) = (z^{\mu}, \bar z^{2\mu})$ is $J$-holomorphic with respect to the structure
\begin{equation}
Q(u_1,u_2) = 
\begin{pmatrix}
0  &  0 \cr
2\bar u_1& 0 
\end{pmatrix}
\end{equation}
In the representation \eqqref(norm-form-map) we have for this example $P=(1,0)$ - a constant
vector polynomial, $v(z) = \bar z$. From here one sees that \eqqref(2.2) cannot be improved.
\end{exmp}

\begin{rema} {\bf (a)} \rm
The fact that for $J\in\calc^{1}$ zeroes of a differential of a
$J$-holomorphic map are isolated was first proved by J.-C. Sikorav in
 \cite{Sk1}.

\smallskip\noindent{\bf (b)}  We shall crucially need this fact for Lipschitz-continuous 
structures in this paper. It is stated in Proposition 3 of \cite{Sk2}, but, unfortunately, 
the proof of \cite{Sk2} uses the expression $(dJ.f')$, see the first line after the
formula (2.3) on page 363 of \cite{Sk2}. Here the Author means the pointwise scalar product 
$(dJ(f), f')$. But $dJ(f)$ cannot be defined for Lipschitz $J$ and no explanations of how one 
might give the sense to this expression are given. Therefore, in our opinion, the proof of
Proposition 3 of \cite{Sk2} goes through  only in the case $J\in \calc^1$ and this was already 
achieved in \cite{Sk1}. Remark that this problem - absence of the chain rule for Lipschitz maps -
makes the issue quite delicate. 
\end{rema}

\newsect[sectLOC]{Local structure of $J$-holomorphic maps}

\newprg[prgLOC.uniq]{Uniqueness for solutions of $\dbar$-inequalities}

We start with a generalization of Lemma 1.4.1 from \cite{IS2}. Recall that for a complex 
structure $J(z)$ in the trivial bundle
$E=\Delta\times \rr^{2n} (=\Delta\times\cc^n)$ over the unit disc we defined 
the operator $\dbar_J$ by the formula  $\dbar_J u = \frac{\d u}{\d x} + 
J(z)\frac{\d u}{\d y}$, see \eqqref(dbar-type).

\begin{lem} Let $J$ be an almost complex structure in the trivial
$\cc^n$-bundle over the disc $\Delta$ which is $L^{1,p}$-regular for
some $2<p<\infty$ and such that $J(0) =J\st$. Suppose that a
function $u\in L^{1,2}_\loc(\Delta, \cc^n)$ is not identically $0$
and satisfies {\sl a.e.} the inequality
\begin{equation}
| \dbar_J u | \leq h\cdot | u | \eqqno(3.1)
\end{equation}
for some nonnegative $h\in L^p_\loc(\Delta)$. Then:

\smallskip
\sli $u\in L^{1, p}_\loc(\Delta)$, in particular $u\in
\calc^{\alpha}_\loc (\Delta)$ with $\alpha\deff 1-\frac{2}{p}$;

\slii for any $z_0\in\Delta$ such that $u(z_0)=0$ there exists
$\mu\in\nn$---the~multiplicity of zero of $u$ in $z_0$---such that
$u(z)=(z-z_0)^\mu \cdot g(z)$ for some $g\in L^{1, p}_\loc (\Delta)$
with $g(z_0)\not=0$.
\end{lem}

\proof We reduce the case of general $J$ to the special one in which $J=
J\st$. 
For this purpose we fix a $(J\st, J)$-complex bundle isomorphism $F: \Delta\times\cc^n \to
\Delta\times\cc^n$ of regularity $L^{1,p}$, so that $F\inv \circ J \circ F = J\st$. Then any
section $u(z)$ of $\Delta\times\cc^n$ has the form $u(z) = F(v(z))$ and $u(z)$ is
$L^{1,p}$-regular if and only if so is $v(z)$. Moreover,
\[
\begin{split}
\dbar_J u(z) &= (\partial_x + J(z)\partial_y )F(v(z)) =\\
&=F\big( \partial_x + F\inv \cdot J(z) \cdot F \;\partial_y \big)\,v(z)
\;+\; (\partial_x F+ J(z)\partial_y F)\,v(z).
\end{split}
\]
Consequently, \eqqref(3.1) is equivalent to the differential inequality
\[
| \dbar\st v | \leq | F\inv( \dbar_J u(z) )| 
+  |F\inv (\partial_x F+ J(z)\partial_y F)\,v(z) | \leq
\]
\begin{equation}\leq  h\cdot | F\inv ||u| + 
 |F\inv (\partial_x F+ J(z)\partial_y F)|\,v(z)|
 \leq h_1\cdot | v | \eqqno(3.2)
\end{equation}
with a new $h_1\in L^p(\Delta)$.

\medskip

The statement of the lemma  is reduced now to Lemma 1.4.1 from
\cite{IS2}.

\medskip
\qed

\begin{lem}
\label{lem3.2}
Let $J$ be a Lipschitz-continuous almost complex structure in the unit ball
$B$ in $\cc^n$ and $u_1, u_2: \Delta \to B$ two $J$-holomorphic maps such that
$u_1(0) = u_2(0) =0$ and $u_1 \not \equiv u_2 $. Then there exists an integer $\nu
>0$ 
and $v(z)\in L^{1,p}(\Delta, \cc^n)$, $v(0) \neq 0$ such that $u_1(z) -u_2(z) = z^\nu
v(z)$.
\end{lem}

\proof Set $v=u_1-u_2$ and let us compute $\dbar_{J\circ u_1}(v) =
(\partial_x + J(u_1(z))\partial_y)v(z)$:
\[
\dbar_{J \circ u_1}(v) = (\partial_x + J(u_1)\cdot \partial_y) (u_1
-u_2) = (\partial_x + J(u_1)\cdot \partial_y) (u_1 -u_2) +
(\partial_x + J(u_2)\cdot \partial_y) u_2 =
\]
\[
= (J(u_2)\cdot
\partial_y - J(u_1)\cdot
\partial_y)u_2 = (J(u_1 -v) - J(u_1))\cdot
\partial_yu_2.
\]
By the Lipschitz regularity of $J$ and $\partial_yu_2 \in L^p(\Delta)$ we obtain a pointwise
differential inequality
\[
|\dbar_{J \circ u_1}(v)(z)| \leq h(z)\cdot |v(z)|
\]
for some $h\in L^p(\Delta)$. Now we apply  Lemma 3.1.
\qed

\begin{rema}\rm
The statement of this lemma implicitly  appeared for the first time in
\cite{Sk1}. 
Really, his proof of Proposition 3.2.1 (i) clearly goes through under the
assumption of Lipschitz continuity of $J$  only. 
\end{rema}

\newprg[prgLOC.compa]{Proof of the part {\bf (a)} of the Comparison Theorem.}
In  the proof we use the abbreviation ``$L^{1,p}$-regular'' instead of ``$L^{1,p}$-regular
for any $p<\infty$'' and a similar abbreviation ``$L^{2,p}$-regular''.

\medskip\noindent{\bf (a)}
Denote by $E$ the a  trivial
$\cc^n$-bundle over $\Delta$,  $E:=\Delta\times\cc^n$. Equip $E$ with linear complex 
structures $J_i:=J
\circ u_i$ as it was explained at the beginning 
of Section 2. Observe  that
the maps $u_1$ and $u_2$ are sections of
$E$, and $u_i$ satisfy the equation $\dbar_{J_i} u_i = (\partial_x
+ J_i\partial_y) u_i =0$. 
Without loss of generality we suppose that
$u_1$ has no critical points, possibly except $0$.

\smallskip\noindent
{\slsf Claim 1. \it The image $E_1$ of the differential $du_1: T\Delta
\to E$ is a well-defined $J_1$-complex line subbundle of the complex
bundle $(E, J_1)$ over $\Delta\setminus\{0\}$. It extends to a $J_1$-complex
line subbundle of $E$ of regularity $L^{1,p}$ over $\Delta$ such that
$du_1:T\Delta\to E_1$ is $L^{1,p}$-regular.}

\smallskip %
The claim is trivial in the case when $u_1$ is an immersion
and $\mu=1$. Otherwise we use \eqqref(norm-form-map) and write $u_1$ in the
form
\begin{equation}
\eqqno(u1-mu)
u_1(z) = z^\mu P(z) + z^{2\mu-1} v(z) = z^\mu P(z) + z^\mu \cdot
z^{\mu-1} v(z).
\end{equation}
Notice that now $\mu-1 \geq1$ and hence $z^{\mu-1} v(z)$ is
$L^{2,p}$-regular by  Lemma \ref{l2p-reg}.  It follows that $du_1(z)$ has the
form $z^{\mu-1}H(z)$ for some $L^{1,p}$-regular real bundle
homomorphism $H: T\Delta \to E$ with $H(0) = \mu P(0)$. 
For $z\neq0$ consider the homomorphism $H_1:= du_1\circ (z^{1-\mu}):T_z\Delta\to E$ given by
$w\in T_z\Delta \mapsto du_1(z)(z^{1-\mu}\cdot w) \in E_z$. Observe that in the formulas above the
multiplication of a vector $w \in E_z \cong \cc^n$ with $z$ is 
understood as $(x +
J\st y)\cdot w$.  On the other hand, $du_1$ is $J_1$-linear, 
and consequently
\[
H_1(z) = (x+ J\st y)^{\mu-1}\circ (x+ J_1(z) y)^{1-\mu} \circ  H(z).
\]
The proof of the claim will follow if we shall show that 
$(x + J\st y)^{\mu - 1}(x + J_1(z)y)^{1-\mu}$ is sufficiently close to the identity map. 
For this is it sufficient to show that $(x + J\st y)(x + J_1(z)y)^{-1}$ is sufficiently close to the identity map.
More exactly, that it is $1+ O(|z|)$. Really, if that is proved then for every $k>1$ we shall have
\[
(x + J\st y)^{k}(x + J_1(z)y)^{-k} = (x + J\st y)^{k-1}(1 + O(|z|))(x + J_1(z)y)^{-k+1} = 
\]
\[
(x + J\st y)^{k-1}(x + J_1(z)y)^{-k+1} + (x + J\st y)^{k-1}O(|z|)(x + J_1(z)y)^{-k+1}
\]
and the second term is of order $|z|^{k-1}|z||z|^{-k+1} = O(|z|)$. Therefore the induction will do the job.
Now let us turn to $(x + J\st y)(x + J_1(z)y)^{-1}$. It will be easier to estimate an inverse expression 
$(x + J_1(z)y)(x + J\st y)^{-1}$. Here we obtain
\[\begin{split}
&(x + J_1(z) y)(x + J\st y)\inv =
(x + J_1(z) y)(x - J\st y) \cdot (x^2 +y^2)\inv=\\[3pt]
&(x^2 + y^2 + xy(J_1(z) -J\st ) - y^2(\id+  J_1(z)J\st ) \cdot (x^2
+y^2)\inv.
\end{split}
\]
So we can conclude the pointwise estimate
\[
\norm{ (x + J_1(z) y)(x + J\st y)\inv \; - \; \id } \leq C\cdot
\norm{ J\st -J_1(z)} \leq C' \cdot |z|,
\]
and the claim follows. $L^{1,p}$-regularity of $du_1$ is clear because $u_1$ is $L^{2,p}$-regular
by Proposition \ref{morrey}.

\medskip
Fix an $L^{1,p}$-regular $(J_1,J\st)$-linear trivialization $\Phi: (E, J_1)
\xrar{\;\cong\;} (\Delta\times\cc^n,J\st)$
such that $E_1=du_1(T\Delta)$ is mapped to the
subbundle $\Delta\times\cc^1 \subset \Delta\times\cc^n$ with the fiber consisting of vectors of the form
$(a, 0,\ldots,0)$. Then, denoting by $\cc^{n-1}$ the subspace of $\cc^n$ of vectors
the form $(0, a_2,\ldots,a_n)$, we obtain the bundle $E_2 := \Phi\inv (\cc^{n-1})$
which is a complementary bundle to $E_1$ in $E$, \ie$E=E_1\oplus E_2$. Notice that
$E_1$ and $E_2$ are $J_1$-complex subbundles of $E$.

\medskip

The idea of the proof consists of three steps:

\medskip\noindent
{\sl First,}   to represent $u_2$ in the form
\[
u_2(z) = u_1(\psi(z)) + w(z)
\]
where  $w(z)$ is a $L^{1,p}$-regular section of $E_2$ and  $\psi:
\Delta (r) \to \Delta$ an appropriate $L^{1,p}$-regular
``reparameterization map'' defined locally near the origin.

\smallskip\noindent
{\sl Second,}   to show that $w(z)$ satisfies a differential inequality
of the form \eqqref(3.1).

\smallskip\noindent
{\sl Third,}  to prove that $\psi$ can be chosen to be a holomorphic function.

\medskip
Define an ``exponential'' map $\exp: \Delta \times \Delta \times
\cc^{n-1} \to E$  by
\begin{equation}
\eqqno(exp)
\exp : (z, \zeta, w) \mapsto \big(z, u_1(\zeta) + \Phi\inv(z)w\big).
\end{equation}
The map $\exp$ is well-defined,
$L^{1,p}$-regular in $z$, $L^{2,p}$-regular in $\zeta$ and linear in
$w$. In particular $\exp$ is continuous in $(z,\zeta , w)$.
Moreover, for a fixed $z \neq 0\in \Delta$ the linearization of
$\exp_z:= \exp(z, \cdot,\cdot)$ with respect to variables $\zeta ,
w$ at $\zeta , w=0$ 
is an isomorphism between $T_{\zeta =z}\Delta
\oplus \cc^{n-1}$ and $E_z$. Thus for $z\neq 0$ the map $\exp_z$ is
an $L^{2,p}$-regular diffeomorphism of some neighborhood
$U_z\subset\{ z\} \times \Delta\times\cc^{n-1}$ of the point $(z,0)$ onto
some neighborhood $V_z$ of the point $u_1(z)$ in $E_z =\cc^n$.

\smallskip We need to estimate the size of $V_z$.
In order to do so let us consider the rescaled maps
\[
u_1^t(z) := t^{-\mu}\, u_1(t\cdot z) \text{ \ with } t\in (0,1].
\]
{\slsf Claim 2. \it The family $u_1^t(z)$ is uniformly bounded in
$t\in (0,1)$ with respect to the $L^{2,p}$-norm and the limit map
$\lim_{t\searrow 0} u_1^t(z)$ is $u_1^0(z) := v_0 z^{\mu}$, where $v_0=P(0)$. The limit
is taken in $L^{2,p}$-topology.}

The $\calc^0$-convergence $u_1^t\rightrightarrows u_1^0$ is clear from the representation 
\eqqref(u1-mu). To derive from here the $L^{2,p}$-convergence remark that $u_1^t$ is
$J_t$-holomorphic with respect to the structure
$J_t(w):=J(t^{\mu}\cdot w)$, $w\in B$, and that $J_t$ converge to $J\st$ in
the Lipschitz norm. This implies the $L^{2,p}$-convergence.

\smallskip Further, define the rescaled exponential maps
\begin{equation}
\exp^t_z(\zeta, w) := u_1^t(\zeta ) + \Phi\inv(t\cdot z) w \text{ \
with } t\in [0,1]. 
\eqqno(exp-t)
\end{equation}

\smallskip\noindent
{\slsf Claim 3. \it There exist constants $c^*,c_1, \eps >0$ such that for
every $z\in \{ |z|=\eps \}$ and $t \in [0,1]$ $\exp^t_z(\zeta, w)$ is an $L^{2,p}$-regular
diffeomorphism of $U_z= \{ (\zeta ,w): |\zeta -z |<c_1, |w|<c_1\}$
onto a neighborhood $V_z^t$ of $u_1^t(z)$ in $E_z$ which contains the
ball $\{ |\xi - u_1(z)| < c^*\}$. 

Moreover, the inverse maps
$(\exp^t_z)\inv:V^t_z\to U_z$ are $L^{2,p}$-regular, their $L^{2,p}$-norms 
are bounded by a uniform constant independent of  $z\in \{ |z|=\eps \}$ and $t$,
and the dependence of $(\exp^t_z)\inv$ on $z$ is $L^{1,p}$.}

\smallskip This claim readily follows from the facts that
$K=\{ |z|=\eps \}\times \{ t
\in [0,1]\}$ is a compact, the function $\exp^t_z(\zeta,w)$ is a local
$L^{2,p}$-regular diffeomorphism for every fixed $(z,t)\in K$, and depends
continuously on $t\in[0,1]$ and $L^{1,p}$-regular on $z\in\{|z|=\eps\}$ with respect
to the $L^{2,p}$-topology.

\smallskip Without loss of generality we may assume that $\eps = 1$. This can be always
achieved by an appropriate rescaling.

\smallskip\noindent
{\slsf Claim 4. \it For arbitrary $z\in \Delta\bs\{0\}$ there is a neighborhood $V_z\ni
 u_1(z)$ containing the ball $B(u_1(z), c^*\cdot |z|^{\mu})$ with the constants
 $c^*$ from the Claim 3 such that $\exp_z$ is an $L^{1,p}$-regular
 homeomorphism between some neighborhood $U_z$ of $(z,0)$ in the fiber
 $\{z\}{\times}\cc^n$ and $V_z$. In particular, $c^*$ is independent of $z$.}

Here by  $L^{1,p}$-regular homeomorphism we understand a  homeomorphism 
which is  $L^{1,p}$-regular and its inverse is also  $L^{1,p}$-regular.

\smallskip In order to prove this claim fix some $0<|z|<\frac{1}{2}$ and
set $\tilde z=\frac{z}{|z|}, \tilde \zeta = \frac{\zeta}{|z|},
\tilde w = \frac{w}{|z|^{\mu}}, t=|z|$. Then according to Claim 3 we
have a homeomorphism
\[
\exp^t_{\tilde z} : \{ |\tilde\zeta - \tilde z|<c_1, |\tilde
w|<c_1\} \xrar{\;\cong\;} V_{\tilde z}\supset \{ |\xi
-u_1(\tilde z)|<c^*\} .
\]
But $\exp^t_{\tilde z}(\tilde\zeta ,\tilde w) =
t^{-\mu}u_1(t\tilde\zeta) + \Phi^{-1}(t\tilde z)\tilde w = t^{-\mu}[
u_1(\zeta ) + \Phi^{-1}(z)w] = t^{-\mu}\exp_z(\zeta ,w)$ and this
map is a homeomorphism between $\{
|\frac{\zeta}{|z|}-\frac{z}{|z|}|<c_1, \frac{|w|}{|z|^{\mu}}<c_1\}$
and some $V_z$ containing $\{ |\tilde \xi - u_1(\tilde z)|<c^*\}$.
Therefore $\exp_z$ is a homeomorphism between
\[
 \{ |\zeta - z|<c_1|z|, |\tilde w|<c_1|z|^{\mu}\} \leftrightarrow
V_{z}\supset \{ |\xi -u_1(z)|<c^*\cdot|z|^{\mu}\}.
\]

\medskip\noindent{\slsf Claim 5. \it For $z$ sufficiently small,
 $u_2(z)=u_1(\psi(z)) +w(z)$ for some $L^{1,p}$-regular function $\psi(z)$ in $\Delta$
 and some $w\in L^{1,p}(\Delta,E_2)$. 
}

\smallskip Since $u_2(z)-u_1(z)=O(|z|^{\mu +\alpha})$, for $z$ small enough we obtain
$u_2(z)\in B(u_1(z),c^*\,|z|^{\mu})$. Define  $(\zeta(z),W(z)):=\exp_z^{-1}(u_2(z))$
where $\exp_z^{-1}:V_z\to U_z$ is the local inversion of the map $\exp_z$ which
exists by {\sl Claim 4}. Set $\psi (z):=\zeta (z)$, $w(z):=\Phi^{-1}(z)W(z)$.  We obtain the desired
relation
\begin{equation}
u_2(z) = u_1(\psi(z)) + w(z) \eqqno(3.4),
\end{equation}
which holds in some small punctured disc $\Delta_r\bs\{0\}$.  Making an appropriate
rescaling we may assume that \eqqref(3.4) holds in the whole punctured disc
$\Delta\bs\{0\}$.  Moreover, $\psi(z)$ and $w(z)$ are $L^{1,p}_\loc$-regular in
$\Delta\bs\{0\}$.

To estimate the norm $\norm{\psi}_{L^{1,p}(\Delta)}$ we define the rescalings 
$u_2^t(z):= t^{-\mu}{\cdot}u_2(tz)$, $\psi^t(z):= t\inv\psi(tz)$ and
$w^t(z):=t^{-\mu}w(tz)$. Then we obtain $u_2^t(z) = u_1^t(\psi^t(z)) + w^t(z)$
which is the rescaled version of \eqqref(3.4). Consequently, these function
satisfy the rescaled relation
\[
(\psi^t(z), \Phi(tz)w^t(z))=(\exp_z^t)^{-1}(u_2^t(z)).
\]
By {\sl Claim 4} the family of maps $(\exp_z^t)^{-1}$ is continuous in $t$ and
$L^{1,p}$-regular in $z$ with respect to $L^{2,p}$-topology. {\sl Claim 2}
applied to $u_2(z)$ gives us the uniform $L^{2,p}$-boundedness of the  family $u_2^t(z)$
in $t$.
As a consequence, we conclude that the functions $\psi^t(z)$ satisfy the uniform
estimate 
\[
\norm{\psi^t}_{L^{1,p}(\Delta\bs\Delta(\half))} \leq C
\]
with the constant $C$ independent of $t$. Making the reverse rescaling we
conclude the estimate
\[
\norm{d\psi}_{L^p(\Delta(2^{-k-1})\bs\Delta(2^{-k}))} \leq C{\cdot}2^{-2k/p}.
\]
Now the summation over the annuli $\Delta(2^{-k-1})\bs\Delta(2^{-k})$ gives us the desired estimate 
$\norm{d\psi}_{L^p(\Delta)} \leq \left(\frac{4}{3}\right)^{1/p}{\cdot}C$. The
$L^{1,p}_\loc$-regularity of $w(z)$ in $\Delta$ follows from the relation \eqqref(3.4).
The Claim is proved.

\medskip 
Consider the pulled-back bundles $E':=\psi^*E$, $E'_1 := \psi^*E_1$
and $E'_2:= \psi^*E_2$ over $\Delta$. Equip $E'$ with the complex
structure $J_1':=\psi^*J_1 = J(u_1\circ\psi )$. Let $\pr'_2$ be the
projection of $E'$ onto $E'_2$ parallel to $E'_1$. Consider  the
following expression
\begin{equation}\eqqno(pr-strix)
\pr'_2\big( ( \partial_x + J_1' \cdot \partial_y) w(z) \big) =
\pr'_2\big( ( \partial_x + J_1' \cdot
\partial_y) (u_2(z) -u_1(\psi(z))\big).
\end{equation}
Let us treat the terms in \eqqref(pr-strix) separately. The second term
on the right hand side
\[
( \partial_x + J_1' \cdot \partial_y)u_1(\psi(z))
\]
is the Cauchy-Riemann operator $\dbar_{J_1'}$ applied to the
composition $u_1':=u_1 \circ \psi$.

\medskip\noindent{\slsf Claim 6.} {\it Let $\psi :\Delta\to \Delta$ be a
$L^{1,p}$-map and let $u:\Delta\to \cc^n$ be a $J$-holomorphic
curve. Then
\begin{equation}
\dbar_J(u\circ\psi) = du\circ \dbar\psi, \eqqno(3.5)
\end{equation}
where $\dbar\psi$ is the standard $\dbar$-derivative of the function
$\psi$. }

\smallskip The expression
$\dbar_{J}(u \circ \psi)$ computes the $J$-antilinear component of
the differential $d(u \circ \psi) = du \circ d\psi$. Since $du$ is
$J$-linear,  the antilinear part of $du\circ d\psi$ will be $du$ of
the antilinear part of $d\psi$ which is $\dbar\psi$. Therefore we
conclude the relation \eqqref(3.5). The claim is proved.

\smallskip In our case this gives

\begin{equation}
\dbar_{J_1'}(u_i\circ\psi ) = du_1\circ \dbar\psi. \eqqno(3.6)
\end{equation}

\smallskip
Further, observe that $du_1 \circ \dbar\psi$ takes values in the
pulled-back $E'_1 = \psi^*E_1$. So\break $\pr'_2( \dbar_{J_1'}(u_1 \circ
\psi))$ vanishes identically.

The next term to estimate is $( \partial_x + J_1' \cdot
\partial_y)u_2$. Subtracting the 
equation $0 =\dbar_J u_2 =
(\partial_x + J\circ u_2 \cdot \partial_y) u_2$ we obtain
\[
(J_1' - J\circ u_2)\d_yu_2 = (J \circ u_1 \circ \psi -  J\circ u_2)
\cdot
\partial_yu_2.
\]
The $L^{2,p}$-regularity for $u_2$ provides the $L^{1,p}$-regularity
of $\partial_yu_2$, whereas the Lipschitz condition on $J$ yields
the pointwise estimate
\begin{equation}
\eqqno(pr-chast)
\big| J \circ u_1 \circ \psi(z) -  J\circ u_2(z) \big|\leq
Lip(J) \cdot |w(z)|.
\end{equation}
Therefore the right hand side of \eqqref(pr-strix) is estimated by $h\cdot|w|$
with some $h\in L^p(\Delta)$.

\smallskip %
Now, let us rewrite the left hand side $\pr'_2 \big ( ( \partial_x + J \circ u_1
\circ \psi \cdot \partial_y) w\big)$ 
of \eqqref(pr-strix) as a $\dbar$-type operator of
$w$. Consider the restriction $\pr'_2 : E_2 \to E'_2$ of the
projection $\pr'_2$ onto $E_2$. Using the facts that $\psi(z)$ is
continuous and $\psi(0) =0$, we conclude that $\pr'_2(z) : (E_2)_z
\to (E'_2)_z = (E_2)_{\psi(z)}$ 
is a bundle isomorphism over a
sufficiently small disc $\Delta_r$.  So setting $\ti w(z) := \pr'_2
w(z)$ 
we obtain a pointwise estimate 
\begin{equation}
\eqqno(w-til-w)
1/C\cdot |w(z)| \leq |\ti w(z)| \leq C\cdot |w(z)|
\end{equation}
in the disc $\Delta_r \ni z$ with uniform constant $C$.

Similar to $\pr'_2$ define the projection $\pr'_1$ from $E$ onto
$E'_1 = \psi^*E_1$.  Denote $\nabla_x :=  \pr'_2 \circ
\partial_x \circ \pr'_2$
and $\nabla_y:=  \pr'_2 \circ \partial_y \circ \pr'_2$.  Using this
we obtain
\begin{equation}
\eqqno(nabla-h)
\pr'_2 \big ( ( \partial_x + J_1' \cdot
\partial_y) w\big) = ( \nabla_x + J_1' \cdot
\nabla_y) (\pr'_2w) + H( w)
\end{equation}
with some $L^p$-regular endomorphism $H$. Summing up, we conclude a
pointwise differential inequality
\begin{equation}
\eqqno(debar-nab)
\big|( \nabla_x + J_1' \cdot \nabla_y) (\tilde w)\big| \leq h|\tilde w|
\end{equation}
with $J'_1 :=  J \circ u_1 \circ \psi$ and an  $L^p$-regular function $h$. 
If we fix some $L^{1,p}$-trivialization $e_1,...,e_{n-1}$ of 
$E_2^{'}$ and remark that in (any) such trivialization $\nabla_x = \d_x + R_x$
and $\nabla_y = \d_y + R_y$ with some $R_x, R_y\in L^p(\Delta , \endo (E_2^{'}))$.
This gives us the following estimate
\begin{equation}
\big | (\partial_x + J'_1 \partial_y)\ti w \big| \leq h\cdot |\ti w|
\eqqno(3.7)
\end{equation}
Observe that $\ti w(z)$ can not vanish identically
since otherwise the image $u_2(\Delta_r)$ would lie in
$u_1(\Delta)$.

\smallskip Now we can apply  Lemma 3.1 and conclude that $\tilde w(z)$
either vanishes identically or $\tilde w(z) = z^{\nu} f(z)$ for
some $ f(x)\in L^{1,p}(\Delta_r, \cc^{n-1})$ with $f(0) \neq 0$.
The integer $\nu$ must be bigger than $\mu$, because
$u_2(z)-u_1(z)=o(|z|^{\mu +\alpha})$.  Since the projection $w(z)
\mapsto \ti w(z) := \pr'_2(w(z))$ is an $L^{1,p}$-regular
isomorphism, we obtain the same structure for $w$. Finally, observe
that $f(0)$ lies in the fiber $(E_2)_0$ which is
$J(0)=J\st$-transverse to $(E_1)_0 = \cc v_0$.
Therefore we obtain

\begin{equation}
u_2(z) = u_1(\psi (z)) + z^{\nu} w(z), \eqqno(3.8)
\end{equation}
where $\nu >\mu$ and $w(0)$ linearly independent of $v_0$.

\smallskip
\smallskip\noindent{\slsf Claim 7.} {\it There exists a \emph{holomorphic} $\psi $
 satisfying \eqqref(3.8). }

\medskip


Assume that we have $w(z)\equiv0$ and therefore $u_2(z) =
u_1(\psi(z))$.  It follows from \eqqref(3.6) that
\[
du_1\circ \dbar\psi = \dbar_{J_1'}(u_1\circ \psi) =
\dbar_J(u_1\circ\psi) = \dbar_Ju_2=0
\]
and therefore that $\dbar \psi(z)\equiv0$. So $\psi(z)$ is
holomorphic.

\smallskip %
Assume now $w(z)$ is not identically $0$.  In this case we are going to
construct recursively a sequence of complex polynomials  $\phi_i(z)$  and an
increasing sequence $\mu< \nu_1 < \nu_2 <\ldots <\nu_l$ of integers with the
following properties: 

$\bullet$ \ $\phi_i(z)= z + O(z^2)$ 

$\bullet$ \ $u_2(z) = u_1(\phi_i(z)) + z^{\nu_i} v_i(z)$ with some
$v_i(z)\in L^{1,p}(\Delta,\cc^n)$ such
that for $j<l$ the vectors $v_j(0)$ are proportional to $v_0$.

Lemma 3.2 insures the existence of the desired $\nu_1 >\mu$ and $v_(z)$ with
$\phi_1(z) \equiv z$. Assume that we have constructed such sequences
$\mu<\nu_0<\nu_1<\nu_2<\ldots\nu_k$ and $v_1(z), \ldots, v_k(z)$, and that $v_1(0), \ldots, v_k(0)$ are
proportional to $v_0$. Observe that for any integer $m\geq2$ and any $a\in\cc$ we have
\[\begin{split}
u_1(\phi_k(z) + a z^m) &
= u_1(\phi_k(z))  + du(\phi_k(z))\circ d\phi_k ( a z^m)  +
O(z^{m+\mu}) \\[3pt]
& = u_1(\phi_k(z)) + \mu \phi_k(z)^{\mu-1}\cdot \phi_k'(z)\cdot a z^m \cdot  v_0 +
O(z^{m+\mu})\\[3pt]
& = u_1(\phi_k(z)) + \mu z^{m+\mu-1} \cdot a\cdot  v_0 +
O(z^{m+\mu}).
\end{split}
\]
Set $m_k := \nu_k -\mu +1$, defined $a$ from the relation
\[
\mu \cdot a\cdot  v_0 + w_k(0) =0
\]
and put
\[
\phi_{k+1}(z) := \phi_k(z) + a z^{m_k}.
\]
Then $u_2(z) - u_1(\phi_{k+1}(z)) = O( |z|^{m+\mu+\alpha})$.  Now we can apply Lemma 3.2
to $u_2(z)$ and to $u_1(\phi_{k+1}(z))$ and obtain a new $\nu_{k+1} > m +\mu \geq \nu_k$
and a new $v_{k+1}(z)$. 

Compare the obtained presentations $u_2(z) = u_1(\phi_i(z)) + z^{\nu_i} v_i(z)$
with the decomposition \eqqref(3.8). Notice that for a fixed bundle $E_2$ the
decomposition \eqqref(3.8) is unique.  This implies that at some step we
obtain $\nu_l = \nu$ and $v_l(0) = w(0)$ with $\nu$ and $w(z)$ from Comparison
Theorem. At this step $v_l(0) = w(0)$ is not proportional to $v_0$ and the
recursive procedure halts.

\medskip All what is left to prove is \eqqref(log-est). In Section 6
we shall prove that $J$-holomorphic mappings in Lipschitz-continuous $J$
are $\calc^{1, LnLip}$ and therefore the subbundle $E_1=du_1(T\Delta)$
is a $\calc^{LnLip}$-regular. This implies the same regularity of the
projection $\pr_{v_0}$. Since $\pr_{v_0}w(0)=0$ this gives
\eqqref(log-est).

\smallskip
This finishes the proof of the part (a) of the Comparison Theorem. 

\begin{rema} \rm
\label{rem-transvers}
As we claimed in the Introduction the vector $w(0)$ can be chosen in any 
given complex hyperplane transversal to $v_0$. Really, 
if $E_2(0)$ is such a plane, then we chose $\phi$ (after  the end of the proof 
of Step 1) in such a way that $E_2(0) = \Phi^{-1}(\cc^{n-1})$. Then in the 
remaining part of the proof of the Part (a) of the Comparison Theorem we 
established that the vector function $w$ in question takes its values in 
the bundle $E$. Therefore, in particular, $w(0)\in E_2(0)$.
\end{rema}

\newprg[prgLOC.compb]{Proof of the part {\bf (b)} of the Comparison Theorem.}

We continue the proof of the Comparison Theorem.

\medskip\noindent{\slsf Claim 8.} This claim we shall state in the form of a lemma.

\begin{lem}
\label{funk-eq}
Let $d>1$ be an integer and $\eta $ be a primitive root of unity of degree $d$,
and $\psi $ a holomorphic function in the unit disc $\Delta$ of the form
$\psi(z)=z+O(z^2)$.  Then there exists a holomorphic function $\varphi $ of the form
$\varphi(z)=z+O(z^2)$ defined in some smaller disc $\Delta_r$ such that %
$\eta\varphi(z)-\psi(\varphi(\eta z))=z^{d+1}\cdot\gamma(z^d)$ with some function $\gamma$
holomorphic in $\Delta_r$.
\end{lem}
\proof Roughly speaking, the lemma claims that making an appropriate
reparameterization one can eliminate all the terms of the Taylor expansion of
$\psi(z)$ except those of degrees  $kd+1$.

\smallskip We want to apply the implicit function theorem.
For this purpose we
need to fix certain smoothness class of holomorphic functions, the concrete
choice of such a space plays no role in the proof. Denote by
$\calh$ the space of holomorphic functions in $\Delta$ which are $\calc^{\alpha}$
smooth up to boundary.

Replacing the given function $\psi(z)$ by its appropriate rescaling
$\psi^{(t)}(z):=t\inv\psi(tz)$ we may assume that the norm
$\norm{\psi(z)-z}_{\calc^1(\Delta)}$ is small enough.

For $l=0,\ldots,d-1$, denote by $\calh_l$ the subspace of $\calh$ consisting of
functions $\varphi(z)$ of the form $\varphi(z)=z^l\varphi_1(z^d)$. In other words, the Taylor
series of $\varphi(z)\in\calh_l$ contains only monomials of degree $m\equiv l(\mod d)$. The
space $\calh_l$ is the kernel of the operator $\varphi(z)\mapsto\varphi(\eta z)-\eta^l\varphi(z)$. Clearly,
we obtain the decomposition $\calh=\oplus_{l=0}^{d-1}\calh_l$. Denote by
$\calh_1^\perp$ the complement to $\calh_1$ in this sum, \ie
$\calh_1^\perp=\oplus_{l\neq1}\calh_l$, and by $\pi_1^\perp$ the projection on this space
parallel to $\calh_1$. Finally, let $\calv$ be the Banach subspace of
$\calh_1^\perp$ consisting $\varphi(z)\in\calh_1^\perp$ satisfying $\varphi(z)=O(z^2)$ and
$(z+\calv)$ the shift of $\calv$ in $\calh$ by the function $z$. Thus $\varphi(z)$
lies in $(z+\calv)$ if and only if  $\varphi(z)=z+\varphi_1(z)$ with $\varphi_1(z)\in\calv$.

Now consider the map $F_\psi:(z+\calv)\to\calh_1^\perp$ given by
\[
F_\psi : \varphi(z)\in(z+\calv)\; \longmapsto \:\pi_1^\perp\big(\eta\varphi(z)-\psi(\varphi(\eta z))\big),
\]
in which $\psi$ is considered as a parameter, varying in the space of holomorphic
function defined in some larger disc $\Delta_{1+\varepsilon}$, so that $\psi\in\calh(\Delta_{1+\varepsilon})$.  Then
$F_\psi(\varphi)$ takes values in $\calv$, is holomorphic in $\varphi\in(z+\calv)$, continuous
(in fact, also holomorphic) in $\psi(z)$ and its linearization in $\varphi$ at point
$(\psi_0(z)\equiv z,\varphi = 0)$ is
\[
DF_{\psi_0,0}:\dot\varphi \mapsto\Big(\eta\dot\varphi(z)- \dot\varphi(\eta z)\Big).
\]
Then $DF_{\psi_0}$ is an isomorphism on $\calv$ since its restriction on each
$\calh_l$ is the multiplication with the non-zero scalar $\eta-\eta^l$.

Now the implicit function theorem applies and for every $\psi$ close to $\psi_0(z)\equiv z$
gives a function $\varphi \in z + \calv$ (\ie of the form $\varphi (z) = z + O(z^2)$)
such that $\eta \varphi (z) - \psi (\varphi (\eta z)) \in \calh_1$, \ie is
of the form $z\gamma (z^d)$ as required.

\smallskip\qed

\smallskip\noindent{\slsf Claim 9.} We shall prove the following:

\begin{prop}
 \label{compar1} Let $J$ be a Lipschitz-continuous almost complex structure in
 the unit ball $B\subset\cc^n$ with $J(0)=J\st$, $u:\Delta\to B$ a $J$-holomorphic map such
 that $u(z)=v_0z^\mu+ O(|z|^{\mu+\alpha})$ with $v_0\neq0\in\cc^n$, $d\neq1 $ a divisor of
 $\mu$, and $\eta=e^{2\pi\isl/d}$ the primitive root of unity of degree $d$.  Let
 $u(\eta z)=u(\psi(z))+z^{\nu }w(z)$ be the presentation provided by the part (a) of the Comparison
 Theorem.

Then there exists a holomorphic reparameterization $\phi(z)$ of the form
 $\phi(z)=z+O(z^2)$ such that
\begin{itemize}
\item  $u(\phi(\eta z))\equiv u(\phi(z))$ in the case when $w(z)\equiv0$,
\item $u(\phi(\eta z))-u(\phi(z))=w(0) z^\nu+O(|z|^{\nu +\alpha})$ otherwise. Moreover, in this
 case $\nu$ is not a multiple of $d$.
\end{itemize}
\end{prop}

\proof Let $\varphi(z)=z+O(z^2)$ be the function constructed in Lemma \ref{funk-eq},
so that $\eta\varphi(z)=\psi(\varphi(\eta z))) + z^{d+1}\gamma(z^d)$. Substitute the latter relation into the
comparison relation $u(\eta z)=u(\psi(z))+z^\nu w(z)$ and obtain
\begin{equation}
\eqqno(u-psi-phi)
u\big(\psi(\varphi(\eta z)) +z^{d+1}\gamma(z^d) \big)= u(\eta\varphi(z))=
u(\psi(\varphi(z)) + \varphi^\nu(z)\cdot w(\varphi(z)).
\end{equation}
Denote $u(\psi(\varphi(z)))$ by $\ti u(z)$, this is a reparameterization of the map
$u(z)$ in the new coordinate, such that the old one is given by the formula
$\varphi^{-1}(\psi^{-1}(z))$. (Notice that we use the same notation $z$ for both.)

We want to rewrite the \eqqref(u-psi-phi) in this new coordinate.  Let us start from the left hand side.
Assume that
$\gamma(z)$ is not identically $0$ and denote by $k$ the order of vanishing of
$\gamma(z)$ at $z=0$. Then $z^{d+1}\gamma(z^d)=az^{(k+1)d+1}+O(z^{(k+1)d+2})$. Since $\varphi(z)$ and
$\psi(z)$ are reparameterization of the form $z+O(z^2)$, one can rewrite
$\psi(\varphi(\eta z)) +z^{d+1}\gamma(z^d)$ in the form $\psi\big[\varphi\big(\eta z + z^{(k+1)d+1}\ti
\gamma(z)\big)\big]$ %
with holomorphic function $\ti\gamma(z)$ satisfying $\ti\gamma(0)=a$. In the other case
$\gamma(z)\equiv0$ we obtain a similar relation with $\ti\gamma(z)\equiv0$.

\smallskip As for the right hand side  one can rewrite the expression $\varphi^\nu(z)\cdot w(\varphi(z))$ in the form
\[
(\psi^{-1}(z)))^\nu\cdot \ti w(\psi^{-1}(z)))
\]
with a new function $\ti w(z)$ of the same regularity $L^{1,p}$ such that %
$\ti w(0)=w(0)$. So we conclude that the reparameterized map %
$\ti u(z)=u(\psi(\varphi(z)))$ satisfies the relation
\[
\ti u(\eta z+z^{(k+1)d+1}\ti \gamma(z))=\ti u(z)+ z^\nu\cdot \ti w(z)
\]
Finally, using $\ti u(z)=v_0z^\mu+O(|z|^{\mu+\alpha})$, we obtain
$\ti u(\eta z)=\ti u(z)+ z^\nu\cdot \ti w(z)-\mu\eta^{\mu-1}\ti \gamma(0)v_0z^{kd+\mu} +O(|z|^{kd+\mu+\alpha})$.

\smallskip Now assume that $kd+\mu\leq\nu$ or that $\ti w(z)\equiv0$. In this case
$\ti u(\eta z)=\ti u(z)+ w'\cdot z^{kd+\mu} +O(|z|^{kd+\mu+\alpha})$ with some non-zero vector
$w'$. Then, using the equality $\eta^\mu=1$ (since $d$ is a divisor of $\mu$)
\begin{equation}\eqqno(sum-eta-k)
 \textstyle
 0=\sum_{j=1}^d\big(\ti u(\eta^jz)-\ti u(\eta^{j-1}z) \big)=
d\cdot w'\cdot z^{kd+\mu} +O(|z|^{kd+\mu+\alpha})
\end{equation}
which is a contradiction.

Observe that we obtain the same contradiction in the case when $kd+\mu>\nu$
(including the extremal case $\gamma(z)\equiv0$) and $\nu$ is a multiple of $d$, so that
$\eta^\nu=1$.

From this contradiction we can conclude the following:
\begin{itemize}
\item In the case $w(z)\equiv0$ we must have $\gamma(z)\equiv0$, and so the function
 $\phi(z):=\psi(\varphi(z))$ is the desired reparameterization.
\item In the case $w(0)\neq0$ we must have $\nu<kd+\mu$ and $d$ can not be a divisor
 of $\nu$. Again, $\phi(z):=\psi(\varphi(z))$ is the desired reparameterization.
\end{itemize}

\smallskip\qed

\smallskip Comparison Theorem is proved.

\newsect[sectPP]{Primitivity and Positivity of Intersections}

In this section we shall prove the important regularity properties of
$J$-complex curves with  Lipschitz-continuous $J$-s, \ie Theorems A and B
from the Introduction.

\newprg[prgPP.def]{Definitions}

We fix an almost complex manifold $(X,J)$ with $J\in \calc^{Lip}$. Let
$(S_i,j_i)$, $i=1,2$  be two Riemann surfaces with complex structures $j_1$ and $j_2$.

\begin{defi}
\label{dist-map}
Two $J$-holomorphic maps $u_1: (S_1,j_1) \to X$ and $u_2: (S_2,j_2) \to X$ with
$u_1(a_1)=u_2(a_2)$ for some $a_i\in S_i$ are called {\sl distinct at
 $(a_1,a_2)$} if there are no neighborhoods $U_i\subset S_i$ of $a_i$ with
$u_1(U_1)=u_2(U_2)$. We call $u_i: (S_i,j_i) \to X$ {\sl distinct} if they are
distinct at all pairs $(a_1,a_2)\in S_1 \times S_2$ with $u_1(a_1)=u_2(a_2)$.
\end{defi}

A related notion is the following

\begin{defi}
\label{prim-map}
A $J$-holomorphic map $u: (S,j) \to X$ is called {\sl primitive} if there are no
disjoint non-empty open sets $U_1, U_2\subset S$ with $u(U_1) =u(U_2)$.
\end{defi}

Note that a primitive $u$ must be non-constant. Let $B$ be the unit ball in $\cc^2$ and 
$u_1,u_2:\Delta\to B$ two $\calc^1$-regular maps with the following properties: both images
$\gamma_1:=u_1(\partial\Delta)$ and $\gamma_2:=u_2(\partial\Delta)$ of the boundary boundary
circle $\partial\Delta$ are immersed real curves lying on boundary sphere $\partial B=S^3$,
the origin $0$ is the only intersection point of the images $M_1:=u_1(\Delta)$ and
$M_2:=u_2(\Delta)$. Let $\ti M_i$ be small perturbations of $M_i$ making them intersect transversally
at all their common points.

\begin{defi}
\label{int-number}
The {\slsf intersection number} of $M_1$ and $M_2$ at zero is defined to be the algebraic
intersection number of $\ti M_1$ and $\ti M_2$. It will be
denoted by $\delta_0(M_1,M_2)$ or, $\delta_0$ if $M_1$ and $M_2$ are clear from the
context.
\end{defi}

This number is independent of the particular choice of perturbations
$\ti M_i$.  We shall use the fact that the intersection number of $M_1$ and
$M_2$ is equal to the {\sl linking number} $l(\gamma_1, \gamma_2)$ of the curves $\gamma_i$
on $S^3$, see \eg \cite{Rf}. $M_1$ and $M_2$ intersect
transversally at zero if the tangent spaces $T_0M_1$ and $T_0M_2$ are
transverse. In this case $\delta_0(M_1,M_2)=\pm 1$.

\newprg[prgPP.proof]{Proof of Theorems A and B}

We turn now to the proof of Theorems A and B. It is divided into
several steps, some of them will be stated as lemmas for the convenience of the
future references. Let  $M_i$ are $J$-complex discs in $(\cc^2,J)$, $i=1,2$. 
By  \eqqref(norm-form-map) we have the following presentations

\begin{equation}
u_1(z) = z^{\mu_1} v_1(0) + O(|z|^{\mu_1 + \alpha})
\eqqno(4.1)
\end{equation}
\[
u_2(z) = z^{\mu_2} v_2(0) + O(|z|^{\mu_2 + \alpha})
\]
with non-zero vectors $v_1, v_2 \in T_0 B =\cc^2$ and integers
$\mu_i >0$ and with some $\alpha>0$. Moreover, by
\eqqref(norm-form-diff) for both curves we have
\begin{equation}
du_i(z) = \mu_iz^{\mu_i-1}v_i(0) + O(|z|^{\mu_i-1+\alpha}).
\eqqno(4.2)
\end{equation}
\eqqref(4.1) and \eqqref(4.2) imply the transversality of small $J$-complex discs
$u_i(\Delta (\rho))$ to small spheres $S^3_r$. More precisely, there
exist radii $\rho>0$ and $R>0$ such that for any $0<r<R$ the
$J$-curves $u_i(\Delta (\rho))$ intersect the sphere $S^3_r\deff \{
|w_1|^2+ |w_2|^2=r^2 \}$ transversely along smooth immersed circles
$\gamma_i(r)$. In fact, the asymptotic relation (4.1) provides that
for any $\theta \in [0,2\pi]$ there exists at least one solution of
the equation $\bigl|u_i(\rho_{i} e^{\isl\theta}) \bigr| =r$ with
$\rho_i< \rho$, and that for any such solution $\rho_i$ the quotient
$\rho_i/\left( \frac{r}{|v_1(0)|}\right)^{1/\mu_i}$ must be close to one. Then one uses (4.2) to
show that the set $\ti\gamma_i(r) \deff \{ z\in\Delta: \bigl|u_i(z)
\bigr| =r \}$ is, in fact, a smooth immersed curve in $\Delta$,
parameterized by $\theta\in [0,2\pi]$, and that $u_i: \ti\gamma_i(r)
\to S^3_r$ is an immersion with the image $\gamma_i(r)$.

Taking an appropriate small subdisc and rescaling, we may assume
that $\rho=1 =R$. Note that the points of the self- (resp.\ mutual)
intersection of $\gamma_i(r)$ are self- (or resp.\ mutual)
intersection points of $u_i (\Delta)$. Let us call $r \in ]0,1[$
non-exceptional if curves $\gamma_i(r) \subset S^3_r$ are imbedded
and disjoint. Thus $r^* \in ]0,1[$ is exceptional if $S^3_{r^*}$
contains intersection points of $u_i (\Delta)$.

Lemma \ref{discrt-two} and Corollary \ref{discrt-one} provide that any such
intersection point is isolated in the punctured ball $\check B \deff
\{ 0< |w_1|^2 + |w_2|^2 <1 \}$. This implies that either there exist
finitely many exceptional radii $r^*\in ]0,1[$, or that they form a
sequence $r^*_n$ converging to $0$.

Denote $M_i(r) \deff u_i (\Delta) \cap B(r)$. For non-exceptional
$r$ we can correctly define the intersection index of $M_1(r)$ with
$M_2(r)$ as the linking number of $\gamma_1(r)$ and $\gamma_2(r)$.

\medskip\noindent {\slsf Step 1.} In this step we shall prove that two distinct
$J$-complex curves intersect by a discrete set.
\begin{lem}
\label{discrt-two}
Let $J$ be a Lipschitz-continuous almost complex structure on a manifold $X$ and
let $u_1: (S_1, j_1) \to X$ and $u_2: (S_2, j_2) \to X$ be two distinct non-constant
$J$-holomorphic maps. Then:

\sli  The set $\{(z_1,z_2)\in S_1\times S_2:
u_1(z_1) = u_2(z_2)\}$ is discrete in $S_1\times S_2$.

\slii The intersection index at every $p=u_1(z_1)=u_2(z_2)$ is at least $\mu_1\cdot \mu_2$,
where $\mu_i$ is the multiplicity of zero of $u_i$, $i=1,2$.
\end{lem}
\proof The claim is local so we may assume that $(S_1, j_1) = (S_2, j_2) = (\Delta,
J\st)$, $X$ is
the unit ball $B$ in $\cc^2$, $J(0) = J\st$, and $u_1(0) =
u_2(0) =0 \in B$.

\smallskip Write each map in the form $u_i(z)=v_iz^{\mu_i}+O(|z|^{\mu_i+\alpha})$. We
must consider three cases.

\smallskip\noindent
{\sl Case 1.} {\it The~vectors $v_1(0)$ and $v_2(0)$ are not
collinear.}

It is easy to see that, in this case, $0\in\cc^2$ is an~isolated
intersection point of $u_1(\Delta)$ and $u_2(\Delta)$ with
multiplicity exactly $\mu_1\cdot \mu_2$. In particular, intersection index
in every such point is positive. In fact, consider the dilatations:
$J_t(w) = J(t^{\mu_1\mu_2}w)$, $u_1^{t}(z) = t^{-\mu_1\mu_2}u_1(t^{\mu_2}z)$ and
$u_2^{t}(z) = t^{-\mu_1\mu_2}u_2(t^{\mu_1}z)$ for a small $t>0$. $u_i^t$ are
$J_t$-holomorphic and converge to $\mu_i$-times taken disc in the direction
of $v_i(0)$. The rest is obvious.

\smallskip\noindent
{\sl Case 2.} {\it The~vectors $v_1(0)$ and $v_2(0)$ are collinear and $\mu_1=\mu_2=1$.}
In other words,  $u_1(\Delta)$ and $u_2(\Delta)$ are non-singular tangent discs.

Rescaling parameterization of $u_i$ and rotating coordinates in
$\cc^2$ we can suppose that $v_1(0)=v_2(0)=e_1$. Applying the Comparison Theorem we see that
\begin{equation}
\ti u_2(z) - \ti u_1(\psi (z)) = z^{\nu}w(z)
\end{equation}
where $w(0)=e_2$ and $\psi(z)=z+O(z^2)$ is some holomorphic
reparameterization. Considering intersections we see that for $r>0$ small
enough the circles $\gamma_1(r):=u_1(\Delta)\cap S^3_r$ and $\gamma_2(r):=u_2(\Delta)\cap S^3_r$ are
imbedded and, as we go along $\gamma_1(r)$, $\gamma_2(r)$ stays in the tubular
neighborhood of $\gamma_1(r)$ of radius $\rho=2r^\nu$ and winds $\nu$ times around
$\gamma_1(r)$. This shows that the linking number $l(\gamma_1(r),\gamma_2(r))$ is $\nu$.

\smallskip\noindent {\sl Case 3.} {\it The~vectors $v_1(0)$ and
 $v_2(0)$ are collinear and $\mu_1,\mu_2$ are arbitrary.}

The discs $u_1(\Delta)$ and $u_2(\Delta)$ are immersed outside the origin $0$. The
consideration from {\sl Case 2} show that the intersection points of $u_1(\Delta)$
and $u_2(\Delta)$ are discrete in $u_1(\Delta)\bs\{0\}$. In particular, for any
sufficiently small $r>0$ there are finitely many intersection points in the
spherical layer $B_{2r}\bs B_r$. In particular, there exists a sufficiently
small $r>0$  such that the circles $\gamma_1(r):=u_1(\Delta)\cap S^3_r$ and
$\gamma_2(r):=u_2(\Delta)\cap S^3_r$ are immersed and disjoint.

In Theorem \ref{cusp-pert} below we show that for a given $r>0$ small enough
there exists a $J$-holomorphic perturbation $\ti u_2(z)$of the map $u_2(z)$
such that $\ti u_2(z)=\ti v_iz^{\mu_i}+O(|z|^{\mu_i+\alpha})$ with $\ti v_2$ different
from but arbitrarily close to $v_2=v_1$. Moreover, the map $\ti u_2(z)$ is
arbitrarily close to $u_2(z)$. In particular, the circle %
$\ti \gamma_2(r):=\ti u_2(\Delta)\cap S^3_r$ remains disjoint from $\gamma_1(r)$ and homotopic
to $\gamma_2(r)$ in $S^3_r\bs\gamma_1(r)$, the linking number $l(\gamma_1(r),\ti
\gamma_2(r))$ %
remains equal $lk(\gamma_1(r),\gamma_2(r))$, and $\ti u_2(\Delta)\cap B_r$ remains immersed
outside the origin. Now using first two cases we conclude that there are
finitely many intersection points of $u_1(\Delta)$ and $\ti u_2(\Delta)$ in $B_r$,
the intersection index in $0$ is $\mu_1\cdot\mu_2$ and that all other
intersection indices are positive. Since  $l(\gamma_1(r),\ti
\gamma_2(r))$ is the sum of these indices we conclude the part \slii of the lemma.

\qed

\begin{rema}\rm
 Let us point out that the statements (i) and (ii) of Theorem B are proved.
 The proof of (iii) is now obvious. Really, if $\delta_p=1$ then $\mu _1 = \mu _2 =1$,
 \ie $u_i$ are not singular. If they are tangent then $\nu >1$, but we proved
 that $\delta_p =\nu$, contradiction. This finishes the proof of Theorem B.
\end{rema}
We continue the proof, now of Theorem A and, therefore, turn our attention to
a single $J$-holomorphic mapping $u:S\to X$. The following step in the case of
multiplicity $\mu =1$ is trivial and therefore we suppose that $\mu \geq 2$.
We will also use \eqqref(4.1) and \eqqref(4.2) as holding true for $\cc^n$-valued
maps (which is, of course, so in \eqqref(norm-form-map) and \eqqref(norm-form-diff)).

\smallskip\noindent{\slsf Step 2.} Multiple $J$-holomorphic mapping $u$ with
multiplicity of zero equal to $\mu$ can be locally
represented in the form $u(z) = \tilde u(z^d)$ with some $J$-holomorphic $\tilde u$
and some integer $d$.
\begin{lem}
\label{factoriz}
Let $u:S\to X$ be a $J$-holomorphic map with $J\in\calc^{Lip}$ and let
$p \in S$ be a critical point of $u$ of multiplicity $\mu\geq 2$. Then there exist
a neighborhood $W\subset S$ of $p$, a holomorphic map $\pi:W \to \Delta$, and a
$J$-holomorphic map $\ti u: \Delta \to X$ such that
\begin{itemize}
\item $\pi$ is a covering of some degree $1\leq d\leq \mu$, $d|\mu$, with $p$ being a single
branching point ($d=1$ corresponds to the trivial case when $u$ is an imbedding itself);
\item $\ti u \circ \pi = u|_W$;
\item the map $\ti u: \Delta \to X$ has multiplicity $1$ at zero and is a topological
imbedding.
\end{itemize}
\end{lem}
\proof Choose local complex coordinates $(w_1, \ldots, w_n)$ in a
neighborhood of $u(p) \in X$ such that the complex structure $J\st$
defined by $(w_1, \ldots, w_n)$ coincides with $J$ at the point
$u(p)$. Let $z$ be a local complex coordinate on $S$ in a
neighborhood $W \subset S$ of $p$. We may assume that $(w_1, \ldots,
w_n)$ (resp.\ $z$) are centered at $u(p)$ (resp.\ at $p$). By
\eqqref(norm-form-map), after an appropriate rotation and rescaling
of coordinates $w_1,...,w_n$ in $\cc^n$, we can write $u$ in the
form
\begin{equation}
\eqqno(4.5)
u(z) = e_1z^{\mu} + O(|z|^{\mu + \alpha}).
\end{equation}
Let $p_1:\cc^n\to \cc$ be the canonical projection onto the coordinate plane $<e_1>$.
Then $u_1(z)\deff (p_1\circ u)(z)$ is a ramified covering of degree $\mu$ at the origin.
Really, from \eqqref(norm-form-map) and Lemma \ref{l2p-reg} we have that
\[
u_1(z) = z^{\mu}(1 + g(z))
\]
where $g\in L^{2,p}_{loc}$ and $g(z) = O(|z|)$. From this we obtain
that $1 + g(z)$ is an imbedding in a neighborhood of zero and it
admits a root of degree $\mu$, \ie $1 + g(z) = \big(1 +
f(z)\big)^{\mu}$ for some $f\in L^{2,p}_{loc}$ and $f(z) = O(|z|)$.
Therefore we can write
\begin{equation}
\eqqno(4.6)
u_1(z) = z^{\mu}\big(1 + f(z)\big)^{\mu}
\end{equation}
and
\begin{equation}
\eqqno(4.7)
u_1^{\frac{1}{\mu}}: z\to z\big(1 + f(z)\big)
\end{equation}
is an imbedding. If $u|_W$ is an imbedding for some neighborhood $W\ni p$ then the Lemma is trivial
with $d=1$. Suppose now that $u|_W$ is not an imbedding for any neighborhood of $p$.
Take $z_1\not= z_2$ near $p=0$ such that $u(z_1)=u(z_2)$. This implies $u_1(z_1) = u_1(z_{2})$
and the latter reads now as
\begin{equation}
\eqqno(4.8)
z_{1}^{\mu}\big(1+f(z_{1})\big)^{\mu} =  z_{2}^{\mu}\big(1+f(z_{2})\big)^{\mu}.
\end{equation}
We supposed that this happens for any $W$ and therefore we can find two sequences $z_{1,n}\not= z_{2,n}$,
both converging to $0$ and such that $u(z_{1,n})=u(z_{2,n})$. Therefore, after extracting a subsequence,
in view of \eqqref(4.8) we have

\begin{equation}
\eqqno(4.9)
z_{1,n}\big(1+f(z_{1,n})\big) = \eta^kz_{2,n}\big(1+f(z_{2,n})\big)
\end{equation}
for some $0\leq k\leq\mu -1$, where $\eta = e^{\frac{2\pi i}{\mu}}$. If $k=0$ then from
\eqqref(4.7) and \eqqref(4.9) we obtain that $z_{1,n}=z_{2,n}$ for all $n$ and this
is not our case.

\smallskip Therefore we have that $0<k<\mu$. We consider
$z(1+f(z))$ as a new holomorphic coordinate in a neighborhood of $p$. Therefore
for $u$ in this coordinate \eqqref(4.9) means that for some sequence $z_n\to 0$ one has
\begin{equation}
u(\eta^kz_n) = u(z_n).
\end{equation}
From Proposition \ref{compar1} we obtain that there exists a holomorphic reparameterization
$\phi $ such that in new coordinates $u(\eta^kz)\equiv u(z)$, \ie $u$ is multiple of
multiplicity $d=\mu /k$.

\smallskip\qed
Let us remark that proving the last lemma we also proved the following
\begin{corol}
\label{discrt-one}
Let $u:\Delta\to (X,J)$ be a primitive $J$-holomorphic map with
$J\in\calc^{LIp}$. Then for every $0<r<1$ the set $\{ (z_1,z_2)\in
\Delta^2_r: z_1\not=z_2  \text{ and } u(z_1)=u(z_2)\}$
is finite.
\end{corol}
\proof Suppose not. Then there exist two sequences $z_{1,n}\not= z_{2,n}$
converging to $z_1$ and respectively to $z_2$, both in $\Delta$, such
that $u(z_{1,n})=u(z_{2,n})$ for all $n$.

\smallskip\noindent{\sl Case 1. $z_1\not= z_2$.} In that case  the
statement of the Corollary follows from the Theorem B, just proved,
applied to the restrictions of $u$ onto a non intersecting
neighborhoods of $z_1$ and $z_2$. Really, let $V_1\ni z_1$ and
$V_2\ni z_2$ be such neighborhoods. Set $u_1\deff u|_{V_1}$,
$u_2\deff u|_{V_2}$. After translation and rescaling we can suppose
that both $u_1$ and $u_2$ are defined on the unit disc. Theorem B
now applies  and implies that $u_1$ and $u_2$ are not distinct.
Therefore $u$ is not primitive. Contradiction.

\smallskip\noindent{\sl Case 1. $z_1 = z_2$.} This case was considered in the proof of
Lemma \ref{factoriz}. In that case $u$ occurs to be non-primitive.
Contradiction.

\smallskip\qed

\medskip\noindent
{\slsf Step 3}. {\it Construction of the surface $\wt S$ and a
primitive map $\ti u: \wt S \to X$.}\\[3pt]
Consider the set $\scrv$ of pairs $(V, u_V)$ such that $V$ is an
abstract complex curve and $u_V: V \to X$ is a \emph{primitive}
holomorphic map with the image $u_V(V)$ lying in $u(S)$.  We write
$V \in \scrv$ meaning $(V, u_V) \in\scrv$. Take the disjoint union
$\sqcup_{V\in \scrv} V$ and define the following equivalence
relation on $\wt S$: points $p_1 \in V_1 \in\scrv$ and $p_2 \in V_2
\in \scrv$ are identified if there exist $V_3 \in \scrv$, a point
$p_3 \in V_3$ and holomorphic imbeddings $\phi_1: V_3
\hookrightarrow V_1$, $\phi_2: V_3 \hookrightarrow V_2$, such that
$\phi_i(p_3) = p_i$ and the both compositions $u_{V_i} \circ \phi_i$
give $u_{V_3} : V_3 \to X$. Define $\wt S := \sqcup_{V\in
 \scrv} V /\sim$, denote the
natural projections $V \hookrightarrow \wt S$ by $\pi_V$, and equip
the set $\wt S$ with the {\sl quotient topology} whose basis form
the images $\pi_V(V) \subset\wt S$ with $V\in \scrv$. It follows
from the construction of $\wt S$ that there exists a continuous map
$\ti u: \wt S \to X$ such that $\ti u \circ \pi_V = u_V : V \to X$
for any $V\in \scrv$.

The primitivity of the map $\ti u: \wt S \to X$ follows from the
definition of  $\wt S$.

\medskip\noindent
{\slsf Step 4}. {\it $\wt S$ is Hausdorff and there exists a natural
complex structure $\ti\jmath$ on $\wt S$ such that for every $V \in
\scrv$ the projection $\pi: V_V \to \wt S$
is $(j, \ti\jmath)$-holomorphic and such that the map $\ti u: \wt S
\to X$ is $J$-holomorphic.}\\[3pt]
Let $\ti p_1$ and $\ti p_2$ be two distinct points on $\wt S$. Fix
their representatives $p_i \in V_i \in \scrv$. If $\ti u(\ti p_1)
\not= \ti u(\ti p_2)$, then there exist disjoint neighborhoods $\ti
u(\ti p_1) \in W_1 \subset X$ and $\ti u(\ti p_2) \in W_2\subset X$.
Since $\ti u: \wt S \to X$ is continuous, the pre-images $U_i := \ti
u\inv (W_i)$ are open in $\wt S$. Then $U_i$ are desired disjoint
neighborhoods of  $\ti p_1$ and $\ti p_2$.

Now assume that  $\ti u(\ti p_1) = \ti u(\ti p_2)$.  Then by {\sl
Step 2} there exists neighborhoods $\ti p_i \in U_i \subset V_i$
such that  $\ti u(\ti p_1) = \ti u(\ti p_2)$ is the only
intersection point of $\ti u(U_1)$ and $\ti u(U_2)$. It follows from
the definition of the topology on $\wt S$ that $U_i$ are desired
disjoint neighborhoods of $\ti p_1$ and $\ti p_2$.

\smallskip
By the construction, for every $V\in \scrv$ the map $\pi_V : V \to
\wt S$ is an open imbedding so that each $V$ is an open chart for
$\wt S$. We claim that the complex structures on  $V\in \scrv$
induce a well-defined structure $\ti\jmath$ on  $\wt S$. For this
purpose it is sufficient to consider the case $V_1 \subset V_2$.
Since the map $\ti u:V_2 \to X$ is $\calc^1$-regular, the complex
structure on $V_1$ is determined by the structure $J$ on $X$ at each
point $\ti p\in V_2$ with $d\ti u(\ti p) \neq 0$. Thus the inclusion
$V_1 \subset V_2$ is holomorphic outside  the set of critical point
of $\ti u$, which is discrete. Now we use the fact that the
extension of a complex structure over an isolated point is unique
(if exists).

Finally, we observe that  $\ti u: \wt S \to X$ is
$(\ti\jmath,J)$-holomorphic.

\medskip\noindent
{\slsf Step 5.} {\it Construction of the projection $\pi : S \to\wt
S$.} Consider the set $\scrw$ consisting of pairs $(W, \pi_W)$ in
which $W$ is an open subset in $S$ and $\pi_W: W\to \wt S$ is a
holomorphic map such that $\ti u \circ \ti \pi_W = u|_W: W\to X $.
Since $u: S \to X$ is non-constant, it is locally an imbedding
outside the discrete set of critical points of $u$. Using the fact
of the primitivity of $\ti u: \wt S \to X$ we conclude that  $\pi_W:
W\to \wt S$ is unique if exists. In particular, $\pi_{W_1}$ and
$\pi_{W_2}$ must coincide on each intersection $W_1 \cap W_2$ so
that there exists the maximal piece $W_\max := \cup_j W_j$ with the
map $\pi_\max: W_\max \to
 \wt S$. By {\sl Step 1}, $W_\max$ is the whole surface
$S$.

\smallskip\qed

\newprg[prgPP.corol]{Corollaries}

\smallskip
The same proof gives the following variation of Theorem A:

\begin{thm}
Let $(S_1, j_1)$ and $(S_2, j_2)$ be smooth connected complex curves
and $u_i : (S_i,j_i) \to (X,J)$ non-constant $J$-holomorphic maps
with $J \in \calc^{Lip}$. If there are non-empty open sets $U_i \subset S_i$
with $u_1(U_1)= u_2(U_2)$, then there exists a smooth {\slsf
connected} complex curve $(S, j)$ and a $J$-holomorphic map $u:
(S,j) \to (X,J)$ such that $u_1(S_1) \cup u_2(S_2) =u(S)$ and $u:S
\to X$ is primitive.

Moreover, maps $u_i :S_i \to X$ factorize through $u: S \to X$, \ie
there exist holomorphic maps $g_i : (S_i,j_i) \to (S,j)$ such that
$u_i = u \scirc g_i$.
\end{thm}

\smallskip\qed

In the case of closed $J$-complex curves we obtain the following:

\begin{corol}
\label{closed}
Let $(S,j)$ be a connected closed complex curve and let $u:(S,j)\to (X,J)$ 
be a non-constant $J$-holomorphic map into an almost complex manifold $X$
with Lipschitz-continuous almost complex structure $J$. Then there exists
a connected closed complex curve $(\tilde S,\tilde j)$, a ramified
covering $\pi : S\to \tilde S$ and a {\slsf primitive} $J$-holomorphic
map $\tilde u : \tilde S \to X$ such that $u = \tilde u \circ \pi$.
\end{corol}

The following corollary is immediate.

\begin{corol}
Let $u_i: S_i \to (X, J)$, $i=1, 2$ be closed irreducible
$J$-complex curves with $J\in \calc^{Lip}$, such that $u_1(S_1)=M_1\not=
u_2(S_2)=M_2$.  Then
they have finitely many intersection points and the~intersection
index in any such point is strictly positive. Moreover, if $\mu_1$
and $\mu_2$ are the~ multiplicities of $u_1$ and $u_2$ in such
a~point $p$, then the~intersection number of $M_1$ and $M_2$ in $p$
is at least $\mu_1\cdot \mu_2$.
\end{corol}

\newsect[sectOPT]{Optimal Regularity in Lipschitz Structures}

In this section we shall prove the Theorem C from the Introduction.

\newprg[prgOPT.prelim]{Preliminaries}
Consider the Cauchy-Green operator $T_{CG}=\frac{1}{\pi z}*(\cdot
)$:
\begin{equation}
\label{CG-oper}
(T_{CG}u)(z):=\frac{1}{2\pi
i}\int_{\cc}\frac{u(\zeta)}{\zeta-z}d\zeta\land d\bar\zeta .
\end{equation}
$T_{CG}$ is a bounded operator from $\calc^{k,\alpha }(\Delta ,
\cc^n)$ to $\calc^{k+1,\alpha }(\Delta , \cc^n)$ for $0<\alpha <1$.
In particular, there exists $H_{k,\alpha}$ (the norm of $T_{CG}$)
such that
\begin{equation}
\label{CZ-est}
\norm{T_{CG}u}_{\calc^{k+1,\alpha}(\Delta)}\leq
H_{k,\alpha}\norm{u}_{\calc^{k,\alpha}(\Delta)}
\end{equation}
for all $u\in \calc^{\alpha}(\Delta)$.

\smallskip We shall need also the Calderon-Zygmund operator
\begin{equation}
\label{CZ-oper} (T_{CZ}u)(z):=p.v.\frac{1}{2\pi
i}\int_{\cc}\frac{u(\zeta)}{(\zeta-z)^2}d\zeta\land d\bar\zeta .
\end{equation}
It is a bounded operator in spaces $\calc^{k,\alpha}(\Delta)$ and
$L^{k,p}(\Delta)$ and its norm in these spaces will be denoted as
$G_{k,\alpha}$ and $G_{k,p}$ correspondingly.

\smallskip Next consider the Cauchy operator

\begin{equation}
\label{C-oper}
\left( T_Cu\right) (z) = \frac{1}{2\pi
i}\int\limits_{\d\Delta}\frac{u(\zeta)}{\zeta - z}d\zeta.
\end{equation}
$T_C$ is a bounded operator from $\calc^{k,\alpha}(\d \Delta)$ to
$\calc^{k,\alpha}(\Delta)$. For all these facts we refer to
\cite{MP}. We have the following {\slsf Cauchy-Green Formula:} for
$u\in\calc^1(\bar\Delta)$ and $z\in\Delta$

\begin{equation}
\label{CG-form}
u(z) = \left( T_Cu\right) (z) + \left(
T_{CG}\frac{\d u}{\d\bar z}\right) (z).
\end{equation}

Via the Cauchy-Green formula the differential equation \eqqref(J-holo2)
is equivalent to the following integral one:

\begin{equation}
u = T_{C}u +  T_{CG}Q(J_u(z))\overline{\frac{\d u}{\d z}}.
\eqqno(J-holo3)
\end{equation}

\newprg[prgOPT.proof]{Approximation by smooth curves} We fix $J$-holomorphic
$u:\Delta\to\rr^{2n}$ supposing that $u(0)=0$ and that $u$ is
defined in a neighborhood of $\bar\Delta$. Let $B$ be a closed ball
containing the image $u(\bar\Delta)$. Remark that since
$J\in\calc^{\alpha}$ for any $0<\alpha <1$ then by the standard
regularity of $J$-complex curves $u\in
\calc^{1,\alpha}(\bar\Delta)$. We also suppose that $J(0)=J\st$.
Considering dilatations $J_{\delta}(u)=J(\delta u)$ we can suppose
that  $Lip(J)$ is as small as we wish.
Rescaling $u$ by $u_{\delta,\eps}\deff\frac{1}{\delta}u(\eps z)$ be
also can suppose that $\calc^{1,\alpha}(\Delta)$-norm of $u$ is as
small as we wish with $u$ staying to be $J_{\delta}$-holomorphic.
The proof will be achieved via approximation of $J$ in Lipschitz
norm by smooth (of class $\calc^{1,\alpha}$) structures.

\begin{lem}
\label{approx-sol1}
There exists an $\eps>0$ such that if $Lip(J), \norm{u}_{\calc^{1,\alpha}(\Delta)}<\eps$
then for any almost complex structure $\tilde J$ of class $\calc^{1,\alpha}$ on
$B$, standard at origin and such that $||\tilde J - J||_{\calc^{Lip}(B)}<\eps$
there exists a $\tilde{J}$-holomorphic $\tilde u:\bar\Delta \to B$ such that $\tilde
u(0)=0$ and
\begin{equation}
\tilde u(z) = (T_Cu)(z) - (T_Cu)(0) + T_{CG}\left[ Q(\tilde J(\tilde
u)) \overline{\frac{\d \tilde u}{\d z}}\right](z) - T_{CG}\left[
Q(\tilde J(\tilde u))\overline{\frac{\d \tilde u}{\d z}}\right](0).
\eqqno(u-tilde)
\end{equation}
\end{lem}
\proof Actually \eqqref(u-tilde) implies that $\tilde u$ is
$\tilde{J}$-holomorphic and $\tilde u(0)=0$. Therefore all we need is to
construct a solution of \eqqref(J-holo3). In order to do so set $u_0(z) = u(z)$
and define by iteration
\begin{equation}
\label{u_n}
u_{n+1}(z) = (T_Cu)(z) - (T_Cu)(0) + T_{CG}\left[
Q(\tilde J(u_n))\overline{\frac{\d u_n}{\d z}}\right](z) -
T_{CG}\left[ Q(\tilde J(u_n))\overline{\frac{\d u_n}{\d
z}}\right](0).
\end{equation}
We want to prove that $u_n$ converge to a solution of \eqqref(J-holo3).
First we need a uniform bound on
$\norm{u_n}_{\calc^{\alpha}(\Delta)}$ and $\norm{\frac{\d u_n}{\d
z}}_{\calc^{\alpha}(\Delta)}$.

\medskip\noindent{\sl Step 1. Estimate of
$\norm{u_n}_{\calc^{\alpha}(\Delta)}$.}

\smallskip Set $q=\norm{Q}_{\endr{Mat(2n\times 2n,\rr)}}$ and write
\[
\norm{\frac{\d u_{n+1}}{\d z}}_{L^p(\Delta)}\leq
C\norm{u}_{\calc^{1,\alpha}(\Delta)} + \norm{T_{CZ}Q(\tilde
J(u_n))\overline{\frac{\d u_{n}}{\d z}}}_{L^p(\Delta)}\leq
C\norm{u}_{\calc^{1,\alpha}(\Delta)} +
\]
\[
+ q\eps G_p\norm{\frac{\d u_{n}}{\d z}}_{L^p(\Delta)}\leq
C\norm{u}_{\calc^{1,\alpha}(\Delta)} +q\eps
G_pC\norm{u}_{\calc^{1,\alpha}(\Delta)} +(q\eps G_p)^2\norm{\frac{\d
u_{n-1}}{\d z}}_{L^p(\Delta)}\leq
\]
\[
\leq ... \leq C\norm{u}_{\calc^{1,\alpha}(\Delta)}\sum_{k=1}^{n}(q\eps
G_p)^k + (q\eps G_p)^{n+1}\norm{\frac{\d u}{\d z}}_{L^p(\Delta)} \leq
C\norm{u}_{\calc^{1,\alpha}(\Delta)},
\]
if $\eps >0$ was chosen small enough, \ie $q\eps G_p<\frac{1}{2}$.
At the same time from (\ref{u_n}) we see immediately that
$\norm{\frac{\d u_{n+1}}{\d\bar z}}_{L^p(\Delta)}\leq q\eps
\norm{\frac{\d u_{n}}{\d z}}_{L^p(\Delta)}\leq Cq\eps
\norm{u}_{\calc^{1,\alpha}(\Delta)}$. Taking into account the fact
that $u_n(0)=0$ we obtain from the Sobolev Imbedding
$L^{1,p}(\Delta)\subset \calc^{\alpha}(\Delta)$ the desired bound
\begin{equation}
\label{c-alpha}
\norm{u_n}_{\calc^{\alpha}(\Delta)}\leq
C\norm{u}_{\calc^{1,\alpha}(\Delta)}
\end{equation}
with $C$ independent on $n$. Here $p>2$ should be taken at the very
beginning satisfying $\alpha = 1-\frac{2}{p}$.

\medskip\noindent{\sl Step 2. Estimate of $\norm{\frac{\d u_n}{\d
z}}_{\calc^{\alpha}(\Delta)}$.}

\smallskip
\[
\norm{\frac{\d u_{n+1}}{\d z}}_{\calc^{\alpha}(\Delta)} \leq
C\norm{u}_{\calc^{1,\alpha}(\Delta)} + G_{\alpha} \norm{Q(\tilde
J(u_n))\overline{\frac{\d u_{n}}{\d z}}}_{\calc^{\alpha}(\Delta)}
\leq C\norm{u}_{\calc^{1,\alpha}(\Delta)} +
\]
\[
+ G_{\alpha}q\eps C\norm{u}_{\calc^{1,\alpha}(\Delta)}\norm{\frac{\d
u_{n}}{\d z}}_{\calc^{\alpha}(\Delta)} \leq
C\norm{u}_{\calc^{1,\alpha}(\Delta)}\big(1 + G_{\alpha}q\eps
\norm{\frac{\d u_{n}}{\d z}}_{\calc^{\alpha}(\Delta)}\big) \leq
\]
\[
\leq C\norm{u}_{\calc^{1,\alpha}(\Delta)}\big(1 + CG_{\alpha}q\eps
\norm{u}_{\calc^{1,\alpha}(\Delta)} +
C\norm{u}_{\calc^{\alpha}(\Delta)}(G_{\alpha}q\eps )^2\norm{\frac{\d
u_{n-1}}{\d z}}_{\calc^{\alpha}(\Delta)}\big) \leq
\]
\begin{equation}
\label{deu-c-alpha} \leq ... \leq
C\norm{u}_{\calc^{1,\alpha}(\Delta)}\sum_{k=1}^{n+1}\big[ 1 +
C(G_{\alpha}q\eps )^k\big]\leq C\norm{u}_{\calc^{1,\alpha}(\Delta)}
\end{equation}
for $\eps >0$ sufficiently small.

\smallskip\noindent{\sl Step 3. Convergence of approximations.}

\smallskip We proved that $\norm{u_n}_{\calc^{\alpha}(\Delta)}, \norm{\frac{\d u_n}{\d
z}}_{\calc^{\alpha}(\Delta)}\leq C$ if $\eps >0$ and
$\norm{u}_{\calc^{1,\alpha}(\Delta)}$ were taken small enough. Now
we can write
\[
\norm{u_{n+1}-u_n}_{\calc^{1,\alpha}(\Delta)}\leq
2H_{\alpha}\norm{Q(\tilde J(u_n))\overline{\frac{\d u_n}{\d z}}
-Q(\tilde J(u_{n-1}))\overline{\frac{\d u_{n-1}}{\d
z}}}_{\calc^{\alpha}} \leq
\]

\[
\leq 2H_{\alpha}\norm{Q(\tilde J(u_n))\left[ \overline{\frac{\d
u_n}{\d z}} - \overline{\frac{\d u_{n-1}}{\d z}}\right]
}_{\calc^{\alpha}} + 2H_{\alpha}\norm{ \left[ Q(\tilde
J(u_n))-Q(\tilde J(u_{n-1}))\right] \overline{\frac{\d u_{n-1}}{\d
z}}}_{\calc^{\alpha}}\leq
\]

\begin{equation}
\leq 2CH_{\alpha}q\eps\norm{u_n-u_{n-1}}_{\calc^{1,\alpha}(\Delta)} +
2CH_{\alpha}q \eps \norm{u_n-u_{n-1}}_{\calc^{\alpha}(\Delta)}.
\eqqno(5.15)
\end{equation}
For $\eps>0$ small enough we obtain

\begin{equation}
\norm{u_{n+1}-u_n}_{\calc^{1,\alpha}(\Delta)}\leq
r\norm{u_{n}-u_{n-1}}_{\calc^{1,\alpha}(\Delta)}
\eqqno(5.16)
\end{equation}
with some fixed $0<r<1$. Therefore $\{ u_n\} $ converge in
$\calc^{1,\alpha}(\Delta)$ to a solution $\tilde u$ of
\eqqref(u-tilde). Lemma is proved.

\smallskip\qed

\begin{lem}
\label{approx-sol2}
Let $\{ J_n\} $ be a sequence of almost complex structures on $B$ of
class $\calc^{1,\alpha}$, standard at origin, converging to $J$ in
$\calc^{Lip}(B)$. Let $u_n$ be some solution of \eqqref(u-tilde) for
$J_n$. Then $\norm{u_n - u}_{\calc^{1,\alpha}(\Delta)}\to 0$.
\end{lem}
\proof Since $u$ also satisfies \eqqref(J-holo3) we can write 
$\norm{u_n - u}_{\calc^{1,\alpha}} \leq$
\[
\leq 2H_{\alpha}\norm{Q(J_n(u_n))\overline{\frac{\d u_n}{\d z}} -
Q(J(u))\overline{\frac{\d u}{\d z}}}_{\calc^{\alpha}} \leq
2H_{\alpha}\norm{Q(J_n(u_n))}_{\calc^{\alpha}} \norm{\frac{\d
u_n}{\d z} - \frac{\d u}{\d z}}_{\calc^{\alpha}} +
\]
\[
 +  2H_{\alpha}\norm{\frac{\d u}{\d z}}_{\calc^{\alpha}}\norm{Q(J_n(u)) - Q(J(u))}_{\calc^{\alpha}} \leq
2H_{\alpha}q\eps \norm{u_n - u}_{\calc^{1,\alpha}} + C\norm{J_n -
J}_{\calc^{\alpha}}.
\]
And this implies
\[
\norm{u_n - u}_{\calc^{1,\alpha}} \leq \frac{C}{1 -
2H_{\alpha}q\eps}\norm{J_n - J}_{\calc^{\alpha}}\to 0.
\]

\smallskip\qed

\begin{rema} \rm 
Remark that $u_n$ have regularity $\calc^{2,\alpha}$.
\end{rema}

\newprg[prgOPT.loglip]{Log-Lipschitz convergence of approximating sequence}

\begin{lem}
\label{approx-sol3}
Let $u_n$ and $J_n$ be as in Lemma \ref{approx-sol2}. Then $\{ u_n\}$ are uniformly 
bounded in $\calc^{1,LnLip}(\bar\Delta)$.
\end{lem}
\proof To prove this statement we need to recall one useful formula. For a
smooth function $\lambda$ on an almost complex manifold $(X,J)$ the
$1$-form $d^c_J\lambda$ is defined by

\begin{equation}
d^c_J\lambda (v)=-d\lambda (Jv) \eqqno(5.17)
\end{equation}
for every tangent vector $v$. If $J$ is of class $\calc^1$ then
$dd^c_J\lambda$ is then defined by usual differentiation. As usual,
$\Delta =\frac{\partial^2}{\partial x ^2} + \frac{\partial^2}{\partial y
^2}$ will denote the Laplacian on the plane $\cc$. The notation
$d^c=d^c_{J\st}$ is relative to the standard complex structure on
$\cc$. So for a function $\lambda$ defined on an open set of $\cc$:
$d^c\lambda =-\frac{\partial \lambda }{\partial y} dx+ \frac{\partial
\lambda }{\partial x} dy$. Therefore $dd^c\lambda = \Delta\lambda
dx\land dy$.

\smallskip This can be generalized to the functions on $X$ as follows.
Let $J$ be a  $\calc^1$-regular almost complex structure defined on
an open set $\Omega \subset \rr^{2n}$ and let $\lambda$ be a
$\calc^2$ function defined on $\Omega$ (as $\lambda$ we shall take coordinate
functions $u_1,...,u_{2n}$ in $\rr^{2n}$ to deduce the needed regularity of 
$u:\Delta\to\rr^{2n}$). If $u:\Delta \to (\Omega
,J)$ is a $J$-holomorphic map, then:

\begin{equation}
\Delta (\lambda\circ u)(z)~=~[dd^c_J\lambda ]_{u (z)} \Big(\frac{\partial
u}{\partial x}(z),J_{u(z)}\frac{\partial u}{\partial x}(z)\Big).
\eqqno(5.18)
\end{equation}
For the proof we refer to \cite{IR1}.

\smallskip Apply the formula \eqqref(5.18) to the $J_n$-holomorphic
mapping $u_n$ obtained above:

\begin{equation}
\Delta (\lambda\circ u_n)(z)~=~[dd^c_{J_n}\lambda ]_{u_n (z)}
\Big(\frac{\partial u_n}{\partial x}(z),J_n(u_n(z))\frac{\partial u_n}{\partial x}(z)\Big).
\eqqno(5.19)
\end{equation}
Since $J_n$ converge to $J$ in Lipschitz norm their first
derivatives are uniformly bounded on $B$ and therefore the right
hand side of \eqqref(5.19) shows that for any smooth function
$\lambda$ on $\rr^{2n}$ Laplacians $\{ \Delta (\lambda\circ u_n)\}$
are uniformly bounded on $\Delta$ for all $n$. Lemma 1.7 from
\cite{IR1} gives now that $\{ \lambda\circ u_n\}$ are bounded in
$\calc^{1,LnLip}(\Delta)$. As $\lambda$ we can take any coordinate
function $u_k$ on $\rr^{2n}$ and obtain the desired statement.

\medskip To finish the proof of Theorem C all is left is to remark
that if $u_n\to u$ uniformly and $\{ u_n\}$ stay bounded in
$\calc^{1,LnLip}(\Delta)$ then $u\in \calc^{1,LnLip}(\Delta)$.

\medskip\qed

\begin{rema}\rm
\rm If $J$ is of class $\calc^1$ then the following statement holds
true, see \cite{IS2}. Let $u:\Delta\to X$ be a $J$-holomorphic map and
let $(E,J_1\deff u^*J)$ be the induced bundle. Then this complex vector
bundle has a natural structure of a holomorphic bundle and $du$ is a
holomorphic morphism of holomorphic bundles $T\Delta\to E$.

\smallskip The approximations made in the proof of Theorem C permit to
extend this statement to the case of Lipschitz-continuous $J$. Really,
for a given $J$-holomorphic $u:\Delta\to X$ of class $\calc^{1, LnLip}(\Delta)$
we constructed a sequence $J_n$ of smooth structures converging to $J$ in
Lipschitz norm and a sequence of $J_n$-holomorphic $u_n:\Delta\to X$
converging to $u$ in the space $\calc^{1, LnLip}(\Delta)$. That means
that holomorphic structures constructed in \cite{IS2} will converge and
$du_n$ will converge to an analytic morphism of sheaves $du$.
\end{rema}

\newsect[sectPC]{Perturbation of a Cusp}

\newprg[prgPC.inverse]{Inversion of a $\dbar$-type operators}
We shall perturb a cusp of a  $J$-holomorphic map
$u_0:(\Delta ,0)\to (\cc^n,0)$, which we suppose to be given in the
form \eqqref(norm-form-map):
\[
u_0(z) = z^{\mu}v(z), \qquad v_0\deff v(0)\not=0,
\]
where $v\in L^{1,p}_{loc}$, $zv\in L^{2,p}_{loc}$.
We shall use
perturbations of cusps in several different ways in this paper. Our first aim
is to perturb $u_0$ in such a way that the perturbed map $u$ stays to be
$J$-holomorphic and has no cusps. For that aim we should search for a perturbation
$u$ of $u_0$ in the form
\begin{equation}
\eqqno(pert1)
u(z) = u_0(z) + z\cdot w(z),
\end{equation}
where $w(0) = w_0\not= 0$ and $w_0$ is not collinear to $v_0$. Such perturbations will be used
in the following section for the proof of the Genus Formula. Later, for deriving an
``essential part of a Puiseux series'' we will need to perturb $u_0$ adding a term
of an arbitrary degree and along a tangent which may be collinear to $v_0$.  Of course,
we are interested only in $J$-holomorphic perturbations. We start with the following
\begin{prop}
\label{morrey1}
If, under the assumptions of Proposition \ref{morrey} (for $k=1$),  the sum of the  norms
$\norm{J-J\st}_{\calc^{Lip}(\Delta)} + \norm{R}_{L^{1,p}(\Delta)}$ is sufficiently small
then:

\sli There exists a linear bounded operator
$T_{J,R}^0:L^{1,p}(\Delta)\to L^{2,p}(\Delta)$ such that $(\dbar_J + R)\circ T_{J,R}^0\equiv \id$
and $(T_{J,R}^0u)(0) = 0$ for every $u\in L^{1,p}(\Delta)$;

\slii The same operator acts also from $\calc^{\alpha}(\Delta)$ to $\calc^{1,\alpha}(\Delta)$ with
the same properties.
\end{prop}

For $J=J\st$ and $R=0$ the operator in question is $T_{J\st,0}^0(u) = T_{CG}u - (T_{CG}u)(0)$,
where $T_{CG}$ is the standard Cauchy-Green operator. For general $J,R$ the operator
$T_{J,R}^0$ can be constructed as the perturbation series:
\begin{equation}
\eqqno(pert-op)
T_{J,R}^0 \deff \sum_{n=0}^{\infty} (-1)^nT_{J\st,0}^0\circ\big((\dbar_J - \dbar_{J\st} + R\circ
T_{J\st,0}^0\big)^n.
\end{equation}

\newprg[prgPC.main]{Proof of the main result} Let us state and prove the main result of 
this section.

\begin{thm}
\label{cusp-pert}
Let $J$ be a Lipschitz-continuous almost complex structure in the unit ball $B\subset \cc^n$
with $J(0)=J\st$ and let $u_0:\Delta\to B$ be a $J$-holomorphic map. Let $\nu\geq 0$ be an integer
and $w_0\in\cc^n$ be a vector. Then there exists a $J$-holomorphic map $u :\Delta_r\to B$,
defined in a smaller disc $\Delta_r$, such that
\begin{equation}
\eqqno(pert2)
u(z) = u_0(z) + z^{\nu}\cdot w(z),
\end{equation}
with $w(0)=w_0$ and $w\in L^{1,p}_{loc}$ for any $p<\infty$.
\end{thm}
\proof Let us apply the Cauchy-Riemann operator to $u(z)$ in the form \eqqref(pert2):
\[
0 = \dbar_{J\circ u}u = \d_xu + J(u)\d_yu = \d_xu +
\big(J(u)-J(u_0)\big)\d_yu + J(u_0)\d_yu =
\]
\[
= \d_xu_0 + J(u_0)\d_yu_0 + \d_x(z^{\nu}w) + J(u_0)\d_y(z^{\nu}w) +
\big(J(u)-J(u_0)\big)\d_y(u_0 + z^{\nu}w) =
\]
\[
= \dbar_{J\circ u_0}u_0 + \dbar_{J\circ u_0}(z^{\nu}w) +
\big(J(u)-J(u_0)\big)\big(\d_yu_0 + \d_y(z^{\nu}w)\big).
\]
Therefore we need to solve the equation
\begin{equation}
\eqqno(pert3)
\dbar_{J\circ u_0}(z^{\nu}w) =
\big(J(u_0)-J(u)\big)\big(\d_yu_0 + \d_y(z^{\nu}w)\big).
\end{equation}
Multiplying by $z^{-\nu}$ we obtain
\begin{equation}
\eqqno(pert4)
z^{-\nu}\dbar_{J\circ u_0}(z^{\nu}w) = z^{-\nu} \big(J(u_0) - J(u)\big)\big( \d_yu_0 +
\d_y(z^{\nu}w)\big).
\end{equation}
The left hand side of \eqqref(pert4) can be transformed as follows
\[
z^{-\nu}\dbar_{J\circ u_0}(z^{\nu}w) = \big(\d_x + z^{-\nu}J(u_0)z^{\nu}\d_y\big)w
+ \nu z^{-\nu}\big(z^{\nu -1} + J(u_0)J\st z^{\nu -1}\big)w =
\]
\[
= \big(\d_x + z^{-\nu}J(u_0)z^{\nu}\d_y\big)w
+ \nu z^{-\nu}\big(1 + J(u_0)J\st \big)z^{\nu -1}w =: \dbar_{J^{(\nu )}\circ u_0} w +
R^{(\nu )}w,
\]
where $J^{(\nu )}\deff z^{-\nu}J(u_0)z^{\nu}$ is Lipschitz-continuous by Lemma \ref{lipschitz-A}
and $R^{(\nu)}$ admits an obvious pointwise estimate $|R^{(\nu)}(z)|\leq
\nu Lip(J)\norm{u_0}_{\calc^{1,\alpha}}$. 
Therefore the left hand side of \eqqref(pert4) has the
form
\begin{equation}
D_{J^{(\nu)},u_0}(w) \deff \dbar_{J^{(\nu )}} w + R^{(\nu )}w,
\end{equation}
for a Lipschitz-continuous $J^{(\nu )}$ and a bounded $R^{(\nu)}$ with $\norm{R^{(\nu)}}_{L^{\infty}}$
small. The right hand side
\[
F^{(\nu)}(z,w) \deff z^{-\nu} \big(J(u_0) - J(u)\big)\big( \d_yu_0 +
\d_y(z^{\nu}w)\big)
\]
of \eqqref(pert4) admits the following estimates:

\begin{equation}
\eqqno(est1)
\norm{F^{(\nu)}(z,w)}_{L^p(\Delta)}\leq C\cdot Lip(J)\norm{u_0}_{\calc^{1,\alpha}(\Delta)}
\norm{w}_{L^{1,p}(\Delta)};
\end{equation}
\begin{equation}
\eqqno(est2)
\norm{F^{(\nu)}(z,w)}_{\calc^{\alpha}(\Delta)}\leq C\cdot Lip(J)\norm{u_0}_{\calc^{1,\alpha}(\Delta)}
\norm{w}_{\calc^{1,\alpha}(\Delta)};
\end{equation}
\begin{equation}
\eqqno(est3)
\norm{F^{(\nu)}(z,w_1) - F^{(\nu)}(z,w_2) }_{L^p/\calc^{\alpha}(\Delta)}\leq
C\cdot Lip(J)\norm{u_0}_{\calc^{1,\alpha}} \norm{w_1 - w_2}_{L^{1,p}/\calc^{1,\alpha}(\Delta}.
\end{equation}
Our goal is to solve the following equation
\begin{equation}
\eqqno(syst1)
\begin{cases}
D_{J^{(\nu)},u_0}w = F^{(\nu)}(z,w),\cr
w(0) = w_0.
\end{cases}
\end{equation}
We can apply Newton's method of successive approximations by
setting
\begin{equation}
\eqqno(newton)
w_{n+1} = T_{J^{(\nu)},R^{(\nu)}}^0\big[F^{(\nu)}(z,w_n)\big] + w_1,
\end{equation}
where $w_1$ is to be found as a solution of the following system
\begin{equation}
\eqqno(syst2)
\begin{cases}
D_{J^{(\nu)},u_0}w_1 = 0, \cr
w_1(0) = w_0.
\end{cases}
\end{equation}
I.e.,
\[
w_1(z) = w_0 - T_{J^{(\nu)},R^{(\nu)}}^0\big(D_{J^{(\nu)},u_0}w_0\big).
\]
Estimates \eqqref(est1), \eqqref(est2), \eqqref(est3) guarantee the convergence of the
iteration process. The proof is very similar to that of the previous section.
As it was explained there we can suppose that $Lip(J)$
as well as $\norm{u_0}_{\calc^{1,\alpha}(\Delta)}$  are as small as we wish, less then
some $\eps>0$ to be specified in the process of the proof. We can also suppose that
$\norm{R^{(\nu)}}_{L^{1,p}(\Delta)}\leq \eps$.
$\norm{w_0}$ will be supposed also small enough. Finally, we shall suppose
inductively that $\norm{w_n}_{L^{1,p}(\Delta)}\leq \frac{1}{2}$.

\smallskip As in the proof of Lemma \ref{approx-sol1} we start with estimating first the
$L^p$ and then $\calc^{\alpha}$-norms of derivatives.

\smallskip\noindent{\slsf Step 1. There exists a constant $C$, independent of $n$,
such that  $\norm{w_n}_{\calc^{\alpha}(\Delta)}\leq C\norm{w_0}$.}

\[
\norm{\frac{\d w_{n+1}}{\d \bar z}}_{L^p(\Delta)} = \norm{\dbar_{J\st}w_{n+1}}_{L^p(\Delta)}
= \norm{\dbar_{J^{(\nu)}}w_{n+1} + (J\st - J^{(\nu)}(u))\d_yw_{n+1}}_{L^p(\Delta)} =
\]
\[
= \norm{(\dbar_{J^{(\nu)}} + R^{(\nu)})w_{n+1} - R^{(\nu)}w_{n+1} + (J\st - J^{(\nu)}(u))
\d_yw_{n+1}}_{L^p(\Delta)} \leq
\eps \norm{w_{n+1}}_{L^p(\Delta)} +
\]
\[
+ \norm{F^{(\nu)}(z,w_n)}_{L^{p}(\Delta)}
+ Lip(J^{(\nu)})\norm{u_0+z^{\nu}w_n}_{L^{\infty}(\Delta)}\norm{\nabla
w_n}_{L^{p}(\Delta)} \leq
\]
\[
\leq \eps \norm{w_{n+1}}_{L^p(\Delta)} + C\eps \norm{w_{n}}_{L^{1,p}(\Delta)}.
\]
Further
\[
\norm{\frac{\d w_{n+1})}{\d z}}_{L^p(\Delta)} \leq  \norm{\frac{\d w_{1}}{\d z}}_{L^p(\Delta)}
+ \norm{\frac{\d}{\d z}T^0_{J^{(\nu)},R^{(\nu)}}[F^{(\nu)}(z,w_n)]}_{L^p(\Delta)} =
\]
\[
=  \norm{\frac{\d}{\d z}T^0_{J^{(\nu)},R^{(\nu)}}[D_{J^{(\nu)},u_0}w_0]}_{L^p(\Delta)} +
\norm{\frac{\d}{\d z}T^0_{J^{(\nu)},R^{(\nu)}}[F^{(\nu)}(z,w_n)]}_{L^p(\Delta)}\leq
\]
\[
\leq C \norm{D_{J^{(\nu)},u_0}w_0}_{L^p(\Delta)} + \norm{F^{(\nu)}(z,w_n)}_{L^p(\Delta)}
\leq C \norm{(\dbar_{J^{(\nu)}} + R^{(\nu)})w_0}_{L^p(\Delta)} +
\]
\[
+ C\eps\norm{w_n}_{L^{1,p}(\Delta)}
\leq C\eps\norm{w_0} + C\eps\norm{w_n}_{L^{1,p}(\Delta)}.
\]
Taking into account that $w_{n+1}(0)=w_0$ we obtain
\begin{equation}
\eqqno(6.12)
\norm{w_{n+1}}_{L^{1,p}(\Delta)}\leq C\big(\norm{w_0} + \eps\norm{w_n}_{L^{1,p}(\Delta)}\big),
\end{equation}
From \eqqref(6.12) we obtain
\begin{equation}
\eqqno(6.13)
\norm{w_{n+1}}_{L^{1,p}(\Delta)}\leq C\norm{w_0}\sum_{k=0}^n(C\eps )^k +
(C\eps )^{n+1}\norm{w_0} \leq C\norm{w_0},
\end{equation}
with $C$ independent of $n$. This justifies our inductive assumption that
$\norm{w_n}_{L^{1,p}(\Delta)}\leq \frac{1}{2}$ and implies in its turn the needed estimate
\begin{equation}
\eqqno(6.14)
\norm{w_{n+1}}_{\calc^{\alpha}(\Delta)}\leq C\norm{w_0}.
\end{equation}

\smallskip\noindent{\slsf Step 2. There exists a constant $C$ independent of $n$ such that
$\norm{\nabla(w_n)}_{\calc^{\alpha}(\Delta)}\leq C\norm{w_0}$.}

\smallskip\noindent Using computations of the Step 1  write
\[
\norm{\frac{\d(w_{n+1})}{\d\bar z}}_{\calc^{\alpha}(\Delta)} \leq
\norm{F^{(\nu)}(z,w_n)}_{\calc^{\alpha}(\Delta)} + \norm{R^{(\nu)}w_n}_{\calc^{\alpha}(\Delta)}
+ \norm{(J\st - J^{(\nu)}(u))\d_yw_n}_{\calc^{\alpha}(\Delta)}
\]
\begin{equation}
\eqqno(6.15)
\leq Lip(J)\norm{u_0}_{\calc^{\alpha}(\Delta)}\norm{w_n}_{\calc^{1,\alpha}(\Delta)} +
\eps\norm{w_n}_{\calc^{\alpha}(\Delta)} + Lip(J)\norm{w_n}_{\calc^{1,\alpha}(\Delta)} \leq
\eps\norm{w_n}_{\calc^{1,\alpha}(\Delta)}.
\end{equation}

Analogously to $L^p$-case write further
\[
\norm{\frac{\d w_{n+1})}{\d z}}_{\calc^{\alpha}(\Delta)} \leq
\norm{\frac{\d w_{1}}{\d z}}_{\calc^{\alpha}(\Delta)}
+ \norm{\frac{\d}{\d z}T^0_{J^{(\nu)},R^{(\nu)}}[F^{(\nu)}(z,w_n)]}_{\calc^{\alpha}(\Delta)} =
\]
\[
=  \norm{\frac{\d}{\d z}T^0_{J^{(\nu)},R^{(\nu)}}[D_{J^{(\nu)},u_0}w_0]}_{\calc^{\alpha}(\Delta)} +
\norm{\frac{\d}{\d z}T^0_{J^{(\nu)},R^{(\nu)}}[F^{(\nu)}(z,w_n)]}_{\calc^{\alpha}(\Delta)}\leq
\]
\[
\leq C \norm{D_{J^{(\nu)},u_0}w_0}_{\calc^{\alpha}(\Delta)} +
\norm{F^{(\nu)}(z,w_n)}_{\calc^{\alpha}(\Delta)}
\leq C \norm{(\dbar_{J^{(\nu)}} + R^{(\nu)})w_0}_{\calc^{\alpha}(\Delta)} +
\]
\begin{equation}
\eqqno(6.16)
+ C\eps\norm{w_n}_{\calc^{1,\alpha}(\Delta)}
\leq C\eps\norm{w_0} + C\eps\norm{w_n}_{\calc^{1,\alpha}(\Delta)}.
\end{equation}

From \eqqref(6.15) and \eqqref(6.16) we obtain
\begin{equation}
\eqqno(6.17)
\norm{\nabla (w_{n+1})}_{\calc^{\alpha}(\Delta)}\leq
C\eps\big(\norm{w_0} + \norm{\nabla (w_n)}_{\calc^{\alpha}(\Delta)}\big)
\end{equation}
and therefore
\[
\norm{\nabla (w_{n+1})}_{\calc^{\alpha}(\Delta)}\leq
C\norm{w_0} + C\eps\norm{w_0} + \eps^2\norm{\nabla (w_{n-1})}_{\calc^{\alpha}(\Delta)}
\leq ... \leq
\]
\begin{equation}
\eqqno(6.18)
\leq C\norm{\nabla (w_{0})}_{\calc^{\alpha}(\Delta)} \leq C\norm{w_0}.
\end{equation}
We conclude these two steps with the estimate:
\begin{equation}
\label{c1-alpha}
\norm{\nabla (w_{n})}_{\calc^{\alpha}(\Delta)}\leq C\norm{w_0}
\end{equation}
with $C$ independent on $n$, provided $Lip(J)$ and
and $\norm{w_0}$ are small enough.

\smallskip\noindent{\slsf Step 3. Convergence of approximations.}

\smallskip Write
\[
\norm{w_{n+1}-w_n}_{\calc^{1,\alpha}(\Delta)}\leq C
\norm{F^{(\nu)}(z,w_n) - F^{(\nu)}(z,w_{n-1})}_{\calc^{\alpha}(\Delta)}\leq
\]
\begin{equation}
\eqqno(6.20)
 \leq  C\cdot\eps \norm{w_{n}-w_{n-1}}_{\calc^{1,\alpha}(\Delta)}
\end{equation}
by \eqqref(est3) with $\eps >0$ as small as we wish. This gives us the desired convergence of
$w_n$ to a solution $w$ of \eqqref(syst1).

\smallskip\qed

\begin{rema}\rm
 Let $w'$ and $w''$ be solutions of \eqqref(syst1) with initial data
 $w'(0)=w^{'}_0$ and $w''(0)=w^{``}_0$.  Then, as in \eqqref(6.20), we have
\[
\norm{w' -w''}_{\calc^{1,\alpha}(\Delta)}\leq \norm{w^{'}_0-w^{``}_0} + C
\norm{F^{(\nu)}(z,w') - F^{(\nu)}(z,w'')}_{\calc^{\alpha}(\Delta)}\leq
\]
\[
\leq  \norm{w^{'}_0-w^{``}_0} + \eps \norm{w' -w''}_{\calc^{1,\alpha}(\Delta)}.
\]
And therefore
\begin{equation}
\label{cont-dep1}
\norm{w' - w''}_{\calc^{1,\alpha}(\Delta)}\leq \frac{1}{1-\eps} \norm{w'_0-w''_0},
\end{equation}
\ie a solution $zw$ of \eqqref(syst1) continuously depend on the initial data
$w(0)=w_0$.  In particular we have
\begin{equation}
\label{cont-dep2}
\norm{w}_{\calc^{1,\alpha}(\Delta)}\leq \frac{1}{1-\eps} \norm{w_0},
\end{equation}
for the solution with $w(0)=w_0$.
\end{rema}

In the following lemma we suppose that $u_0(z) = z^{\mu}v(z)$ with $v(0)=v_0=e_1$.
\begin{lem}
Let $Lip(J)$ and $a\in\cc$ be small enough. Set $w_0=ae_2$. Then the $J$-holomorphic curve 
$u(z)=u_0(z) + zw$ has no cusps, where $zw$ is a solution of \eqqref(syst1) 
with $\nu =1$ and initial data $w(0)=w_0$.
\end{lem}
\proof In an appropriate coordinates we have $u_0(z) = z^{\mu}e_1 +
z^{2\mu -1}v(z)$. After making a dilatations $J_{\delta}(z)\deff
J(\delta^{\mu}z)$ and $u_0^{\delta}(z)=
\frac{1}{\delta^{\mu}}u_0(\delta z)$ we can suppose that
$Lip(J)\leq\eps$ - as small as we wish. Moreover, \eqqref(2.2) gives
us the behavior of $\norm{zv}_{L^{2,p}(\Delta (r)}$ and therefore we
can estimate the differential of $u_0$ in the following way:
\begin{equation}
du (z) = \mu z^{\mu -1}e_1 + R(z) \text{ with } |R(z)|\leq \eps |z|^{\mu -1 +\alpha}
\text{ for } |z|\leq \frac{1}{2}.
\end{equation}
Let $zw$ be a solution of \eqqref(syst1) with $w(0)=w_0e_2$. 
Remark that it satisfies (\ref{cont-dep2}), \ie
\begin{equation}
\norm{zw}_{\calc^{1,\alpha}(\Delta)}\leq C\norm{w_0}
\end{equation}
and therefore its differential can be written as
\begin{equation}
d(zw) = w_0e_2 + P(z),
\end{equation}
where $\norm{P(z)}\leq C\norm{w_0}|z|^{\alpha}$.

\smallskip
With these data we need to show that the differential 
of the $J$-holomorphic map $u(z)= u_0(z) + zw(z)$
does not  vanishes in $\Delta_{\frac{1}{2}}$. Let us  write this differential:
\begin{equation}
\label{diff-du}
du (z) = \mu z^{\mu -1}e_1 + R(z) + w_0e_2 + P(z).
\end{equation}

We use the following notations: $R(z) \deff R_1(z)e_1 +
R_2(z)$, 
where $R_2(z)$ takes values in the subspace {\sl span}$\{e_2,...,e_n\}$ of $\cc^n$.
And the same for $P(z)$.
First of all, since $\norm{P(z)}\leq C\norm{w_0}|z|^{\alpha}$ we see that
there exists $0<r_0<\frac{1}{2}$ such that
\begin{equation}
\norm{w_0 + R_2(z) + P_2(z)} \geq \norm{w_0}(1-C|z|^{\alpha} - C|z|^{\mu - 1 + \alpha})> 0
\end{equation}
for all $|z|\leq r_0$ (independently of $w_0$!). This gives us that 
the second coordinate of the differential
is not vanishing for $|z|\leq r_0$. At the same time
\begin{equation}
|\mu z^{\mu - 1} + R_1(z) + P_1(z)| \geq \mu r_0^{\mu -1} - \eps - C\norm{w_0} > 0
\end{equation}
if $w_0$ and $\eps$  where taken sufficiently small. Therefore the first
coordinate of the differential does not
vanishes for $|z|\geq r_0$.

\smallskip\qed

This lemma permits us to define the cusp index of a cusp point of a
$J$-complex curve.

\begin{defi}
\label{cusp-ind}
Let $u_0(z) = z^{\mu}v_0 + O(|z|^{\mu + \alpha})$,  $\mu\geq 1$, $v_0\not=0$ be a
$J$-complex curve and let $u$ be a
small perturbation of $u_0$ as in Lemma 6.1 which has no cusps. The cusp index
$\varkappa_0$ of $u_0$ at zero is
defined as the sum of intersection indices of self-intersection points of such
a perturbation $u$.
\end{defi}

In the following section we shall see that this number doesn't depend on a perturbation (provided it
is sufficiently small).

\newsect[sectGF]{Genus Formula in Lipschitz Structures}

\newprg[prgGF.num]{Local numeric invariants} 
To state the Genus Formula we need to define local numeric invariants of 
$J$-complex curves and to insure that these invariants are positive (otherwise such
``formula'' will be useless). The problem is that we need to do this  in the case when $J$ is only
Lipschitz-continuous. The first invariant --- the local intersection number --- was introduced in
Definition \ref{int-number} and in Theorem B it was proved that this number is always positive
and is equal to $1$ if and only if the local intersection in question is transverse.
The second --- the cusp index for a cusp point --- was defined at the end of the previous section,
see Definition \ref{cusp-ind}.
Now we need to prove that it doesn't depend on perturbation. It will be done by relating it
to the Bennequin index. Much more details of this approach can be found in
\cite{IS1} and we suggest that the interested reader has the latter preprint in his hands
while reading this section.

\newprg[prgGF.ben]{The Bennequin index of a cusp}

Let $u : (\Delta , 0)\to (\cc^2, 0) $ be a germ of a non-constant $J$-complex
curve at zero (and $J$ is Lipschitz). Without loss of generality we always suppose that
$J(0)=J\st$. Taking into account that zeros of $du$ are isolated, we can suppose
that $du$ vanishes only at zero. Furthermore, let $w_1, w_2$ be the standard
complex coordinates in $(\cc^2, J\st)$. We already used several times in this paper the
following presentation of $u$ and its differential $du$:
\begin{equation}
u(z) = z^\mu\cdot a + O(| z|^{\mu +\alpha}) \quad \text{ and } \quad
du (z) = \mu z^{\mu -1} a +  O(| z|^{\mu -1 +\alpha}).
\eqqno(7.1)
\end{equation}
Here $a$ is a non-zero vector in $\cc^2$, $\mu \geq 2$ and $0<\alpha <1$.

\smallskip For $r>0$ define $F_r\deff TS^3_r \cap J(TS^3_r)$ to be the distribution
of $J$-complex planes in the tangent bundle $TS^3_r$ to the sphere of radius
$r$. $F_r$ is trivial, because $J$ is homotopic to $J\st = J(0)$. By $F$
we denote the distribution $\cup_{r>0}F_r\subset \cup_{r>0}TS^3_r \subset
TB^*$, where  $TB^*$ is the tangent bundle to the punctured ball in $\cc^2$.
Set $M=u(\Delta )$.
\begin{lem}
\label{transvers}
The (possibly not connected) curve $\gamma_r = M\cap S^3_r$ is
transverse  to $F_r$ for all sufficiently small $r>0$.
\end{lem}

\proof Since $J\approx J\st$ for $r$ sufficiently small, $T\gamma_r$ is close to $J\st
n_r$, where $ n_r$ is the field of normal vectors to $S^3_r$.
On the other hand, for sufficiently small $r$, the distribution $F_r$ is
close to the one of $J\st$-complex planes in $TS^3_r$, which is orthogonal
to $J\st n_r$.

\smallskip\qed

\smallskip
This fact permits us to define the Bennequin index of $\gamma_r$. Namely, take
any non-vanishing on $S^3_r$ section $\vec v$ of $F_r$ and move $\gamma_r$ along the vector
field $\vec v$ to obtain a curve $\gamma'_r$. We can make this move for a small
enough time, so that $\gamma'_r $ does not intersect $\gamma_r$.
\begin{defi}
\label{ben-index}
The {\slsf Bennequin index $b(\gamma_r)$} is the linking number of $\gamma_r$
and $\gamma'_r$.
\end{defi}
This number does not depend on $r>0$, taken sufficiently small, because
$\gamma_r$ is homotopic to $\gamma_{r_1}$ for $r_1<r$ within the curves
transverse to $F$, see \cite{Bn}. It is also independent of the particular choice
of the field $\vec v$. For the standard complex structure $J\st$ in
$B\subset \cc^2$ we use $\vec v\st(w_1,w_2)\deff  (-\bar w_2, \bar w_1)$
for calculating the Bennequin index of the curves
on sufficiently small spheres. For an arbitrary almost complex structure $J$
with $J(0)=J\st$ we can find the vector field $\vec v_J$, which is defined
in a small punctured neighborhood of the origin, is a small perturbation
of $\vec v\st$, and lies in the~distribution $F$ defined by $J$.

\smallskip The following statement is crucial for proving that the quantity
$\varkappa_0 = \frac{b(\gamma_r)+1}{2}$ is a well defined and non-negative numerical
invariant of a cusp.  Denote by $B_{r_1, r_2}$ the spherical shell $B_{r_2}\bs
\bar B_{r_1}$ for $r_1<r_2$.

\begin{lem}
\label{ben=int}
Let $\Gamma $ be an immersed $J$-complex  curve
in a neighborhood of $\barr B_{r_1, r_2}$ such that all self
intersection points of $\Gamma $  are contained in $B_{r_1, r_2}$ and 
all components of the curves $\gamma_{r_i} := \Gamma \cap
S^3_{r_i}$ are transverse to $F_{r_i}$ for $i=1,2$. Then
\begin{equation}
b(\gamma_{r_2}) = b(\gamma_{r_1}) + 2\cdot \sum\nolimits _{x\in
Sing(\Gamma) }\delta_x, \eqqno(7.2)
\end{equation}
where the sum is taken over self-intersection points of $\Gamma $.
\end{lem}

\proof Move $\Gamma $ a little along $\vec v_J$ to obtain $\Gamma^\epsi$. By
$\gamma^\epsi_{r_1}, \gamma^\epsi_{r_2}$ denote the intersections
$\Gamma^\epsi \cap S^3_{r_1},\Gamma^\epsi\cap S^3_{r_2}$, which are of
course the moves of $\gamma_{r_j}$ along $v_J$. We have $l(\gamma_{r_2},\gamma
_{r_2}^\epsi) - l(\gamma_{r_1},\gamma_{r_1}^ \epsi) = \sf{int}(\Gamma,\Gamma^\epsi)$, 
where $l(\cdot ,\cdot )$ is the~linking number and $\sf{int}(\cdot ,\cdot )$ is the 
intersection number, see \cite{Rf}.

Now let us calculate $\sf{int}(\Gamma ,\Gamma^\epsi)$. From Theorem B we know that there
are only a finite number $\{ p_1,\ldots, p_N\} $ of self-intersection points of $\Gamma
$. Take one of them, say $p_1$. Let $M_1, \ldots, M_d$ be the discs on $\Gamma $ with a
common point $p_1$ and otherwise mutually disjoint.  More precisely we take
$M_j$ to be irreducible components of $\Gamma \cap B_{\rho}(p_1)$ for $\rho >0$ small
enough. Remark that $M_j$ are transverse to $v_J$, so their moves $M_j^\epsi$
do not intersect them, \ie $M_j\cap M_j^\epsi = \emptyset $. Note also that for $k\not=
j$
we have $\sf{int}(M_k, M_j) = \sf{int}(M_k, M^\epsi_j)$ for $\epsi >0$
sufficiently small.  Therefore $\sf{int}(\Gamma \cap B_{\rho}(p_1), \Gamma^\epsi\cap B_{\rho}(p_1))
= \sum_{1\leq k<j\leq d} \sf{int}(M_k, M^\epsi_j) + \sf{int}(M^\epsi_k, M_j)= 2\sum_{1\leq
 k<j\leq d} \sf{int}(M_k, M^\epsi_j) = 2\delta_{p_1}$. This means that $\sf{int}(\Gamma
,\Gamma^\epsi) = 2\cdot \sum_{j=1}^N\delta_{p_j}$.

\smallskip\qed

Now we are ready to describe the cusp-index in a different way. Let $u$ be a
parameterization of $M$ near $p$. Take a small ball $B_r(p)$
around $p=u(0)$ and a small disc $\Delta$ centered at zero such that $u(\Delta) =
M\cap B_r(p)$. More precisely $u(\Delta)$ is the irreducible component of
$M\cap B_r(p)$ containing the cusp $p$. We take $r>0$ small enough, such that
$\Delta$ contains no other critical points of $u|_{\Delta}$ then the origin and that
$u(\Delta)$ has no self-intersections. Let $\gamma_r\deff u(\Delta ) \cap \d B_r(p)$
and $b_p$ be the~Bennequin index of $\gamma_r$, defined in Definition  \ref{ben-index}.
Lemma \ref{ben=int} together with the obvious fact that the Bennequin
index of a smooth point is $-1$ tell us that the number 
\[
\varkappa_p\deff (b_p +1)/2
\] 
is well defined, non-negative and is equal to the number of double points of a generic
perturbation.

\smallskip

Let us define the local invariants of a $J$-complex curve $M$ in an almost
complex surface. From Theorem B and Corollary 5.1 it follows readily that
a compact  $J$-complex curve with a finite number of irreducible components $M =
\bigcup_{i=1}^d M_i$ has only 
a finite number of local self-intersection points, provided $J$ is Lipschitz-continuous.

For each such point $p$ we can introduce, according to Definition
\ref{int-number}, the self-in\-ter\-sec\-tion number $\delta_p(M)$
of $M$ at $p$. Namely, let $S_j$ be a parameter curve for $M_j$, \ie
$M_j$ is given as an image of the $J$-holomorphic map $u_j : S_j\to
M_j$. We always suppose that the 
parameterization $u_j$ is primitive, \ie they cannot be decomposed like $u_j=v_j\scirc r$
where $r$ is a nontrivial covering of $S_j$ by another Riemann
surface. Denote by $\{ x_1,\ldots, x_N\}$ the set of all pre-images
of $p$ under $u : \bigsqcup_{i=1}^d S_j\to X$, and take mutually disjoint discs 
$\{ D_1,\ldots, D_N\}$ with centers $x_1,\ldots, x_N$
such that their images have no other common points different from
$p$. For each pair $D_i, D_j$, $i\not= j$, define an intersection
number as in Definition \ref{int-number} and take the sum over all different
pairs to obtain $\delta_p(M)$.

Now put $\delta = \sum_{p\in D(M)}\delta_p(M)$, where the sum is taken over the
set $D(M)$ of all local intersection points of $M$, \ie points which have at least two
pre-images. Consider now the set $\{ p_1,\ldots, p_L\} \subset \bigcup_{j=1}^d
S_j = S$ of all cusps of $M$, \ie points 
where the differential of the appropriate parameterization vanishes. Set $\varkappa\deff 
\sum_{i=1}^L \varkappa_i$.

\smallskip
Numbers $\delta$ and $\varkappa$ are the local numerical invariants of $M$ involved in
the Genus Formula.

\newprg[prgGF.genus]{Genus Formula for $J$-Complex Curves}

Denote by $c_1(X, J)$ the first Chern class of $X$ with respect to
$J$. Since, in fact, $c_1(X, J)$ does not depend on continuous changes of $J$
we usually omit the dependence of $c_1(X)$ on $J$.

\begin{lem}
\label{genus-immers}
Let $M=\bigcup_{j=1}^d M_j$ be a compact immersed $J$-complex curve in a four-dimensional
almost complex manifold $(X, J)$ with Lipschitz continuous $J$. Then
\begin{equation}
\sum_{j=1}^d g_j = \frac{[M]^2 - c_1(X)[M]}{2} + d - \delta .
\eqqno(genus-immers)
\end{equation}
\end{lem}
For the proof see, ex. \cite{IS1},  \cite{MW}, or any other text.

\smallskip The proof of the general Genus Formula \eqqref(genus-form) will be
reduced to the immersed case via perturbations. To do so we need
the following ``matching'' lemma from \cite{IS1}.
Let $B(r)$ be a ball of radius $r$ in $\rr^4$ centered at zero, and $J_1$ a
Lipschitz continuous almost complex structure on $B(2)$, $J_1(0)=J\st$. 
Further, let $M_1 =
u_1(\Delta )$ be a closed primitive $J_1$-complex disc in $B(2)$ such that
$u_1(0)=0$ and $M_1$ transversely meet $S^3_r$ for $r\geq 1/2$. Here $S^3_r = \d
B(r) $ 
and transversality are understood with respect to both $TS^3_r$ and
$F_r$.

By $B(r_1,r_2)$ we shall denote the spherical shell $\{ x\in \rr^4:r_1< \norm{
x} <r_2 \}$. 
In the lemma below denote by $D_{1+\delta }$ the pre-image of
$B(1+\delta )$ by $u_1$.

\begin{lem}
\label{matching}
For any positive $\delta >0$ there exists an $\epsi >0$ such
that if an almost complex structure $J_2$ in $B(1+\delta )$ and a closed
$J_2$-holomorphic curve $M_2$ parameterized by $u_2 : D_{1+\delta }\to B(1+\delta )$ satisfy
$\norm{J_2 - J_1}_{\calc^1(\bar B(1+\delta ))} < \epsi $ and
$\norm{u_2 - u_1}_{L^{1,p}(D_{1+\delta })} < \epsi $, then
there exists an almost complex structure $J$ in $B(2)$ and $J$-holomorphic
disc $M$ in $B(2)$ such that:

a) $J\ogran_{B(1-\delta )} = J_2\ogran_{B(1-\delta )}$ and
$J\ogran_{B(1+\delta ,2)} = J_1\ogran_{B(1+\delta ,2)}$.

b) $M\ogran_{B(1-\delta )} = M_2\cap B(1-\delta )$ and
$M\cap B(1+\delta ,2) = M_1\cap B(1+\delta ,2)$.
\end{lem}

\proof We have chosen the parameterization of $M_1$ to be primitive. Thus, $u_1$ is
an imbedding on $D_{1-\delta ,1+\delta } = u_1^{-1}(B_{1-\delta ,1+\delta })$. Let us
identify a neighborhood $V$ of $u(D_{-\delta ,\delta })$ in $B_{1-
\delta ,1+\delta }$ with the neighborhood of the zero-section in the
normal bundle $N$ to $u(D_{1-\delta ,1+\delta })$. Now $u_2\ogran_{D_{-\delta
,\delta}}$ can be viewed as a section of $N$ over $u(D_{-\delta ,\delta })$
which is small \ie contained in $V$. Using an appropriate smooth function $\phi
$ on $D_{1-\delta ,1+\delta }$ (or equivalently on $u(D_{1-\delta ,1+\delta} )$),
$\phi \ogran _{B(1-\delta ) \cap D_{1+\delta }}\equiv 1$, $\phi\ogran_{\d
D_{1+\delta}}\equiv 0$, $0\leq \phi \leq 1$ we can glue $u_2$ and $u_1$ to obtain a
{\slsf symplectic} surface $M$ which satisfies (b).

Patching $J_1$ and $J_2$ and simultaneously making $M_1$  complex
can be done in an obvious way.

\smallskip\qed

\medskip\noindent{\sl Proof of the Genus Formula.} Using Lemma \ref{cusp-pert}
we perturb every irreducible component $M_j$ near each of its cusp and using Lemma
\ref{matching} we glue perturbed pieces back to compact curves and denote them again
by $M_j$. The perturbed structure will be still denoted as $J$. The sum $\delta$
of local intersection indices did not change and by Lemma \ref{ben=int} each cusp $p$
with cusp-index $\varkappa_p$ produces a finite set of intersection points with the
sum of intersection indices equal to $\varkappa_p$. Now the Lemma \ref{genus-immers}
gives us the proof of the general case.

\smallskip\qed

\newsect[sectSING]{Structure of singularities of pseudoholomorphic curves}

In this section we define an analogue of the {\slsf Puiseux series} for
{\it primitive} $J$-holomorphic curves in Lipschitz-continuous almost complex structure $J$.

\newprg[prgSING.puis]{Puiseux series of holomorphic curves} %

\smallskip%
It is known that for a germ of an irreducible complex curve $C$ in $\cc^n$ at
the origin $0$ there exist a local holomorphic reparameterization of $\cc^n$
and a parameterization of $C$ by a non-multiple holomorphic map $u:\Delta\to\cc^n$
such that the first component of $u(z)$ is $z^{p_0}$, whereas all remaining
components have order $>p_0$. In other words $u(z)$ writes as

\begin{equation}
\eqqno(puis1)
\textstyle u(z) =
(z^{p_0\vph}, v_1z^{p_1\vph}+v_2z^{p_2\vph}+\cdots),
\end{equation}
with some non-vanishing $v_i\in\cc^{n-1}$ and $p_{i+1}>p_i$ for $i\geq 0$.
Introducing a new variable $t:=z^{p_0}$ we can
write $u(t)=(t, f_2(t^{1/p_0\vph}),\ldots,f_n(t^{1/p_0\vph}))$, or simply
\begin{equation}
\eqqno(puis0)
u(t)=(t, f(t^{1/p_0\vph})),
\end{equation}
where $f$ is a holomorphic function with values in $\cc^{n-1}$.  The representation \eqqref(puis0) is
called the {\slsf Puiseux series} of $u$ at
$0\in\Delta$. We refer to  \cite{Co}, Book II, Chapter II for a nice exposition on
Puiseux series.
Another reference is  \cite{Fi}, Chapter 7.

\medskip%
The following consideration explains the idea for the generalization of the
notion of Puiseux series to the case of pseudoholomorphic curves.  The
exponents $(p_0,p_1,\ldots)$ of the non-vanishing terms $v_iz^{p_i}$ determine the
topological type of the singularity of $C$ at $0$. In particular, making
non-vanishing deformations of the coefficients $v_i$ we obtain an {\slsf
 equisingular deformation} of the curve $C=u(\Delta)$ such that $0=u(0)$ remains
the only singular point and the cusp index $\varkappa_0$ persists.  However, some of
exponents $p_i$ are non-essential for the singularity type. That means that
the type and the cusp index $\varkappa_0$ remains unchanged if the corresponding $v_i$
vanishes and the term $v_iz^{p_i}$ disappears. The other exponents, called
{\slsf characteristic} or {\slsf essential exponents of the singularity} $0\in
C$, admit the following two criteria.

The first criterion is: $p_i$ is a characteristic exponent in a sequence
$p_0<p_1<\cdots<p_l$ if and only if the sequence $d_j:=\gcd(p_0,\ldots,p_j)$ decreases
after $d_i$, \ie $d_{i+1}<d_i$. The second criterion is as follows:
consider approximations of the parameterizing map $u(z)$ of the form
$u(z)-\tilde u(z^d)=O(z^p)$ such that $\tilde u(z):\Delta_r\to\cc$ is a primitive holomorphic map
in some (small) disc and $p\geq p_0,d>1$ are integers. In particular, $\tilde u(z^d)$ is a
$d$-multiple holomorphic map. Call such an approximation {\slsf extremal} if
there exist no other approximation $u(z)-\tilde u'(z^d)=O(z^{p'})$ with the same
{\slsf multiplicity} $d$ and higher {\slsf degree} $p'>p$ and no other approximation
$u(z)-\tilde u''(z^{d''})=O(z^p)$ with the same degree $p$ and higher multiplicity $d''>d$.  
It is not difficult
to show that the degree $p_i$ of such an extremal approximation is exactly one
of the characteristic exponents, and then the corresponding multiplicity is
$d_{i-1}=\gcd(p_0,\ldots,p_{i-1})$. This second characterization follows
immediately from the Puiseux series.

\begin{rema} \rm Strictly speaking it is not immediately clear that extremal
 approximations do exist.  We shall prove their existence in the following
 subsection, see Lemmas \ref{easy-apprx} and \ref{mult-apprx}.
\end{rema}

\newprg[prgSING.type]{Multiple approximations $J$-complex curves}

\smallskip%
We use the second criterion for extremal exponents of the Puiseux series as a
model for our constructions in pseudoholomorphic case. Till the end of this
section $J$ will be a Lipschitz-continuous almost complex structure in the
unit ball $B\subset\cc^n$ with $J(0)=J\st$ and $u(z):\Delta\to B$ a {\it primitive}
$J$-holomorphic map, written in the form
\begin{equation}
\eqqno(8.3)
u(z)=v_0z^\mu+O(|z|^{\mu+\alpha}) \text{ with } \mu\geq2 \text{ and } v_0\neq 0\in\cc^n.
\end{equation}
Further, the relation $f(z)=O(|z|^{\mu +\alpha})$ will be understood as ``$f(z)=O(|z|^{\mu +\alpha})$
for every $0<\alpha<1$''. Similarly, notation ``$w(z)\in L^{1,p}(\Delta,\cc)$'' will mean
``$w(z)\in L^{1,p}(\Delta,\cc)$ for every $p<\infty$''. We start with the following easy statement.
\begin{lem}
\label{easy-apprx}
Let $\ti u:\Delta_r\to B$ be a $J$-holomorphic map such that
\[
u(\phi(z))-\ti u(z^d)=z^p\tilde w+O(|z|^{p+\alpha}),
\]
for some holomorphic function $\phi$ of the form $\phi (z) = z+ O(z^2)$, some $d>1$,
and some $\tilde w\in\cc^n$.  If $p>\mu$ then $d$ is a divisor of $\mu$. In
particular, $d\leq \mu$.
\end{lem}

\proof Without loss of generality we may assume that $\phi(z)\equiv z$. Really, we can
consider $u_1\deff u\circ \phi$ instead of $u$. Remark that $v_0$ in \eqqref(8.3) for
$u_1$ will be the same.  Let $\eta:=e^{2\pi\isl/d}$ be the primitive root of unity
of degree $d$. Then $u(\eta z)-u(z)=\ti u(\eta^dz^d)-\ti u(z^d) + O(z^{p})=
O(z^{p})$.  On the other hand, %
$u(z)=v_0z^\mu+O(|z|^{\mu +\alpha})$ and hence $u(\eta z)-u(z)=v_0(\eta^\mu-1)z^\mu + O(|z|^{\mu
 +\alpha})$.  Since $p>\mu$, %
this implies that $\eta^\mu=1$. Therefore $d$ is a divisor of $\mu$.

\smallskip\qed

\begin{lem}
\label{mult-apprx}
Let $d$ be a divisor of $\mu$, $\eta$ a primitive root of unity of order $d$, and
\begin{equation}
\eqqno(8.4a)
u(\varphi(\eta z))-u(\varphi(z)) =w(0)z^\nu + O(|z|^{\nu+\alpha})
\end{equation}
the presentation given by Part {\bf (b)} 
of the Comparison Theorem. Further, let $\ti u:\Delta_r\to B$ (for some $r>0$) be a 
$J$-holomorphic map and $\phi(z)$ a holomorphic function in a neighborhood of zero
of the form $\phi(z)=z+O(z^2)$. Assume that
\begin{equation}
\eqqno(8.5a)
u(\phi(z))-\ti u(z^d)=z^p\tilde w+O(|z|^{p+\alpha}),
\end{equation}
with some $\tilde w\not= 0\in\cc^n$ and some $p\in\nn$. Then $p\leq \nu$.

Moreover, in the case $p<\nu$ either the vector $\ti w$ is proportional to $v_0$
or $p$ is a multiple of $d$.
\end{lem}

\proof Recall that by the Comparison Theorem $\nu >\mu$. This gives the proof in the
case $p\leq\mu$. 

Thus we may suppose that $p>\mu$.  In this case Lemma \ref{easy-apprx} says that
that $\lambda=\frac{\mu}{d}$ and $\ti u(z)=z^\lambda v_0+ O(|z|^{\lambda + \alpha})$.

Recall that $\eta$ is the primitive root of unity of degree $d$.  Then
\begin{equation}
\eqqno(8.6a)
u(\phi(\eta z))-u(\phi(z)) =(\eta^p-1)\ti wz^p + O(|z|^{p+\alpha}).
\end{equation}
We want to compare this relation with \eqqref(8.4a). The assertion of the lemma
holds if $\phi(z)\equiv\varphi(z)$ so we assume that this is not the case. Define $\gamma(z)$
from the relation $\phi(z)=\varphi\big( z(1+\gamma(z))\big)$. Then $\gamma(z)$ is given by the
formula $\gamma(z)=(\varphi\inv\circ\phi(z)-z)/z$, where $\varphi\inv(z)$ is the inverse of $\varphi(z)$,
$\varphi\inv{\circ}\varphi(z)\equiv z$. It follows that $\gamma(z)$ is a holomorphic function in some
disc $\Delta_r$ ($r>0$) which is not identically zero and satisfies $\gamma(z)=O(z)$.

Consider first the case $\gamma(z)=\gamma_1(z^d)$. Set $\zeta=z(1+\gamma(z))$. Then $\zeta=z+O(z^2)$,
$\phi(z)=\varphi(\zeta)$, and $\phi(\eta z)=\varphi\big( \eta z(1+\gamma(\eta z)) \big)=\varphi(\eta\zeta)$ since
$\gamma(\eta z)=\gamma_1(\eta^dz^d)=\gamma_1(z^d)=\gamma(z)$. Consequently, 
\[
u(\phi(\eta z))-u(\phi(z)) = u(\varphi(\eta\zeta))-u(\varphi(\zeta))= w(0)\zeta^\nu + O(|\zeta|^{\nu+\alpha})=
 w(0)z^\nu + O(|z|^{\nu+\alpha})
\]
since $\zeta=z+O(z^2)$. Comparing this relation with \eqqref(8.4a) we conclude the
desired inequality $p\leq\nu$. Moreover, in the case $p<\nu$ we also conclude the
relation $\eta^p-1=0$. The latter means that $p$ is a multiple of $d$ which gives
us the second assertion of the lemma

Consider the remaining case. Then $\gamma(z)= \gamma_1(z^d)+bz^k+O(z^{k+1})$ with some
holomorphic $\gamma_1(z)=O(z)$, some $b\neq0\in\cc$, and some $k>0$ which is not a
multiple of $d$. The latter fact is equivalent to $\eta^k-1\neq0$. As above, set
$\zeta=z(1+\gamma(z))$. Then again $\zeta=z+O(z^2)$ and $\phi(z)=\varphi(\zeta)$. On the other hand,
\[
\eta z(1+\gamma(\eta z))-\eta z(1+\gamma(z))=\eta zb((\eta z)^k-z^k)+O(z^{k+2})=\eta b(\eta^k-1)z^{k+1}+O(z^{k+2}),
\]
and hence 
\[
 \phi(\eta z)=\varphi\Big(\eta\zeta +\eta b(\eta^k-1)\zeta^{k+1}+O(\zeta^{k+2}) \Big).
\]
(Here we use the fact that all three functions $\phi(z),\varphi(z)$ and $\zeta(z)$ behave
like $=z+O(z^2)$.) 

\smallskip %
At this point we use the following 

\state Claim. {\it Let $u:\Delta\to B$ be a $J$-holomorphic map of the form
 \eqqref(8.3) and $a(z)$ some function such that $a(z)=O(|z|^{1 +\alpha})$ with
 $\alpha>0$. Then $u(z+a(z)) - u(z) = v_0\mu z^{\mu -1} a(z)+ O(|a|{\cdot}|z|^{\mu - 1 +
  \alpha})$.}

The claim follows from \eqqref(norm-form-diff). Really
\[
u(z+a) - u(z) = a\int_0^1\nabla u(z+ta)dt = av_0\mu z^{\mu -1} + O(|z|^{\mu - 1+\alpha}|a|).
\]

\smallskip %
Let us apply this claim to $u\circ\varphi(\zeta)$ instead of our original map $u(z)$. This
gives us
\[
u(\phi(\eta z))-u(\phi(z)) = u\Big( \varphi\big(\eta\zeta +\eta b(\eta^k-1)\zeta^{k+1}+O(\zeta^{k+2}) \big) \Big) 
- u(\varphi(\zeta))
\]
\[
=  u\Big( \varphi\big(\eta\zeta \ +\ \eta b(\eta^k-1)\zeta^{k+1}  \ +\ O(\zeta^{k+2})\big) \Big) - u( \varphi(\eta\zeta))
\ \ + \ \ u( \varphi(\eta\zeta))-u( \varphi(\zeta)) =
\]
\begin{equation}
\eqqno(8.7)
v_0\mu z^{\mu-1}{\cdot}\eta b(\eta^k-1)z^{k+1} +   w(0)z^\nu  +O(|z|^{k+\mu+\alpha}) + O(|z|^{\nu+\alpha}). 
\end{equation}
We see that \eqqref(8.4a), \eqqref(8.6a), and \eqqref(8.7) are contradictory in
the case $p>\nu$, because $w_0$ is non-zero and orthogonal to $v_0$.

Moreover, in the case $p<\nu$ we conclude the equality of the terms
\[
z^p\ti w = v_0\mu z^{\mu-1}{\cdot}\eta b(\eta^k-1)z^{k+1}
\]
which gives us the desired proportionality $\ti w = \mu \eta b(\eta^k-1){\cdot}v_0$.

\qed

\begin{defi}
\label{mult-ap}
Let $1<d\leq \mu$ be a divisor of $\mu$. A {\slsf multiple approximation} 
of $u$ of {\slsf multiplicity}  $d$  is a primitive $J$-holomorphic map 
$\tilde u:\Delta_r\to B$ such that
\begin{equation}
\eqqno(mult-ap-f)
u(\phi (z)) - \tilde u(z^d) = z^p\tilde w + O(|z|^{p+\alpha}),
\end{equation}
for some holomorphic reparameterization $\phi$ of the form $\phi (z) = z + O(z^2)$
and such that $p>\mu $.
\end{defi}
The degree $p$ in \eqqref(mult-ap-f) depends, in general, on $\phi$ but by Lemma 
\ref{mult-apprx} doesn't exceed $\nu$. Therefore we can give the following:

\begin{defi}
\label{degree-ap}
The maximal possible $p$ in \eqqref(mult-ap-f) is called the {\slsf degree} of the
multiple approximation $\tilde u$.
\end{defi}

Now let us define the principal notion in our approach.

\begin{defi}
\label{extr-apprx}
An approximation $\ti u$ of multiplicity $d$ and degree $p$ is called 
{\slsf extremal} if there exists no other
approximation of the same multiplicity $d$ and higher degree $p_+>p$, and no
other approximation of the same degree $p$ and higher multiplicity $d_+>d$.
\end{defi}
From Lemmas \ref{easy-apprx} and \ref{mult-apprx} it is clear that extremal
approximations do exist.

\smallskip In the case of integrable $J$ the map $u(z)$ itself and any its
multiple approximation $\ti u(z)$ are holomorphic and thus are given by
converging power series. Moreover, making local coordinate change one can
eliminate some non-characteristic terms in the expansion \eqqref(puis1). 
In the case of non-integrable $J$, for two
given maps $u_1(z),u_2(z)$ one can in general define solely one term of their
difference $u_1(z)-u_2(z)=z^\nu v+o(z^\nu)$. In particular, setting $u_2(z)\equiv0$,
one should expect that at most first non-trivial term of the expansion of
$u_1(z)$ is well-defined. Our construction insures that all
characteristic terms are still well-defined.

\newprg[prgSING.proof]{Proof of Theorem E}

\smallskip The proof of the theorem is based on the following lemma, 
which explains how one constructs extremal approximations explicitly.

\begin{lem}\label{pui-apprx} 
 Under hypotheses of Theorem E, let $d>1$ be a divisor of $\mu$, $\eta$ the
 primitive root of unity of degree $d$, and $\nu>\mu$ the number given by the part
 {\bf(b)} of Comparison Theorem. Then there exist $r>0$ and a multiple
 approximation $\ti u :\Delta_r\to B$ such that
\begin{equation}\eqqno(pui-apprx-d)
u(\varphi(z))-\ti u(z^d) =\ti wz^\nu+ O(|z|^{n+\alpha})
\end{equation}
with some reparameterization $\varphi(z)$ of the form $\varphi(z)=z+O(z^2)$ and some
non-zero vector $\ti w\in\cc^n$ orthogonal to $v_0$. In particular, $\ti u(z)$
is a multiple approximation of $u(z)$ of multiplicity $d$ and degree $\nu$.
\end{lem}

\proof Recall that $u(z)=v_0z^\mu+O(|z|^{\mu+\alpha})$. By Theorem \ref{cusp-pert},
there exist $r>0$ and a $J$-holomorphic map $\ti u_0:\Delta_r\to B$ satisfying
$u_0(z)=v_0z+O(|z|^{1+\alpha})$. Set $\ti u_1(z):= \ti u_0(z^{\mu/d})$.  Then
$u(z)-\ti u_1(z^d) =\ti wz^q+ O(|z|^{q+\alpha})$ with some $q>\mu$.

Consider the following more general situation. Let $\ti u_j:\Delta_r\to B$ be a
$J$-holomorphic map such that
\begin{equation}\eqqno(u-phi-p)
u(\varphi_j(z))-\ti u_j(z^d)=\ti w_jz^q+O(|z|^{q+\alpha})
\end{equation}
with the same divisor $d$, some $q$ with $\mu<q\leq\nu$, some $\ti w_j\neq0\in\cc^n$ and
some holomorphic $\varphi_j(z)$ of the form $\varphi_j(z)=z+O(z^2)$. Assume that
\eqqref(u-phi-p) does not satisfy assertion of the lemma. We are going to
describe the construction which shows that an appropriate deformation of %
$\ti u_j$ and $\varphi_j$ refines the situation, such that the iteration of this
construction yields the desired result.

Consider $\gamma(z)=z\cdot(1+a\cdot z^{q-\mu})$ with $a\in\cc$ and set $\varphi_{j+1}(z):=\varphi_j(z\cdot(1+a\cdot
z^{q-\mu}))$. %
Then $u(\varphi_j(\gamma(z)))-u(\varphi_j(z))=v_0a\mu z^q+O(|z|^{q+\alpha})$, see the Claim in proof of Lemma
\ref{mult-apprx}.  Consequently, for an appropriate choice of $a\in\cc$ we
obtain $u(\varphi_{j+1}(z))-\ti u_j(z^d)=\ti w'_{i}z^{q}+O(|z|^{q+\alpha})$ %
where $\ti w'_{i}$ is either vanishing or non-zero and orthogonal to $v_0$.
Notice that the corresponding $a\in\cc$ and $\varphi_{j+1}(z)$ are defined uniquely.

In the case when $\ti w'_j$ vanishes we obtain a new approximation
$u(\varphi_{j+1}(z))-u'(z^d)=\ti w_{j+1}z^{q'}+O(|z|^{q'+\alpha})$ with with $q'>q$. In
this case we repeat the above procedure.

In the case when $\ti w'_j$ is non-zero and orthogonal to $v_0$ and $q=\nu$ our
approximation $u(\varphi_{j+1}(z))-\ti u_j(z^d)=\ti w'_{i}z^{q}+O(|z|^{q+\alpha})$ %
has the desired form.

It remains to consider the case when $\ti w'_j$ is non-zero and orthogonal to
$v_0$, and $q<\nu$. In this case by the second assertion of Lemma
\ref{pui-apprx} $q$ must be a multiple of $d$, $q=d{\cdot}l$. Then by Theorem
\ref{cusp-pert} there exists a $J$-holomorphic map $\ti u_{j+1}:\Delta_{r'}\to B$
which is defined in some (possibly) smaller disc $\Delta_{r'}$ and satisfies %
$\ti u_{j+1}(z)-\ti u_j(z)= \ti w'_j{\cdot}z^{q/d}+ O(|z|^{q/d+\alpha})$. Then
$u(\varphi_{j+1}(z))-\ti u_{j+1}(z^d)=O(|z|^{q+\alpha})$ and hence $u(\varphi_{j+1}(z))-\ti
u_{j+1}(z^d)=\ti w_{j+1}z^{q'}+O(|z|^{q+\alpha})$ with $q'>q$. So this time also we
can repeat our procedure.

Since $q$ is bounded from above by $\nu$, after several repetitions of the
procedure we obtain the desired approximation of the form
\eqqref(pui-apprx-d).

\qed

\state Remark. Notice that for given $\ti u_j(z)$ and $\varphi_j(z)$ satisfying
\eqqref(u-phi-p) the construction $\varphi_{j+1}(z)$ is unique, whereas %
$\ti u_{j+1}(z)$ is unique up to a higher order term
$O(|z|^{q/d+\alpha})$. Furthermore, modifying $\varphi_j(z)$ at each step we add a term
$a\cdot z^{q_j-\mu+1}$, whose degree increases at each step. Consequently, in all
approximations \eqqref(u-phi-p) we can replace all $\varphi_j(z)$ by the final
function $\varphi(z)$ without decreasing the degree of the approximation.

\smallskip \noindent %
{\it Proof of Theorem E.} Let $u:\Delta\to B$ satisfies the hypotheses of
Theorem E.  In particular, $u(z)=v_0z^\mu+O(|z|^{\mu+\alpha})$. Set $p_0:=d_0:=\mu$ and
let $\eta_0:=e^{2\pi\isl/d_0}$ be the primitive corresponding root of unity. Let
$\nu_0$ be the exponent given by Part {\bf(b)} of Comparison Theorem. Set
$p_1:=\nu_0$.  Then by the previous lemma there exist $J$-holomorphic map
$u_0(z)$ of the form $u_0(z)=v_0z+O(|z|^{1+\alpha})$ and a holomorphic function
$\phi_0(z)$ such that $u(\phi_0(z))-u_0(z^{d_0})=v_1z^{p_1}+O(|z|^{p_1+\alpha})$.

From this moment we proceed recursively constructing at each step an
approximation of the form
\begin{equation}\eqqno(apprx-i)
u(\phi_i(z))-u_i(z^{d_i})=v_{i+1}z^{p_{i+1}}+O(|z|^{p_{i+1}+\alpha})
\end{equation}
satisfying the assertion of Theorem E. Since the starting case $i=0$ is
already obtained, we need only to establish the recursive step
$(i)\,\Rightarrow\,(i+1)$. For the divisor $d_i>1$ of $\mu=p_0=d_0$ let $\nu_i$ be the
number given by Comparison Theorem. Then by Lemma \ref{pui-apprx} the number
$\nu_i$ equals the exponent $p_{i+1}$ in the $i$-th approximation
\eqqref(apprx-i). Further, by the Comparison Theorem $p_{i+1}=\nu_i$ is not a
multiple of $d_i$.  Put
$d_{i+1}:=\gcd(d_i,p_{i+1})=\gcd(p_0,p_1,\ldots,p_{i+1})$. In the case $d_{i+1}=1$
we put $l:=i+1$, $\phi(z):=\phi_i(z)$ (the function obtained in the previous step
$(i)$), put $u_i(z):=u(\phi(z))$, and terminate the recursive procedure.

Otherwise we have $d_{i+1}>1$. Let $\eta_{i+1}=e^{2\pi\isl/d_{i+1}}$ be the
corresponding root of unity, and let $\nu_{i+1}$ be the number given by
Comparison Theorem for the divisor $d=d_{i+1}$. Then $\nu_{i+1}$ is not a
multiple of $d_{i+1}$, and we set $p_{i+2}:=\nu_{i+1}$. Then Lemma
\ref{pui-apprx} provides the desired approximation
$u(\phi_{i+1}(z))-u_{i+1}(z^{d_{i+1}})=v_{i+2}z^{p_{i+2}}+O(|z|^{p_{i+2}+\alpha})$. 
This gives us the recursive step of the procedure.

Let us notice that applying the recursive construction in the proof of Lemma
\ref{pui-apprx} we can start from the $i$-th approximation
\eqqref(apprx-i). As we have notice above, constructing $\phi_{i+1}(z)$ from
$\phi_i(z)$ we add only higher order terms, and hence
$\phi_{i+1}(z)-\phi_i(z)=O(z^{p_{i+1}+1-\mu})$. As the result we can conclude that we
can replace all $\phi_i(z)$ by the final function $\phi(z)=\phi_{l-1}(z)$ without
destroying the approximations \eqqref(apprx-i). 

This finishes the proof of Theorem E.
\qed

\state Remark. The sequence of the maps $u_i(z)$ constructed in Lemma
\ref{pui-apprx} is essentially an analogue of Puiseux series. Indeed, in the
case of integrable $J$ there is no need to apply Theorem \ref{cusp-pert} in
order to obtain a deformation with desired properties: one could simply add an
appropriate monomial $vz^k$ to the perturbed map. As the result, each
successive approximation $u_i(z^{d_i})$ will be a polynomial consisting of
certain initial part of the Puiseux series of the holomorphic map $u(z)$.

\begin{prop}
\label{extr-approx}
Under the hypotheses of Theorem E for any extremal approximation $\ti u(z)$ of
multiplicity $d$ and degree $p$ one has $d=d_i$ and $p=p_{i+1}$ for the
uniquely defined $i=0,\ldots l-1$, and then $\ti u(\varphi(z))-u_i(z^{d_i})=wz^{p_{i+1}}+
O(|z|^{p_{i+1} + \alpha})$ %
for an appropriate $w\in\cc^{n+1}$ and an appropriate holomorphic function
$\varphi(z)$.
\end{prop}

\proof Let $d>1$ be a divisor of $\mu$. Let $\nu>\mu$ be the number given by the
part {\bf(b)} of Comparison Theorem. Then  $\nu$ is not a
multiple of $d$. Further, Lemmas \ref{mult-apprx} and \ref{pui-apprx} ensure
that this number $\nu$ is the best possible approximation degree for the
multiple approximations of multiplicity $d$. In particular, for any element
$d_i>1$ in the sequence of divisors $(d_0=\mu,d_1,\ldots,d_l=1)$ there exists an
\emph{extremal} approximation of multiplicity $d_i$ and degree $p_{i+1}$.

Now assume that $d>1$ is a divisor of $\mu$ such that there exists an extremal
approximations of multiplicity $d$ and degree $p$. Find the smallest $d_i$
from the sequence of divisors $(d_0=\mu,d_1,\ldots,d_l=1)$ which is a multiple of
$d$. Such $d_i$ existsts because $d$ and all $d_i$ are divisors of $\mu$ and
$d_0=\mu$. The case $d_i=d$ was considered above, so we assume the contrary. Then
$d_i=d{\cdot}l$ with some integer $l>1$. Then $u^*(z):=u_i(z^l)$ is an
approximation of multiplicity $d$ having some degree $p$. By our
extremality assumption $p\geq p_{i+1}$. The equality case $p=p_{i+1}$ is
impossible since then $d<d_i$ would be not extremal. Consequently,
$p>p_{i+1}$. 

Now let $\eta_i=e^{2\pi\isl/d_i}$ be the primitive root of unity of degree
$d_i$. Then $\eta_i^l$ is the primitive root of unity of degree $d$. Besides
$u(\phi(\eta_iz))-u(\phi(z))=v_i(\eta_i^{p_{i+1}}-1)z^{p_{i+1}} + O(|z|^{p_{i+1}+\alpha})$ by
Theorem E. Put $w_i:= v_i(\eta_i^{p_{i+1}}-1)$ for simplicity. Let us consider
$u(\phi(\eta z))-u(\phi(z))$. Since $\eta=\eta_i^l$, we obtain 
\[\textstyle
u(\phi(\eta z))-u(\phi(z))= \sum_{j=0}^{l-1}u(\phi(\eta_i^{j+1}z))-u(\phi(\eta_i^jz))
=\Big(\sum_{j=0}^{l-1}\eta_i^{jp_{i+1}} \Big) w_iz^{p_{i+1}} + O(|z|^{p_{i+1}+\alpha}).
\]
Observe that $w_i$ is orthogonal to $v_0$ and $p_{i+1}$ is not a multiple of
$d$, since otherwise $d_{i+1}=\gcd(d_i,p_{i+1})$ would be also a multiple of
$d$ in contradiction to the choice of $d_i$. Now the second assertion of Lemma
\ref{mult-apprx} implies that $\sum_{j=0}^{l-1}\eta_i^{jp_{i+1}}$ vanishes. This
can occur only in the case when $l{\cdot}p_{i+1}$ is a multiple of $d_i$. But then
$p_{i+1}$ must be a multiple of $d=\frac{d_i}{l}$, and again
$d_{i+1}=\gcd(d_i,p_{i+1})$ would be a multiple of $d$. The obtained
contradiction shows that we must have $d=d_i$.

\qed

\newprg[prgSING.TYPE]{Singularity type of pseudoholomorphic curves}
Finally, we define the notion of singularity type of pseudoholomorphic curves
and show that every such  singularity type can be realized by an appropriate
$J$-holomorphic curve.

\begin{defi}
\label{sing-type}
A {\slsf singularity type} of pseudoholomorphic curves is a finite sequence of
integers $(p_0,\ldots,p_l)$ with the following properties: 
$1<p_0<p_1<\cdots<p_l$, the sequence $d_i:=\gcd(p_0,\ldots,p_{i})$ is strictly
decreasing, $d_{i-1}>d_i$, and  $d_l=\gcd(p_0,\ldots,p_l)=1$. The numbers
$d_i=\gcd(p_0,\ldots,p_{i})$ are called {\slsf associated divisors} of the
singularity type. 
\end{defi}

\state Remarks.~1.
More precisely, the notion in our definition is the \emph{topological} 
singularity type. For a finer notion \emph{analytic} 
singularity type of analytic or algebraic curves (especially for plane ones)
see \eg \cite{GLSh} and \cite{Ka}.

\state 2. In the higher-dimensional case $n\geq3$ there is some additional part
of the topological structure of a cuspidal curve $C=u(\Delta)$ not covered by the
characteristic exponents $p_0<p_1<\cdots$.  For
example, the condition ``{\sl $v_2$ and $v_1$ are linearly
dependent}'' is left behind. Since our primary interest lies in almost complex
surfaces we leave this topic to the interested reader.

\begin{prop}
 \label{type-exist} Let $J$ be a Lipschitz-continuous almost complex structure
 in the unit ball $B$ in $\cc^n$ and $(p_0,\ldots,p_l)$ a singularity type of
 pseudoholomorphic curves.  Then for any sequence of vectors $v_0,\ldots,v_l\in\cc^n$
 there exists a sequence of $J$-holomorphic maps $u_i:\Delta_r\to B$ defined in the
 disc $\Delta_r$ of some radius $r>0$ such that $u_0(z)=v_0z+O(|z|^{1+\alpha})$ and
 $u_i(z) = u_{i-1}(z^{d_{i-1} /d_i})+ v_i\cdot z^{p_i /d_i} + O(|z|^{p_i /d_i
  +\alpha})$ for $i=1,\ldots,l$. %
 In particular, if $v_1,\ldots,v_l$ are orthogonal to $v_0$, then $u_l(z)$ has
 singularity type $(p_0,\ldots,p_l)$.
\end{prop}

\proof The existence of $u_i(z)$ with the desired properties follows from
Theorem \ref{cusp-pert}.

\qed

\newprg[prgSING.examp]{An Example} Let us consider the following example to
the Theorem E.
\begin{exmp} \rm
Consider a (usual) holomorphic map $u:\Delta\to\cc^2$ given by
\[
u(z)=(z^{12}+z^{30}, z^{24}+z^{30}+z^{36}+z^{42}+z^{46}+z^{47}).
\]
Then the $v_0=(1,0)$ is the tangent vector at $z=0$ and $\mu=p_0=12$ is the
multiplicity. Further, its characteristic exponents are
$\mu=p_0=12,p_1=30,p_2=44,p_3=47$, and the corresponding divisors are
$d_0=p_0=12$, $d_1=6$, $d_2=2$, and $d_0=3$. On the other hand, the map $u(z)$
is a finite series polynomial which includes also the exponents $q=24,q=36$
and $q=44$, however they are non-essential (non-characteristic). A Puiseux
approximation sequence for $u(z)$ is:\\[4pt]
$\bullet\ $ $u_0(z)=(z,z^2) =v_0{\cdot}z+O(z^2)$ with $u(z)-u_0(z^{12})=O(z^{36})$ and
$v_0=(1,0)$,
\\[4pt]
$\bullet\ $ $u_1(z)=(z^2+z^5,z^4+z^5+z^6+z^7)=u_0(z^{d_0/d_1}) +v_1{\cdot}z^{p_1/d_1} +
O(z^{p_1/d_1+1})$ with $v_1=(1,1)$ and $u(z)-u_1(z^{d_1})=O(z^{46})$,
\\[4pt]
$\bullet\ $ $u_2(z)=(z^6,z^{12}+z^{15}+z^{18}+z^{21}+z^{23})=u_1(z^{d_1/d_2})
+v_2{\cdot}z^{p_2/d_2} + O(z^{p_2/d_2+2})$ with $v_2=(0,1)$ and
$u(z)-u_1(z^{d_2})=O(z^{47})$,
\\[4pt]
$\bullet\ $ $u_3(z)=u(z)=u_2(z^{d_2/d_3}) +v_3{\cdot}z^{p_3/d_3}$ with $v_3=v_2=(0,1)$.
\end{exmp}

\newprg[prgSING.index]{Equisingular deformations and cusp index formula.} %

\smallskip%
In this subsection we prove the formula expressing the cusp index of a
planar pseudoholomorphic curve via characteristic exponents at the singular
points. Let $B$ be the ball in $\cc^2$, $J$ a Lipschitz almost complex
structure in $B$ with $J(0)=J\st$, $u:\Delta\to B$ a $J$-holomorphic map with
$u(z)=v_0z^\mu+O(|z|^{\mu+\alpha})$ such that $\mu\geq2$ and $v_0\neq0\in\cc^2$, and
$(p_0=\mu,p_1,\ldots,p_l)$ the topological type of $u$ at $0$.

\begin{lem}\label{deform-type}Let $J_s$ be a family of Lipschitz-continuous
 almost complex structures in $B$ depending continuously on $s\in[0,1]$ such
 that $J_0=J$ and $J_s(0)=J\st$. Then there exists a family of
 $J_s$-holomorphic maps $u_s:\Delta_r\to B$ defined in some smaller disc of radius
 $r>0$ depending continuously on $s\in[0,1]$ such that $u_0(z)=u(z)$,
 $u_s(z)=v_0z^\mu+O(|z|^{\mu+\alpha})$, and such that $(p_0=\mu,p_1,\ldots,p_l)$ is the common
 singularity type for each $u_s(z)$ at $z=0$.
\end{lem}

\proof Let $d_i:=\gcd(p_0,\ldots,p_i)$ be the sequence of associated divisors. In
Theorem E we have constructed a sequence $u_i(z)$ of multiple
approximations of $u(z)$ such that $u_l(z)=u(\phi(z))$ and %
$u(\phi(z)) - u_i(z^{d_i}) = v_{i+1}z^{p_{i+1}} + O(|z|^{p_{i+1}+\alpha})$ for
$i=0,\ldots,l-1$ with an appropriate holomorphic reparameterization $\phi(z)$ and
$v_i\in\cc^2$. We are going to include these maps in a sequence of families of
$J_s$-polymorphic maps $u_{i,s}(z)$ satisfying similar relations with the same
numerical and vector-valued parameters $p_i,d_i\in\nn$, $v_i\in\cc^2$. Let us
formally set $u_{-1,s}(z)\equiv0$, this is the constant family of constant maps
$u_{-1,s}:\Delta\to B$. Then each family $u_{i,s}(z)$, $i=0,\ldots,l$ can be considered as
a solution of the equation
\[
\dbar_{J_s}\big(u_{i-1,s}(z^{d_{i-1}/d_i})+w_{i,s}(z)z^{p_i/d_i}\big) =0
\]
on the family of unknown functions $w_{i,s}(z)\in L^{1,p}(\Delta,\cc^2)$ satisfying
$w_{i,s}(0)=v_{i,s}$. As is the proof of Theorem \ref{cusp-pert}, we want to
obtain the needed functions as the limit of the Newton's successive
approximation procedure of the form \eqqref(newton).

To insure the convergence of the Newton's procedure we need to make our
initial data sufficiently small. For this purpose we make the rescaling
(dilatation) as in the proof of Lemma 6.1. Thus we may assume that
$\norm{J_s-J\st}_{\calc^{Lip}(B)}\leq\varepsilon$ and $\norm{u_i(z)}_{\calc^{1,\alpha}(\Delta)}\leq\varepsilon$
with some $\varepsilon\ll1$.

Now let us fix some $i$ in the interval $1,2,\ldots,l$. Set
$\nu_i:=\frac{p_i}{d_i}$. Define the structures, operators, etc $J_{i,s}^{(\nu)},
R_{i,s}^{(\nu)},T^0_{J_{i,s}^{(\nu)},R_{i,s}^{(\nu)}}, F^{(\nu)}_{i,s}(z,w)$ by the
same formulas as in the proof of Theorem \ref{cusp-pert} substituting $J_s$
instead of $J$, $u_{i-1,s}(z^{d_{i-1}/d_i})$ instead of $u_0(z)$, $\nu_i$ instead
of $\nu$, and so on. Use index $n$ for numeration of successive approximations
$w_{i,s,n}(z)$ in the procedure. In this way we obtain the formula
\begin{equation}\eqqno(shreck1)
w_{i,s,n+1} = T_{J_{i,s}^{(\nu)},R_{i,s}^{(\nu)}}^0\big[F^{(\nu)}_{i,s}(z,w_{i,s,n}(z))\big] + w_{i,s,1}(z),
\end{equation}

The only difference from the procedure \eqqref(newton), which is the key idea
of the proof of the present lemma, lies in the choice of the initial data
$w_{i,s,1}(z)$ of the approximation. Recall that by Lemma \ref{lem3.2},
$u_i(z)=u_{i-1}(z^{d_{i-1}/d_i})+w_i(z)z^{p_i/d_i}$ with some function
$w_i(z)\in L^{1,p}(\Delta,\cc^2)$ with $w_i(0)=v_i$. We use this function instead of
the constant function $\equiv w_0$ in \eqqref(syst2). This means that now
$w_{i,s,1}(z)$ is defined by
\begin{equation}\eqqno(shreck2)
w_{i,s,1}(z) := w_i(z) - T_{
 J_{i,s}^{(\nu)},R_{i,s}^{(\nu)}}^0\big(D_{J^{(\nu)}_{i,s}, u_{i-1,s}(z^{d_{i-1}/d_i})}  w_i(z) \big).
\end{equation}

Since our initial data \eqqref(shreck2) were made small enough, the Newton's
approximation procedure \eqqref(shreck1) converges for every $s\in[0,1]$.
Moreover, for $s=0$ the iteration \eqqref(shreck1) is constant,
$w_{i,0,n+1}(z)=w_{i,0,n}(z)=\ldots=w_{i,0,1}(z)=w_i(z)$ since such was our choice
of the initial data $w_{i,s,1}(z)$.

Now, substitute the  limit functions $w_{i,s,\infty}(z)$ in the relations
$u_{i,s}(z)=u_{i-1,s}(z^{d_{i-1}/d_i})+w_{i,s,\infty}(z)z^{p_i/d_i}$ successively
for $i=0,1,\ldots,l$, and set $u_s(z):=u_{l,s}(\phi\inv(z))$. The obtained family
$u_s(z)$ fulfills the requirements of the lemma.
\qed

\medskip

\begin{prop}\label{index-form} Let $(X,J)$  be an almost
 complex surface with  Lipschitz-continuous structure $J$ and $u:\Delta\to X$ a
 $J$-holomorphic with a singularity at $z=0$ of the type $(p_0,\ldots,p_l)$. Then
 the cusp index of $u(\Delta)$ at $u(0)$ is given by the formula
\begin{equation}
\eqqno(form-kappa) \varkappa = \frac12\sum_{i=1}^m (d_{i-1}- d_i)(p_i -1),
\end{equation}
\end{prop}

\proof Without loss of generality we may assume that $X$ is the unit ball in
$\cc^2$, $u(0)=0\in B$, and $J(0)=J\st$. Define a family of Lipschitz structures
$J_s$ in $B$, $s\in[0,1]$, by the formula $J_s(w_1,w_2)=J(sw_1,sw_2)$. Then
$J_s$ depends continuously on $s\in[0,1]$, $J_1=J$, and $J_0=J\st$. By Lemma
\ref{deform-type}, there exists a family $u_s(z)$ of $J_s$-holomorphic maps
defined in some small disc $\Delta_r$ and depending continuously on $s\in[0,1]$, such
that each $u_s$ has the singularity type $(p_0,\ldots,p_l)$ at $z=0$. Fix some
sufficiently small $\rho>0$ and denote by $\gamma_s$ the intersection $u_s(\Delta_r)$ with
the sphere $S^3_\rho$ of radius $\rho$. Then $\gamma_s$ is an isotopy of knots in $S^3_\rho$
each transverse to the induced contact structure $F_s:=TS^3_\rho\cap J_s(TS^3_\rho)$ on
$S^3_\rho$. Consequently, each $\gamma_s$ has the same Bennequin index $b$ related to
the cusp index of each $u_s(\Delta_r)$ by the formula $\varkappa=(b+1)/2$.  In particular,
each $u_s(\Delta_r)$ has the same cusp-index at $0$.

Since $J_0=J\st$, it is sufficient to consider the case of integrable
structures. In this case the formula is well-known, see \cite{Wl}, page 85 and
Exercise 6.7.2.

\qed

\smallskip  %
\state Remark. It is proved in \cite{Wl}, Section 5, that the
Alexander polynomial of the link $\gamma_\rho= S^3_\rho\cap u(\Delta)$ of the singularity
determines the whole set of characteristic exponents of the singularity. In
particular, $u(\Delta)$ is non-singular at $0$ if and only if the corresponding
link is unknot. Notice that the Alexander polynomial of a knot in an invariant
of the \emph{smooth} isotopy class; in contrary, the Bennequin index is an
invariant of \emph{transversal} isotopy class of a knot.

\newsect[sectEXMP]{Examples and Open Questions}

\begin{exmp}\rm
\label{exmp9.1}
There exists an almost complex structure $J$ on a domain  $X\subset \rr^4$ which belong
to $\bigcap_{1<p<\infty }L^{1,p}(X, \End_{\rr}TX)$
and two $J$-complex curves $M_i$ which coincide by an non-empty but proper open part.
\end{exmp}
The first curve will be the coordinate plane $M_1\deff\rr^\subset \rr^4$ with coordinates
$x_1,y_1$. The second --- $M_2$ --- is defined by equations
\begin{equation}
\eqqno(second)
y_2=0, \quad
x_2 =
\begin{cases}
e^{-\frac{1}{x_1^k}} \text{ if } x_1\geq 0,\cr
0 \text{ if } x_1\leq 0.
\end{cases}
\end{equation}
Since $x_2^{'}(x_1)=(e^{-\frac{1}{x_1^k}})'= - \frac{k}{x_1^{k+1}}e^{-\frac{1}{x_1^k}}$, we see that
the vector $(1,0,- \frac{k}{x_1^{k+1}}e^{-\frac{1}{x_1^k}},0)=$  $(1, 0,- x_2(-\ln x_2)^\frac{k+1}{k},0)$
is tangent
to $M_2$ at every point $(x_1, y_1, e^{-\frac{1}{x_1^k}},0)\in M_2$. Extend it to a vector field
\begin{equation}
v(x_1,y_1,x_2,y_2) =
\begin{cases}
(1, - x_2(-\ln x_2)^\frac{k+1}{k}) \text{ if } x_1\geq 0,\cr
(1,0,0,0) \text{ if } x_1\leq 0,
\end{cases}
\end{equation}
in $X\deff \rr^2\times (-\infty , 1)\times\rr$. The structure $J$ is now defined by
\begin{equation}
\begin{cases}
J\frac{\d}{\d x_2} = \frac{\d }{\d y_2},\cr
Jv = \frac{\d }{\d y_1}.
\end{cases}
\end{equation}
Both $M_1$ and $M_2$ are clearly $J$-convex. The regularity of $J$ is that of $v$, \ie
is $\mathsf{Ln}^{1+\frac{1}{k}}\mathsf{Lip}$. Since,  obviously $\calc^{Ln^{1+\frac{1}{k}}Lip}\subset
\bigcap_{1<p<\infty }L^{1,p}$, we are done.

\begin{exmp}\rm
\label{exmp9.2}
We shall construct an example of a Lipschitz-continuous almost
complex structure $J$ in $\rr^4$ and a $J$-holomorphic map
$u:\Delta\to\rr^4$ which is exactly from $\calc^{1,LnLip}(\Delta)$.
\end{exmp}

Consider the following function $u(z)=z^2\ln (|z|^2)$. Set $v(z) =
\dbar u(z) = \frac{z^2}{\bar z}\in \calc^{Lip}(\Delta)$. Remark that
$\d u(z) = 4z\ln |z| + z\in \calc^{LnLip}(\Delta)$. Therefore $u\in
\calc^{1,LnLip}(\Delta)$. Let us interpret the vector function
$(u,z)$ as a $J$-holomorphic curve for certain Lipschitz $J$. Namely
let us  take

\begin{equation}
J = \left(
\begin{array}{rrrr}
0&-1 &v_2&-v_1\\
1 &0&-v_1&-v_2\\
0&0&0&-1 \\
0&0&1 &0\\
\end{array}\right),
\end{equation}
where $v_1+iv_2=v$ constructed above. One readily checks that $J$ is
an almost complex structure and $(u,1)$ is $J$-holomorphic. $J$ has
the same regularity as $v$, \ie is Lipschitz-continuous.

\bigskip We would like to finish with an open question close to the topics considered in this
paper. For an arbitrary (continuous) $\rr$-linear endomorphism $A = A(z)$ of
the trivial $\cc^n$-bundle over $\Delta$, define the operator
$\dbar_A$ on $L^{1,1} _\loc$-sections of $\cc^n$ by the usual
formula
\begin{equation}
\eqqno(9.5)
\dbar_A u := (\partial_x + A\cdot \partial_y) u.
\end{equation}

\begin{rema} \rm 
One can rewrite this example using operator $Q$ as in \eqqref(q-bar-j). The corresponding 
$Q$ has the form
\begin{equation}
Q(u_1,u_2) = 
\begin{pmatrix}
0 & \frac{u_2^2}{\bar u_2}\cr
0 & 0\cr 
\end{pmatrix}.
\end{equation}

\end{rema}

\begin{quest} Let $A$ be a continuous endomorphism of the trivial
$\cc^n$-bundle over $\Delta$ such that $|A(z) - J\st| \leq c \cdot
|z|^\beta$ with some $c<1$ and $0<\beta<1$. Let $u\in L^{1,1}
_\loc(\Delta, \cc^n)$ be not identically $0$ and satisfy in the weak
sense the inequality
\begin{equation}
| \dbar_A u | \leq h\cdot | u |.
\eqqno(9.6)
\end{equation}
for some nonnegative $h\in L^p_\loc(\Delta)$ with $2< p <
\frac{2}{1-\beta}$. Prove that there exists $\mu\in\nn$ such that
$u(z)=z^\mu \cdot g(z)$ for some $g\in L^{1, p}_\loc (\Delta)$ with
$g(0) \neq 0$.
\end{quest}

This time let $A$ be a Lipschitz-continuous $Mat(2n,\rr)$-valued function on
the unit disc  $\Delta$ and let $\bar\d_A$ be defined by \eqqref(9.5).
We suppose that $\bar\d_A$ is uniformly elliptic, \ie its spectrum $s(A)$ is separated from
$\rr$ in $\cc$. Let  $u$ be a solution of a differential inequality

\begin{equation}
\norm{\d_Au} \leq C\norm{u}.
\eqqno(9.7)
\end{equation}

\begin{quest}
Suppose that for some $z_n\to 0$ one has $u(z_n)=0$. Does it implies that
$u\equiv 0$?
\end{quest}

\medskip If $n=1$, \ie for $\cc$ valued function this is so and it follows
from Theorem 35 of \cite{B} via the trick explained on the page 101.

\medskip And the last question, which closely related to the first and second ones.

\begin{quest}
Let $J$ be an almost complex structure in $\rr^{2n}$ of class $\calc^{\alpha}$ for some $0<\alpha <1$
and let $u:\Delta\to \rr^{2n}$ be $J$-holomorphic. Suppose that for some sequence $z_n\to 0$ one has
$u(z_n)=0$. Does it implies that $u\equiv 0$?
\end{quest}

\begin{rema} \rm 
Very recently a considerable progress in the direction of these questions was made by J.-P. Rosay, \cite{Ro}.
\end{rema}

\ifx\undefined\bysame
\newcommand{\bysame}{\leavevmode\hbox to3em{\hrulefill}\,}
\fi

\def\entry#1#2#3#4\par{\bibitem[#1]{#1}
{\textsc{#2 }}{\sl{#3} }#4\par\vskip2pt}

\bigskip

\bigskip


\end{document}